\documentclass[12pt]{article}
\usepackage{amssymb}
\usepackage{latexsym}
\usepackage[all]{xy}
\newenvironment{conjalt}[1]{\vspace{1ex}
\par\noindent{\bf Conjecture} #1 \em}{\vspace{1ex}\par}

\newtheorem{cor}{Corollary}
\newtheorem{defn}{Definition}
\newtheorem{hyp}{Hypothesis}
\newtheorem{lemma}{Lemma}
\newtheorem{prop}{Proposition}
\newtheorem{propdefn}{Proposition/Definition}
\newtheorem{thm}{Theorem}

\newcounter{example}
\newcounter{rem}[section]

\newcommand{\barbbQ}{\ensuremath{\bar{\bbQ}}}

\newcommand{\bbC}{\ensuremath{\mathbb{C}}}

\newcommand{\bbF}{\ensuremath{\mathbb{F}}}

\newcommand{\bbN}{\ensuremath{\mathbb{N}}}

\newcommand{\bbQ}{\ensuremath{\mathbb{Q}}}

\newcommand{\bbR}{\ensuremath{\mathbb{R}}}
\newcommand{\bbZ}{\ensuremath{\mathbb{Z}}}
\newcommand{\bbZp}{\ensuremath{\mathbb{Z}_p}}

\newcommand{\beq}{\begin{equation}}
\newcommand{\beql}[1]{\begin{equation}\label{#1}}
\newcommand{\bPf}{\noindent \textsc{Proof\ }}

\newcommand{\cA}{\ensuremath{\mathcal{A}}}
\newcommand{\cB}{\ensuremath{\mathcal{B}}}

\newcommand{\cH}{\ensuremath{\mathcal{H}}}
\newcommand{\cI}{\ensuremath{\mathcal{I}}}

\newcommand{\Cl}{{\rm Cl}}

\newcommand{\cN}{\ensuremath{\mathcal{N}}}
\newcommand{\cO}{\ensuremath{\mathcal{O}}}

\newcommand{\cR}{\ensuremath{\mathcal{R}}}
\newcommand{\cS}{\ensuremath{\mathcal{S}}}

\newcommand{\cupdot}{\stackrel{.}{\cup}}

\newcommand{\cX}{\ensuremath{\mathcal{X}}}

\newcommand{\displaymapdef}[5]
{\[
\begin{array}{rcrcl}
 #1 &:& #2 &\longrightarrow& #3 \\
    & &    &                    \\
    & & #4 &\longmapsto    & #5
\end{array}
\]}
\newcommand{\dotlist}[3]{\ensuremath{#1_{#2},\ldots,#1_{#3}}}
\newcommand{\eeq}{\end{equation}}
\newcommand{\eg}{\emph{e.g.}}
\newcommand{\eitheror}[4] 
{\left\{                  
\begin{array}{ll}         
#1& \mbox{#2}\\
#3& \mbox{#4}
\end{array}
\right.}
\newcommand{\ePf}{\hspace*{\fill}~$\Box$\vertsp\par}
\newcommand{\eps}{\ensuremath{\varepsilon}}
\newcommand{\etc}{\emph{etc.}}
\newcommand{\fa}{\ensuremath{\mathfrak{a}}}

\newcommand{\fD}{\ensuremath{\mathfrak{D}}}
\newcommand{\ff}{\ensuremath{\mathfrak{f}}}

\newcommand{\fm}{\ensuremath{\mathfrak{m}}}
\newcommand{\fM}{\ensuremath{\mathfrak{M}}}

\newcommand{\fp}{\ensuremath{\mathfrak{p}}}
\newcommand{\fP}{\ensuremath{\mathfrak{P}}}
\newcommand{\fq}{\ensuremath{\mathfrak{q}}}
\newcommand{\fQ}{\ensuremath{\mathfrak{Q}}}

\newcommand{\fs}{\ensuremath{\mathfrak{s}}}
\newcommand{\fS}{\ensuremath{\mathfrak{S}}}

\newcommand{\Gal}{{\rm Gal}}

\newcommand{\half}{\frac{1}{2}}
\newcommand{\Hom}{{\rm Hom}}
\newcommand{\ie}{\emph{i.e.}}
\newcommand{\im}{\ensuremath{{\rm im}}}
\newcommand{\inv}{^{-1}}
\newcommand{\listsub}[2]{\ensuremath{{#1_1,\ldots,#1_#2}}}
\newcommand{\ndiv}{\nmid}
\newcommand{\nin}{\not\in}
\newcommand{\ord}{{\rm ord}}

\newcommand{\rem}{\refstepcounter{rem}\noindent{\sc Remark \therem}}

\renewcommand{\therem}{\thesection.\arabic{rem}}

\newcommand{\vertsp}{\vspace{1ex}}

\newcommand{\wedgesub}[2]{\ensuremath{{#1_1\wedge\ldots\wedge#1_#2}}}
\newcommand{\ZmodZ}[1]{\ensuremath{{\bbZ/#1\bbZ}}}
\newcommand{\aKkSm}{\ensuremath{a^-_{K/k,S}}}
\newcommand{\badS}{\ensuremath{{\rm Bad}(S)}}
\newcommand{\barLdU}{\ensuremath{\overline{\bigwedge^d_{\bbZ \bar{G}} U_{S}(K^+)}}}

\newcommand{\oridPhi}[1]{\ensuremath{[#1]_{\fa_\Phi}}}

\newcommand{\PhiKk}{\ensuremath{\Phi_{K/k}}}

\newcommand{\piKKp}{\ensuremath{\pi_{K/K^+}}}
\newcommand{\pnpo}{\ensuremath{{p^{n+1}}}}
\newcommand{\Qf}{\ensuremath{\Qxi{f}}}
\newcommand{\Qxi}[1]{\ensuremath{\bbQ(\xi_{#1})}}

\newcommand{\SKkS}{\ensuremath{\fS_{K/k,S}}}
\newcommand{\sKkS}{\ensuremath{\fs_{K/k,S}}}
\newcommand{\Sram}{\ensuremath{S_{\rm ram}}}

\newcommand{\tDelta}{\ensuremath{\tilde{\Delta}}}
\newcommand{\te}{\ensuremath{\tilde{e}}}

\newcommand{\tf}{\ensuremath{\tilde{f}}}

\renewcommand{\th}{\ensuremath{\tilde{h}}}
\newcommand{\ThetaKkS}{\ensuremath{\Theta_{K/k,S}}}
\newcommand{\ThetaKk}{\ensuremath{\Theta_{K/k}}}
\newcommand{\ThetaKpkS}{\ensuremath{\Theta_{K^+/k,S}}}

\newcommand{\WUoKpm}{\ensuremath{{\textstyle \bigwedge^d_{\bbZ_pG}U^1(K_p)^-}}}
\newcommand{\Zbp}{\ensuremath{\bbZ_{(p)}}}
\newcommand{\ZpnpoZ}{\ensuremath{\bbZ/\pnpo \bbZ}}
\newcommand{\ZpnpoZG}{\ensuremath{(\bbZ/\pnpo \bbZ) G}}
\newcommand{\ZpnpoZst}{\ensuremath{(\ZpnpoZ)^\times}}
\begin{document}
\title{Abelian $L$-Functions at $s=1$ and\\
Explicit Reciprocity for Rubin-Stark Elements}
\author{D. Solomon
\\King's College, London\\
david.solomon@kcl.ac.uk}
\maketitle
\begin{center}
\bf \large Abstract\vertsp\\
\end{center}
Given an abelian, CM extension $K$ of any totally real number field $k$, we restate and generalise two conjectures `of Stark type' made in~\cite{twizo}.
The Integrality Conjecture concerns the image of a $p$-adic map $\sKkS$ determined by the minus-part of the $S$-truncated equivariant $L$-function for $K/k$ at $s=1$.
It is connected to the Equivariant Tamagawa Number Conjecture of Burns and Flach.
The Congruence Conjecture says that $\sKkS$ gives an explicit reciprocity law for the element predicted by the corresponding
Rubin-Stark Conjecture for $K^+/k$. We then study the general properties of these conjectures and prove one or both of them under various hypotheses, notably
when $p\ndiv [K:k]$, when $k=\bbQ$ or when $K$ is absolutely abelian.
\section{Introduction}
Throughout this paper $k$ will be a number field of finite degree $d$ over $\bbQ$ and $K$ will be a finite, Galois extension of $k$ such
that the group $G:=\Gal(K/k)$ is abelian. We denote by $S_\infty=S_\infty(k)$ and $S_{\rm ram}=S_{\rm ram}(K/k)$ the sets consisting respectively of
the infinite places of $k$ and those which are finite and ramify in $K$, and we set $S^0=S^0(K/k)=S_{\rm ram}\cup S_\infty$. If $S$ is any finite set of places
containing $S^0$ and $s$ a complex number with ${\rm Re}(s)>1$, we define a convergent Euler product in the complex group ring of $G$ (denoted simply $\bbC G$) by
\beq\label{eq: introduction}
\Theta_{K/k, S}(s):=\prod_{\fq\nin S}\left(1-N\fq^{-s}\sigma_\fq\inv\right)\inv\
\eeq
The product ranges over those places $\fq$ of $k$  which are not in $S$. (Here and henceforth, finite places are identified with prime ideals.)
$\sigma_\fq=\sigma_{\fq,K}$ denotes the Frobenius element of $G$ for $\fq$. If $k$ is totally real and $K$ a $CM$ field with complex
conjugation $c\in G$,
it can be shown that the `minus part' $\Theta_{K/k, S}^-(s):=\half(1-c)\Theta_{K/k, S}(s)$ extends to an entire function $\bbC\rightarrow \bbC G$. This paper concerns
two conjectures of a $p$-adic nature about the element $a^-_{K/k, S}:=(i/\pi)^d\Theta_{K/k, S}^-(1)$ (whose coefficients turn out to be algebraic).
For any  number $p$ we denote by $U^1(K_p)$ the $p$-semilocal principal units of $K$ and define
a $p$-adic regulator on the exterior power $\bigwedge^d_{\bbZ_p G}U^1(K_p)$. By combining this with $a^-_{K/k, S}$ we obtain a map
$\sKkS:\bigwedge^d_{\bbZ_p G}U^1(K_p)\rightarrow\bbQ_p G$.
Assuming for the rest of this Introduction that $p\neq 2$  and $S$ contains all the places above $p$ in $k$, our first conjecture
(the `Integrality Conjecture' or `IC') states simply that the image $\SKkS$ of
$\sKkS$ is contained in $\bbZ_p G$. Recall now that if $K^+$ denotes the maximal totally real subfield of $K$  then (the meromorphic continuation of)
$\Theta_{K^+/k, S}(s)$ has a zero of order at least $d$ at $s=0$. Furthermore, a well known conjecture of Stark
(reformulated and refined by Rubin) states that the coefficient of $s^d$ in the Taylor expansion of $\Theta_{K^+/k, S}(s)$
is given by evaluating an ($\bbR G$-valued) regulator map on
an element of a certain  exterior power of the (global) $S$-units of $K^+$.
Imposing further natural conditions makes this element the unique `Rubin-Stark Element' of the title, here denoted $\eta_{K^+/k,S}$.
Our second conjecture (the `Congruence Conjecture' or `CC') assumes and refines the IC, and in so doing
links the minus part of $\Theta_{K/k, S}(s)$ at $s=1$ to its
plus part at $s=0$. It says that
if $K$ contains the $\pnpo$-th roots of unity for some $n\geq 0$,
then, very roughly speaking, the reduction of $\sKkS$ modulo $\pnpo$ gives the explicit reciprocity law for $\eta_{K^+/k,S}$.

The idea for these conjectures came from the results of~\cite{sconzp}. We shall
not elaborate on the precise connection in the present paper beyond saying that if $p$ splits in $k$ then certain rather strong hypotheses
considered in~\cite{sconzp} imply a weak form of the CC at each
level in a cyclotomic $\bbZ_p$-tower containing $K$. The IC and the CC first appeared explicitly as Conjectures~5.2 and~5.4 at the end of~\cite{twizo} in a form
less general and  more awkward than the present versions.
Conjectures~5.2 and~5.4 also used twisted zeta-functions at $s=0$ where the IC and the CC use
$\Theta_{K/k, S}^-(1)$ and so will, perhaps, prove more accessible.

The remainder of this paper is organised as follows. Section~\ref{sec: dram pers} contains
the precise definitions and basic properties of the main players: the elements $a^-_{K/k, S}$ and $\eta_{K^+/k,S}$, the map $\sKkS$ and
the pairing $H_{K/k,n}$ (a determinant of additive, equivariant Hilbert symbols in terms of which our conjectural reciprocity law is couched).
Section~\ref{section:conjectures} contains the precise statements of the two conjectures.
Section~\ref{sec: evidence} surveys the current evidence in their favour -- now quite considerable -- and
includes the statements of the three main results of this paper which were announced in~\cite{twizo}: Firstly, in the case $p\ndiv |G|$,
we give a complete characterisation of the $\SKkS$ in terms of $L$-functions of odd characters of $G$ at $s=0$.
In this case the IC then follows, thanks to a result of Deligne-Ribet and Pi.~Cassou-Nogu\`es. Secondly, we prove the conjectures in the case $k=\bbQ$, using
an explicit reciprocity law due to Coleman. Thirdly we prove the conjectures when $K/\bbQ$ is abelian (but $k$ is
not necessarily $\bbQ$) by `base-change' from the previous result. In this case,
we require a relatively mild technical
hypothesis on $K/k$, $S$ and $p$. We also discuss briefly A.~Jones' recent results on a rather different refinement
of the IC which follows from a special case of the Equivariant Tamagawa Number Conjecture (ETNC) of Burns and Flach.
(On the other hand, we should mention that the CC currently has no known connection with the ETNC.)
Section~\ref{sec: functoriality} examines the behaviour of the conjectures as
$S$, $K$ and $n$ vary. Sections~\ref{sec: pndivG},~\ref{sec: k equals Q} and~\ref{sec: Kabsab} contain the proofs of the three main results referred to above.

Jones' refinement of the IC mentioned above states that $\SKkS$ is contained in the Fitting ideal
(as $\bbZ_p G$-module) of the minus part of the $p$-part of a certain ray-class group of $K$.
This creates the possibility of links between $\SKkS$ and recent work of Greither on Fitting ideals of
{\em duals} of class groups (see~\cite{Greither}).
Another perspective comes from recent work of the author showing that, in the case $k=\bbQ$ at least,
the CC creates a link between the map $\sKkS$ and some new maps and ideals in Iwasawa Theory.
The latter has connections with the Main Conjecture and applications  to the
{\em plus} part of the class group of abelian fields.

In addition to those introduced above, we use the following basic notations and conventions.
If $\cR$ is a  commutative ring and $H$ a finite abelian group, we write $\cR H$ for the group-ring and if $M$ is a $\bbZ G$-module we shall sometimes
abbreviate $\cR\otimes_\bbZ M$ to $\cR M$ (considered as a $\cR H$-module in the obvious way).
For any subgroup $D\subset H$, we write $N_D$
for the norm element $\sum_{d\in D}d\in\cR H$. If $m$ is a positive integer, we denote by
$\mu_m(\cR)$ the group of  all $m$th roots of unity in $\cR$ and for any prime number $p$
we set $\mu_{p^\infty}(\cR)=\bigcup_{i=0}^\infty\mu_{p^i}(\cR)$.
All number fields in this paper are supposed of finite degree over $\bbQ$ and are considered as subfields
of $\barbbQ$ which is \emph{the algebraic closure of $\bbQ$ within $\bbC$}. We shall write
$\xi_m$ for the particular generator $\exp(2\pi i/m)$ of  $\mu_m(\barbbQ)$.
For any number field $F$ and any integer $r$ we shall write $S_r(F)$ for the set of places (prime ideals) of $F$ dividing $r$.
If $S$ is a set of places of $F$ and $L$ is any finite extension of $F$ we shall write $S(L)$ for the set of places of $L$ lying above those in $S$.
If $S$ contains $S_\infty(F)$ (see above) then the
group $U_S(F)$ of $S$-units of $F$ consists of those elements of $F^\times$ which are local units at every place not in $S$
and we shall often write simply $U_S(L)$ in place of $U_{S(L)}(L)$.
({\em Caution}: $U_S$ and related modules will sometimes be written additively). If $H$ is abelian and
$v$ is any place of $F$ we shall write $D_v(L/F)$ for the
decomposition subgroup of $H$ at any prime dividing $v$ in $L$ and similarly $T_v(L/F)$ for the inertia subgroup
(if $v$ is finite). Suppose $L\supset F\supset M$ are three number fields such that $L/M$ and $F/M$ are Galois extensions. Then
the restriction  map $\Gal(L/M)\rightarrow\Gal(F/M)$ will be denoted $\pi_{L/F}$ and extended $\cR$-linearly to
a ring homomorphism  $\cR\Gal(L/M)\rightarrow\cR\Gal(F/M)$ for
any commutative ring $\cR$. We also write $\nu_{L/M}$ for the $\cR$-linear `corestriction' map
$\cR\Gal(F/M)\rightarrow\cR\Gal(L/M)$ which sends $g\in \Gal(F/M)$ to the sum
of its preimages under $\pi_{L/F}$ in $\Gal(L/M)$.

It is my pleasure to thank Andrew Jones, Cristian Popescu and Jonathan Sands for useful discussions about
the contents of this paper.
I wish also to thank Cristian and UCSD for their hospitality during the sabbatical year in
which the majority of this paper was written.

\section{Dramatis Person\ae}\label{sec: dram pers}
 \subsection{The Function $\ThetaKkS$ and the Element $\aKkSm$}
Let $\hat{G}$ denote the dual group of $G$, namely the group
of all (irreducible) complex characters $\chi:G\rightarrow\bbC^\times$ with identity element $\chi_0$, the trivial character.
For any $\chi\in \hat{G}$ we write $e_{\chi,G}$
for the associated idempotent in the complex group-ring $\bbC G$, namely $e_{\chi,G}:=\frac{1}{|G|}\sum_{g\in G}\chi(g)g\inv$.
Expanding the Euler product~(\ref{eq: introduction}), we get
\beq\label{eq:thetaStozetaandL}
\Theta_{K/k, S}(s)=\sum_{g\in G}\zeta_{K/k, S}(s;g)g\inv
=\sum_{\chi\in \hat{G}}L_{K/k, S}(s,\chi)e_{\chi\inv,G}
\eeq
for ${\rm Re}(s)>1$. Here, $\zeta_{K/k, S}(s;g)$ and $L_{K/k, S}(s,\chi)$ denote respectively
the `$S$-truncations' of the partial zeta-function attached to $G$
and the $L$ function attached to $\chi$. In particular
\beq\label{eq:Lfns}
L_{K/k, S}(s,\chi)=\prod_{\fq\nin S}\left(1-N\fq^{-s}\chi(\sigma_\fq)\right)\inv=
\prod_{\fq\in S\setminus S_\infty \atop \fq\ndiv \ff_\chi}\left(1-N\fq^{-s}\hat{\chi}([\fq])\right)L(s,\hat{\chi})
\eeq
where $\ff_\chi$ and $L(s,\hat{\chi})$ denote respectively the conductor of $\chi$ and the $L$-function of its associated primitive ray-class
character $\hat{\chi}$ modulo $\ff_\chi$.\vertsp\\
\rem\label{rem: artin L functions}
The second (but not the first) expression for $L_{K/k, S}(s,\chi)$ in~(\ref{eq:Lfns}) makes sense when $S$
is any finite set of places of $k$, containing $S_\infty(k)$ but \emph{not necessarily} $\Sram(K/k)$.
In fact it agrees with the definition of the $S$-truncated \emph{Artin $L$-function} attached to $\chi$ considered as a character of $G$ (see for example~\cite[p.~23]{Tate}).\vertsp\\
\noindent
The analytic behaviour of  $L(s,\hat{\chi})$ is well-known. Its (in general)
meromorphic continuation means that we may use equations~(\ref{eq:thetaStozetaandL}) and~(\ref{eq:Lfns}) to continue $\Theta_{K/k,S}$ to a
meromorphic, $\bbC G$-valued function
on $\bbC$. These equations then hold as identities between meromorphic functions on $\bbC$.
Similarly, if  $S\supset S'\supset S^0$, then the obvious identity
\beq\label{eq:thetaS'toS}
\Theta_{K/k,S}(s)=\prod_{\fq\in S\setminus S' }\left(1-N\fq^{-s}\sigma_\fq\inv\right)\Theta_{K/k,S'}(s)\
\eeq
for ${\rm Re}(s)>1$ also holds for all $s$.
In fact, the function $L(s,\hat{\chi})$, hence also the function
$\chi(\Theta_{K/k,S}(s))=L_{K/k, S}(s,\chi\inv)$, is analytic on $\bbC\setminus\{1\}$ and
\beq\label{eq:order_of_chi_theta_at_s_equals_one}
\ord_{s=1}\chi(\Theta_{K/k,S}(s))=
\left\{
\begin{array}{rl}
0& \mbox{if $\chi\neq \chi_0$}\\
-1& \mbox{if $\chi= \chi_0$}
\end{array}
\right.
\eeq
Moreover, the residue of  $\chi_0(\Theta_{K/k,S}(s))=
\prod_{\fq\in S\setminus S_\infty} \left(1-N\fq^{-s}\right)\zeta_k(s)$ at $s=1$ is well-known (see \eg~\cite[Th\'eor\`eme I.1.1]{Tate}).

Using the well known functional equation relating
the \emph{primitive} $L$-function $L(s,\hat{\chi}\inv)$ to $L(1-s,\hat{\chi})$ one might expect to
derive a natural relation between $\ThetaKkS(s)$ and $\ThetaKkS(1-s)$ by means of~(\ref{eq:thetaStozetaandL}) and~(\ref{eq:Lfns}).
There are however two obstacles to this:
firstly, the dependence on $\ff_\chi$ of the second product in~(\ref{eq:Lfns})
and secondly, the presence of (Galois) Gauss sums in the functional equations. Instead, in~\cite{twizo} we used these functional equations   to
give a precise relation
between $\Theta_{K/k,S^0}(s)$ (there denoted $\ThetaKk(s)$) and $\PhiKk(1-s)$, where the function $\PhiKk:\bbC\rightarrow \bbC G$ was defined by means of
twisted zeta-functions and studied, together with its  $p$-adic analogues, in~\cite{twizas,zetap1, sconzp, twizo}.)
For each $v\in S_\infty$, we write $c_v$ for the unique generator of $D_v(K/k)$ so that $c_v=1$ unless $v$ is real and one (hence every)
place $w$ of $K$ above $v$ is complex, in which case
$c_v$ is the complex conjugation associated to any such $w$. We define an entire, $\bbC D_v(K/k)$-valued function
\[
C_v(s)=\left\{
\begin{array}{ll}
\exp(i\pi s).1-\exp(-i\pi s).c_v=2i\sin(\pi s).1&
\mbox{if $v$ is complex}\\
\exp(i\pi s/2).1+\exp(-i\pi s/2).c_v&
\mbox{if $v$ is real}
\end{array}
\right.
\]
Then  Theorem~2.1 of~\cite{twizo}, combined with~(\ref{eq:thetaS'toS}) for $S'=S_0$, gives
\beq
\label{eq:phioneminusstothetas}
i^{r_2(k)}\sqrt{|d_k|}\prod_{\fq\in S\setminus S^0} \left(1-N\fq^{-s}\sigma_\fq\inv\right)\PhiKk(1-s)
=((2\pi)^{-s}\Gamma(s))^d\left(\prod_{v\in S_\infty}C_v(s)\right)\ThetaKkS(s)
\eeq
where $r_2(k)$ denotes the number of complex places of $k$ and $d_k$ its absolute discriminant.
Let $\Theta^{\rm n.t.}_{K/k, S}(s)$ be the function $(1-e_{\chi_0,G})\ThetaKkS(s)$, which is regular at $s=1$
by~(\ref{eq:order_of_chi_theta_at_s_equals_one}). So~(\ref{eq:phioneminusstothetas}) gives
\beq\label{eq:phi_zero_to_theta_one}
\sqrt{|d_k|}\prod_{\fq\in  S\setminus S^0} \left(1-N\fq^{-1}\sigma_\fq\inv\right)(1-e_{\chi_0,G})\PhiKk(0)
=(2\pi)^{-d}i^{|S_\infty|}\left(\prod_{v\in S_\infty}(1-c_v)\right) \Theta^{\rm n.t.}_{K/k, S}(1)
\eeq
from which it follows that $(1-e_{\chi_0,G})\PhiKk(0)$ vanishes unless $k$ is totally real and $K$ is totally complex. On the other hand,
multiplying~(\ref{eq:phioneminusstothetas}) by $e_{\chi_0,G}$ and letting $s\rightarrow 1$, we see that $e_{\chi_0,G}\PhiKk(0)$
vanishes unless $|S_\infty|=1$,
\ie\ $k$ is $\bbQ$ or an imaginary quadratic field, in which case it may easily be calculated from ${\rm res}_{s=1}\zeta_k(s)$.
Thus $\PhiKk(0)$ has little
interest unless $k$ is totally real and $K$ is totally complex. Even then, $\prod_{v\in S_\infty}(1-c_v)$
vanishes unless there is a (unique) CM-subfield $K^-$
of $K$ containing $k$, in which case we lose little but complication by replacing $K$ by $K^-$.
(See Remark~3.1(i) of~\cite{twizo} for further explanations).
For these reasons \emph{we shall henceforth make the}
\begin{hyp}\label{hyp:K_CM_k_totreal}
$k$ is totally real and $K$ is a CM field.
\end{hyp}
This means that $d_k$ is a positive integer and
$c_v=c$, the unique complex conjugation in $G$,
for all $v\in S_\infty$. Let $e^{\pm}$  denote the two idempotents $\half(1\pm c)$ of $\bbC G$ and let $\Theta^-_{K/k,S}(s)$ be the entire function
$e^-\ThetaKkS(s)=e^-\Theta^{\rm n.t.}_{K/k, S}(s)$. The above remarks, together with
a simple calculation of $e_{\chi_0,G}\PhiKk(0)$ when $k=\bbQ$, show that equation~(\ref{eq:phi_zero_to_theta_one}) may be rewritten as
\beq\label{eq:defn_aKkS}
a^-_{K/k,S}:=\left(\frac{i}{\pi}\right)^d \Theta^-_{K/k,S}(1)=
\left\{
\begin{array}{ll}
{\displaystyle \prod_{\fq\in S\setminus S^0}} \left(1-N\fq^{-1}\sigma_\fq\inv\right)\sqrt{d_k}\PhiKk(0)& \mbox{if $k\neq \bbQ$}\\
{\displaystyle \prod_{q\in S\setminus S^0}} \left(1-q^{-1}\sigma_q\inv\right)\Phi_{K/\bbQ}(0)+
{\displaystyle \half\prod_{q\in S\setminus\{\infty\}}}\left(1-q^{-1}\right)e_{\chi_0,G}& \mbox{if $k=\bbQ$}
\end{array}
\right.
\eeq
If $R$ is a commutative ring in which $2$ is invertible and $M$ any $R\langle c\rangle$-module then we shall write
$M^+$ (\emph{resp.} $M^-$) for the $R$-submodule $e^+ M$ (\emph{resp.} $e^-M$), so that $M=M^+\oplus M^-$.
In this notation,
$a^-_{K/k,S}$ clearly lies in $\bbC G^-$ and multiplying~(\ref{eq:defn_aKkS}) by $e^-$ gives
\beq\label{eq:aKkS_intermsof_PhiKko}
a^-_{K/k,S}=e^-a^-_{K/k,S}=e^-\prod_{\fq\in S\setminus S^0} \left(1-N\fq^{-1}\sigma_\fq\inv\right)\sqrt{d_k}\PhiKk(0)
\eeq
whether or not  $k=\bbQ$, but if $k\neq\bbQ$ then the term $e^-$ may be omitted on the R.H.S.
In fact, $a^-_{K/k,S}$ has algebraic coefficients: Let $\ff(K)$ be the integral ideal of $\cO_k$ which is the conductor
of $K/k$ in the sense of class-field theory and let
$f(K)$ be the positive generator of the ideal $\ff(K)\cap\bbZ$. The
product in~(\ref{eq:aKkS_intermsof_PhiKko}) lies in $\bbQ G^\times$, then~(\ref{eq:aKkS_intermsof_PhiKko})
and~\cite[ Prop.~3.1]{twizo} show that
$a^-_{K/k,S}$ has coefficients in $\sqrt{d_k}\bbQ(\mu_{f(K)})$ and that
\beq\label{eq:generating QmufGminus}
a^-_{K/k,S}\bbQ(\mu_{f(K)})G=\sqrt{d_k}\bbQ(\mu_{f(K)})G^-
\eeq
Integrality properties of the coefficients of $a^-_{K/k,S}$ are given in~\cite{Solomon-Roblot2} where it
is shown that they also lie in the Galois closure of $K$ over $\bbQ$. (See {\em ibid.}, Proposition~2 and  Remark~6).

\subsection{Rubin-Stark Elements for $K^+/k$}\label{sec: Rubin-Stark}
Let us write $\bar{G}$ for $\Gal(K^+/k)\cong G/\langle c\rangle$ so that $\piKKp:\bbC G\rightarrow\bbC \bar{G}$ induces an ring isomorphism
$\bbC G^+\rightarrow \bbC \bar{G}$  sending $e^+\Theta_{K/k,S}(s)$ onto $\Theta_{K^+/k,S}(s)$. To study $\Theta_{K^+/k,S}(s)$ at $s=0$.
we define an integer $r_S(\phi)$
for each $\phi\in\widehat{\bar{G}}$  by
\beq\label{eq:defn_rschi}
r_S(\phi):=
\left\{
\begin{array}{ll}
d+|\{\fq\in S\setminus S_\infty\,:\,\phi(D_\fq(K^+/k))=\{1\}\}| &
\mbox{if $\phi$ is non-trivial}\\
d+|S\setminus S_\infty|-1=|S|-1 &
\mbox{if $\phi$ is trivial}
\end{array}
\right.
\eeq
Since $k$ and $K^+$ are totally real, the functional equation of $L(s,\widehat{\phi})$ for $\phi\in\widehat{\bar{G}}$ shows that, for any such $\phi$
we have
\beq\label{eq:ordLKpluskSatzero}
\ord_{s=0}\phi(\ThetaKpkS(s))=\ord_{s=0}L_{K^+/k,S}(s,\phi)=r_S(\phi)
\eeq
(see \eg~\cite[Ch.~I, \S 3]{Tate}). We shall assume until further notice
\begin{hyp}\label{hyp:S_contains_d+1_places}
$|S|\geq d+1$ (\ie\ $S$ contains at least one finite place.)
\end{hyp}
This implies that $r_S(\phi)\geq d$ for $\phi$ trivial hence for every $\phi\in\widehat{\bar{G}}$, so we may define
\[
\Theta_{K^+/k, S}^{(d)}(0):=\lim_{s\rightarrow 0}s^{-d}\Theta_{K^+/k, S}(s)
\]
(an element of $\bbC \bar{G}$ which is easily seen to lie in $\bbR \bar{G}$). Conjectures of Stark, as refined by Rubin~\cite{Rubin},
predict that $\Theta_{K^+/k, S}^{(d)}(0)$ is given by a certain $\bbR \bar{G}$-valued regulator of $S$-units of
$K^+$ defined as follows.
We fix once and for all a set $\tau_1,\ldots,\tau_d$ of left coset representatives for $\Gal(\barbbQ/k)$ in $\Gal(\barbbQ/\bbQ)$
and we define $\bbQ \bar{G}$-linear, real logarithmic
maps for $i=1,\ldots,d$:
\displaymapdef{\lambda_{K^+/k,i}}{\bbQ U_S(K^+)}{\bbR \bar{G}}
{\lambda_{K^+/k,i}(a\otimes \eps)}{a\sum_{g\in \bar{G}}\log|\tau_i(g\eps)|g\inv\in \bbR \bar{G}}
The above-mentioned regulator is the $\bbQ \bar{G}$-linear map
uniquely defined by:
\displaymapdef{R_{K^+/k}}{\textstyle \bigwedge^d_{\bbQ \bar{G}}\bbQ U_{S}(K^+)}{\bbR \bar{G}}{x_1\wedge\ldots\wedge x_d}{\det(\lambda_{K^+/k,i}(x_t))_{i,t=1}^d}
The following definition generalises the above
construction and will be useful later.
\begin{propdefn}\label{propdefn:determinantal maps}
\begin{enumerate}
\item Suppose $\cS$ is a commutative ring, $\cR$ a commutative $\cS$-algebra and that $M$ is any (left) $\cR H$ module for a finite group $H$.
There is an isomorphism from  ${\rm Hom}_{\cR}(M,\cS)$ to ${\rm Hom}_{\cR H}(M,\cS H)$ given by $f\mapsto f^H$ where $f^H$ is defined to be the map
$m\mapsto\sum_{h\in H}f(h\inv m)h$.

\item Suppose $H$ is abelian and $l\in\bbN$. Then  for every $l$-tuple $(f_1,\ldots,f_l)\in {\rm Hom}_{\cR}(M,\cS)^l$ there is a
$\cR H$-linear determinantal map $\Delta_{f_1,\ldots,f_l}$ uniquely defined by
\displaymapdef{\Delta_{f_1,\ldots,f_l}}{\textstyle \bigwedge^l_{\cR H}M}{\cS H}{m_1\wedge\ldots\wedge m_l}{\det(f^H_i(m_t))_{i,t=1}^l}
$\Delta_{f_1,\ldots,f_l}$ is $\cS$-multilinear and alternating as a function of $(\listsub{f}{d})$. Moreover for each $i=1,\ldots,l$ and  $h\in H$ we have
$\Delta_{f_1,\ldots, f_i\circ h,\ldots ,f_l}(\mu)=\Delta_{f_1,\ldots,f_l}(\mu)h$ for all $\mu\in \textstyle \bigwedge^l_{\cR H}M$.\ePf
\end{enumerate}
\end{propdefn}
For instance, taking $\cR=\bbQ$, $\cS=\bbR$,  $M=\bbQ U_S(K^+)$ and $H=\bar{G}$ gives $R_{K^+/k}=\Delta_{f_1,\ldots,f_d}$ where $f_i$ is the map
sending $a\otimes \eps\in \bbQ U_S(K^+)$ to its logarithmic embedding
$a\log|\tau_i(\eps)|$ in $\bbR$. If instead we take $\cR=\cS=\bbQ$, then any $d$ elements
$f_1,\ldots,f_d$ of $\Hom_{\bbQ }(\bbQ U_S(K^+),\bbQ)$ give rise to a
$\bbQ \bar{G}$-linear map $\Delta_{f_1,\ldots,f_d}: \bigwedge^d_{\bbQ \bar{G}}\bbQ U_{S}(K^+)\rightarrow \bbQ G$.
Let us identify  $\Hom_{\bbZ}(U_S(K^+),\bbZ)$ with the lattice in $\Hom_{\bbQ}(\bbQ U_S(K^+),\bbQ)$ which is its image under the map
$f\rightarrow 1\otimes f$. We can then
define a $\bbZ \bar{G}$ submodule $\Lambda_{0,S}=\Lambda_{0,S}(K^+/k)$ of $\bigwedge^d_{\bbQ \bar{G}}\bbQ U_{S}(K^+)$ by
\[
\Lambda_{0,S}(K^+/k):=\left\{
\eta\in\textstyle
\bigwedge^d_{\bbQ \bar{G}}\bbQ U_{S}(K^+)\,:\,
\Delta_{f_1,\ldots,f_d}(\eta)\in\bbZ \bar{G}\ \forall\,
f_1,\ldots,f_d\in \Hom_{\bbZ}(U_S(K^+),\bbZ)
\right\}
\]
This coincides with `$\Lambda_0^d U_S(K^+)$' as defined by Rubin's `double dual' construction in~\cite[\S 1]{Rubin}.
It is clear that $\Lambda_{0,S}$ contains the lattice which is the natural image
of $\bigwedge^d_{\bbZ \bar{G}} U_{S}(K^+)$ in
$\textstyle \bigwedge^d_{\bbQ \bar{G}}\bbQ U_{S}(K^+)$ (we denote this $\overline{\bigwedge^d_{\bbZ \bar{G}} U_{S}(K^+)}$)
but the two are not necessarily equal. In fact, Prop.~1.2 of~\cite{Rubin} implies
\begin{prop}\label{prop:index of barLdU in Lambda}
If $d=1$ (\ie\ $k=\bbQ$) then $\Lambda_{0,S}=\overline{\bigwedge^1_{\bbZ \bar{G}} U_{S}(K^+)}=\overline{U_{S}(K^+)}$. In general,
The index $|\Lambda_{0,S}\,:\,\barLdU |$ is finite and supported on primes dividing $|\bar{G}|$.
\ePf
\end{prop}
Let us define an idempotent $e_{S,d,\bar{G}}$, \emph{a priori} in $\bbC \bar{G}$, by setting
$
e_{S,d,\bar{G}}:=\sum_{\phi\in\widehat{\bar{G}}\atop r_S(\phi)=d}e_{\phi,\bar{G}}
$
. This is  the unique element element $x$ of $\bbC \bar{G}$ such that $\phi(x)=1$ or $0$ according as $r_S(\phi)=d$ or $r_S(\phi)>d$.
It follows easily from this description and the formula~(\ref{eq:defn_rschi}) that
\beq\label{eq: formula for eSdGbar}
e_{S,d,\bar{G}}=
\left\{
\begin{array}{ll}
{\displaystyle \prod_{\fq\in S\setminus S_\infty}\left(1-\frac{1}{|D_\fq(K^+/k)|}N_{D_\fq(K^+/k)}\right)}&
\mbox{if $|S|>d+1$}\\
{\displaystyle \left(1-\frac{1}{|D_\fq(K^+/k)|}N_{D_\fq(K^+/k)}\right)+e_{\chi_0,\bar{G}}}&
\mbox{if $|S|=d+1$, \ie\ $S=\{\fq\}\cup S_\infty$}
\end{array}
\right.
\eeq
Thus  $e_{S,d,\bar{G}}$ is an idempotent of $\bbQ \bar{G}$, so lies in $|G|\inv\bbZ \bar{G}$. We also deduce easily:
\begin{prop}\label{prop: equivalent `eigenspace' conditions} Let $M$ be any $\bbQ \bar{G}$-module and $m\in M$. The following are equivalent
\begin{enumerate}
\item $m\in e_{S,d,\bar{G}}M$
\item $m=e_{S,d,\bar{G}}m$
\item \label{part3: equivalent `eigenspace' conditions}
For all $\fq\in S\setminus S_\infty$,
\beq\label{eq: espace cond for q}
N_{D_\fq(K^+/k)}m\in
\left\{
\begin{array}{ll}
\{0\}&
\mbox{if $|S|>d+1$}\\
M^{\bar{G}}&
\mbox{if $|S|=d+1$, \ie\ $S=\{\fq\}\cup S_\infty$}
\end{array}
\right.
\eeq
\item \label{part4: equivalent `eigenspace' conditions}
$e_{\phi,\bar{G}}(1\otimes m)=0$ in $\bbC\otimes_\bbQ M$ for all $\phi\in\widehat{\bar{G}}$ such that $r_S(\phi)>d$\ePf
\end{enumerate}
\end{prop}
For brevity, we shall sometimes refer to any of these conditions as the \emph{eigenspace condition} on $m$ w.r.t.\ $(S,d,\bar{G})$.
Now, given any subring $\cR$ of $\bbQ$, we formulate a version of the Rubin-Stark conjecture `over $\cR$':
\begin{conjalt}{$RSC(K^+/k,S;\cR)$}\label{conj:weak rubin_stark}\\
Let $K/k$, $S$ and be as above, satisfying Hypotheses~\ref{hyp:K_CM_k_totreal} and~\ref{hyp:S_contains_d+1_places}. Then there exists
an element
$\eta\in\textstyle \bigwedge^d_{\bbQ \bar{G}}\bbQ U_{S}(K^+)$ satisfying the eigenspace condition w.r.t.\ $(S,d,\bar{G})$ and such that
\beq\label{eq:rubin_stark}
\Theta_{K^+/k, S}^{(d)}(0)=R_{K^+/k}(\eta)
\eeq
and
\beq\label{eq: eta lies in halfRLambda}
\eta\in{\textstyle \half}\cR\Lambda_{0,S}(K^+/k)
\eeq
\end{conjalt}
\noindent Notice that $\Theta_{K^+/k, S}^{(d)}(0)$ lies in the
ideal $e_{S,d,\bar{G}} \bbR \bar{G}$  (and in fact generates it) by
equation~(\ref{eq:ordLKpluskSatzero}). Thus if $\eta\in\textstyle
\bigwedge^d_{\bbQ \bar{G}}\bbQ U_{S}(K^+)$ is any solution
of~(\ref{eq:rubin_stark}) then $e_{S,d,\bar{G}}\eta$ is a solution
satsifying the eigenspace condition. On the other hand, it can be
shown that $R_{K^+/k}$ is injective  on $e_{S,d,\bar{G}}\textstyle
\bigwedge^d_{\bbQ \bar{G}}\bbQ U_{S}(K^+)$ (this follows from~\cite[Lemma 2.7]{Rubin})
so a solution of~(\ref{eq:rubin_stark}) satisfying the eigenspace condition is unique. For this reason, we call such an element {\em `the Rubin-Stark element
for $K^+/k$ and $S$'} and denote it $\eta_{K^+/k,S}$ since it is independent of $\cR$.
Of course, Condition~(\ref{eq: eta lies in halfRLambda}) is redundant if $\cR=\bbQ$ and for any prime number $p$ we have
\[
RSC(K^+/k,S;\bbZ)\Rightarrow RSC(K^+/k,S;\Zbp)\Rightarrow RSC(K^+/k,S;\bbQ)
\]
(where $\Zbp$ denotes the localisation $\{a/b\in\bbQ:p\ndiv b\}$).
Moreover $RSC(K^+/k,S;\bbZ)$ is equivalent to the conjunction of
$RSC(K^+/k,S;\Zbp)$ for all primes $p$. We shall mainly be interested in $RSC(K^+/k,S;\Zbp)$ when $p\neq 2$, in which case~(\ref{eq: eta lies in halfRLambda}) reduces to
$\eta\in \Zbp\Lambda_{0,S}$.\vertsp\\
\rem\label{rem: changing taus}
Since $R_{K^+/k}$ depends on the choice
(and ordering) of the $\tau_i$'s, so will $\eta_{K^+/k,S}$, but in a simple way. For example,
if one $\tau_i$ is replaced by $\tau_i\tau\inv$ for some $\tau\in\Gal(\barbbQ/k)$ then we must replace
$\eta_{K^+/k,S}$ by $\tau|_{K^+}\eta_{K^+/k,S}$ where $\tau|_{K^+}\in \bar{G}$.\vertsp\\
\rem\label{rubin_stark_remark 2}
$RSC(K^+/k,S;\bbQ)$ and $RSC(K^+/k,S;\bbZ)$ follow from certain special cases of  Conjectures~$A'$ and~$B'$ of~\cite{Rubin} respectively.
Indeed, if we choose the extension `$K/k$' of Rubin's paper to be our $K^+/k$,
his `$S$' to be ours, his `$r$' to be $d$ and his chosen places
`$w_1,\ldots,w_r$' to be the real places of $K^+$ defined by $\tau_1,\ldots,\tau_d$. Then Rubin's
Hypotheses~2.1.1-2.1.4 are satisfied. His conjectures also require an auxiliary set
$T$ of finite places of $k$ satisfying certain conditions, although for Conjecture~$A'$ the precise choice of such $T$ does not affect the truth
of the conjecture.
For simplicity we take $T=\{\fq\}$ for some prime
$\fq\nin S$ not dividing $2$ and splitting in  $K^+$ (infinitely many of these exist by \v{C}ebotarev's theorem).
Then Rubin's Hypothesis~2.1.5 certainly holds since $U_S(K^+)_{tor}=\{\pm 1\}$. Moreover his
`$\Theta^{(r)}_{S,T}(0)$' is our $(1-N\fq)\Theta_{K^+/k, S}^{(d)}(0)$ and
his `$\Lambda_0^r U_{S,T}$'
is a sublattice of our $\Lambda_{0,S}(K^+/k)$ which also
spans $\textstyle \bigwedge^d_{\bbQ \bar{G}}\bbQ U_{S}(K^+)$ over $\bbQ$. It follows easily that $RSC(K^+/k,S;\bbQ)$ is equivalent to
Rubin's Conjecture~$A'$ with these choices and this (hence any) $T$ . Moreover, if both hold then
Rubin's `$\eps_{S,T}$' equals our $(1-N\fq)\eta_{K^+/k, S}$, by uniqueness. It follows that Rubin's Conjecture~$B'$
with these choices amounts to the further condition that $(1-N\fq)\eta_{K^+/k, S}$
lie in  his $\Lambda^r_0 U_{S,T}$ hence in our $\Lambda_{0,S}(K^+/k)$.
But as $\fq$ varies subject to the above conditions, Lemme~IV.1.1 of~\cite{Tate} says that the
g.c.d.\ of the corresponding integers $1-N\fq$ is $|\mu(K^+)|=2$.
Thus the corresponding cases of Rubin's Conjecture~$B'$ together imply $RSC(K^+/k,S;\bbZ)$.

The connection with Stark's original conjecture in terms of characters
(see~\cite[Conjecture I.5.1]{Tate}) is as follows. Propositions 2.3 and 2.4 of~\cite{Rubin} show that it
holds for $K^+/k$, $S$ and every character $\phi\in \bar{G}$
satisfying  $r_S(\phi)=d$ if and only if Rubin's Conjecture $A'$ holds (for any $T$) which is equivalent to $RSC(K^+/k,S;\bbQ)$, by the above.\vertsp\\
\noindent In the next section we shall be interested in determinantal maps obtained
from a $d$-tuple $(f_1,\ldots,f_d)\in\Hom_\bbZ(U_S(K^+),\bbZ/\pnpo\bbZ)$ for some prime $p$ and $n\geq 0$.
Taking $\cR=\bbZ$, $\cS=\bbZ/\pnpo\bbZ$ and $M=U_S(K^+)$ in Proposition/Definition~\ref{propdefn:determinantal maps} gives
such a map $\Delta_\listsub{f}{d}:{\textstyle \bigwedge^d_{\bbZ \bar{G}}U_{S}(K^+)}\rightarrow (\bbZ/\pnpo\bbZ)\bar{G}$. We shall now show that
{\em provided $p$ is odd,}
this map `extends' naturally to $\Zbp\Lambda_{0,S}$ in a sense to be explained below. First, we have
\begin{lemma}\label{lemma:exatseqforHom}
If $p$ is odd then the following sequence is exact.
\[
0\rightarrow \Hom_\bbZ(U_S(K^+),\bbZ){\stackrel{\pnpo}{\longrightarrow}} \Hom_\bbZ(U_S(K^+),\bbZ)\longrightarrow \Hom_\bbZ(U_S(K^+),\bbZ/\pnpo\bbZ)\rightarrow 0
\]
\end{lemma}
\bPf\ $K^+$ is totally real so $U_S(K^+)/\{\pm 1\}$ is $\bbZ$-free. Thus the sequence is exact if $U_S(K^+)$ is replaced by $U_S(K^+)/\{\pm 1\}$.
But since $\bbZ$ and  $\bbZ/\pnpo\bbZ$ have no $2$-torsion, we may identify $\Hom_\bbZ(U_S(K^+)/\{\pm 1\},\bbZ)$ with $\Hom_\bbZ(U_S(K^+),\bbZ)$ and
$\Hom_\bbZ(U_S(K^+)/\{\pm 1\},\bbZ/\pnpo\bbZ)$ with $\Hom_\bbZ(U_S(K^+),\bbZ/\pnpo\bbZ)$.\ePf
\noindent Thus given $f_1,\ldots,f_d$ in $\Hom_\bbZ(U_S(K^+),\bbZ/\pnpo\bbZ)$ for $p$ odd, we can choose lifts $\listsub{\tf}{d}$ in $\Hom_\bbZ(U_S(K^+),\bbZ)$.
As previously, we may regard these as elements of $\Hom_\bbQ(\bbQ U_S(K^+),\bbQ)$ and use Proposition/Definition~\ref{propdefn:determinantal maps} to
construct $\Delta_{\listsub{\tf}{d}}:{\textstyle \bigwedge^d_{\bbQ \bar{G}}\bbQ U_S(K^+)}\rightarrow \bbQ \bar{G}$. If $\eta\in\Lambda_{0,S}$ then
$\Delta_{\listsub{\tf}{d}}(\eta)$ lies in $\bbZ \bar{G}$ by definition of $\Lambda_{0,S}$ and we write $\tDelta_\listsub{f}{d}(\eta)$
for its image in $(\bbZ/\pnpo\bbZ)\bar{G}$. The latter is independent of the choice of each lift $\tf_i$, as one easily checks using Lemma~\ref{lemma:exatseqforHom},
the linearity of $\Delta_{\listsub{\tf}{d}}$ in $\tf_i$ and the fact that $\eta\in\Lambda_{0,S}$. Consequently we have a well-defined map
$\tDelta_\listsub{f}{d}:\Lambda_{0,S}\rightarrow(\bbZ/\pnpo\bbZ) \bar{G}$
which is linear and so extends uniquely to $\bbZ_{(p)}\Lambda_{0,S}$.
It is now an easy exercise to check the following properties of $\tDelta_\listsub{f}{d}$:
\begin{prop}\label{prop:properties of tilde R} Let $p$ be odd and choose $f_1,\ldots,f_d\in\Hom_\bbZ(U_S(K^+),\bbZ/\pnpo\bbZ)$
\begin{enumerate}
\item\label{part1:properties of tilde R}  The map
$\tDelta_\listsub{f}{d}:\bbZ_{(p)}\Lambda_{0,S}(K^+/k)\rightarrow(\bbZ/\pnpo\bbZ)\bar{G}$ is $\bbZ \bar{G}$-linear.
\item\label{part2:properties of tilde R} It is also
$(\bbZ/\pnpo\bbZ)$-multilinear and alternating as a function of $(\listsub{f}{d})$ and for each $i=1,\ldots,d$ we have
$\tDelta_{f_1,\ldots, f_i\circ g,\ldots ,f_d}(\eta)=\tDelta_{f_1,\ldots,f_d}(\eta)g$ for all $g\in \bar{G}$ and $\eta\in \bbZ_{(p)}\Lambda_{0,S}(K^+/k)$.
\item\label{part3:properties of tilde R} The following diagram commutes
\[
\xymatrix{
{\textstyle \bigwedge^d_{\bbZ \bar{G}}U_{S}(K^+)}\ar[dd]_{\alpha}\ar[drrrr]^{\Delta_\listsub{f}{d}}&&&&\\
&&&&(\bbZ/\pnpo\bbZ)\bar{G}\\
\bbZ_{(p)}\Lambda_{0,S}(K^+/k)\ar[urrrr]_{\tDelta_\listsub{f}{d}}&&&&
}
\]
where $\alpha$ is the natural map
${\textstyle \bigwedge^d_{\bbZ \bar{G}}U_{S}(K^+)}\rightarrow{\textstyle \bigwedge^d_{\bbQ \bar{G}}\bbQ U_S(K^+)}$ with restricted range.
\end{enumerate}\ePf
\end{prop}
\rem\ This shows in particular that $\Delta_\listsub{f}{d}$ vanishes on the kernel of $\alpha$ in the above diagram. One can show that $\ker(\alpha)$
is always finite and supported on primes dividing $2|\bar{G}|=|G|$. Also,
Proposition~\ref{prop:index of barLdU in Lambda} implies that ${\rm im}(\alpha)$
spans $\bbZ_{(p)}\Lambda_{0,S}$ over $\Zbp$ whenever $p\ndiv |G|$. So if  $p\ndiv |G|$ then {\em any} $\bbZ \bar{G}$-linear map
$F:\bigwedge^d_{\bbZ \bar{G}} U_{S}(K^+)\rightarrow(\bbZ/\pnpo\bbZ)\bar{G}$  vanishes on $\ker(\alpha)$ and has a {\em unique} `extension'
$\tilde{F}:\Zbp\Lambda_{0,S}\rightarrow(\bbZ/\pnpo\bbZ)\bar{G}$ satisfying $F=\tilde{F}\circ\alpha$.
\subsection{Hilbert Symbols and the Pairing $H_{K/k,n}$}
Suppose that $L$ is a local field containing $\mu_m$
for some positive integer $m$ coprime to the characteristic of $L$. We recall that
the Hilbert symbol is the map
\displaymapdef{(\cdot,\cdot)_{L,m}}{L^\times\times L^\times}{\mu_m\subset L^\times}{(\alpha,\beta)}{(\beta^{1/m})^{\sigma_{\alpha,L}-1}}
where $\beta^{1/m}$ is any $m$th root of $\beta$ in any abelian closure $L^{\rm ab}$ of $L$ and $\sigma_{\alpha,L}$ denotes the image of $\alpha$ under the reciprocity
homomorphism $(\cdot,L)$ of local class field theory
from $L^\times$ to $\Gal(L^{\rm ab}/L)$. The Hilbert symbol is bilinear and skew-symmetric. For the general theory, see~\cite[Ch.~12]{Artin-Tate},
\cite[V.3]{NeukirchANT} or \cite[Ch.~XIV]{Serre Local fields}. (Note that our notation $(\alpha,\beta)_{L,m}$
is compatible with that of~\cite{Artin-Tate} and
\cite{NeukirchANT} but represents the element denoted $(\beta,\alpha)$ in
\cite{Serre Local fields} and is similarly reversed in the notation of~\cite{Coleman Dilogarithm})

Let $p$ be a prime number and $n\geq 0$ an integer.  We shall assume until further notice that
{\em $K$ contains $\mu_{\pnpo}$} for some $n\geq 0$.
Let $\kappa_n:\Gal(\barbbQ/\bbQ)\rightarrow\ZpnpoZst$ be the cyclotomic character modulo $\pnpo$, determined by
$\tau(\zeta)=\zeta^{\kappa_n(\tau)}\ \forall\,\zeta\in\mu_\pnpo,\ \tau\in \Gal(\barbbQ/\bbQ)$. Since $K$ contains $\mu_{\pnpo}$,
the restriction of $\kappa_n$ to $\Gal(\barbbQ/k)$ factors through a homomorphism $G\rightarrow\ZpnpoZst$ which we denote by the
same symbol. We also use the shorthand $\zeta_n$ for $\xi_\pnpo\in K$.
If $K_\fP$ denotes the completion of $K$ at some prime ideal $\fP$, we may define a bilinear pairing
$[\cdot,\cdot]_{\fP,n}: K_\fP^\times \times K_\fP^\times\rightarrow\bbZ/\pnpo\bbZ$ by setting
\[
\iota_\fP(\zeta_n)^{[\alpha,\beta]_{\fP,n}}=(\alpha,\beta)_{K_\fP,\pnpo}
\ \ \ \mbox{for all $\alpha,\beta\in K_\fP^\times$}
\]
where $\iota_\fP:K\rightarrow K_\fP$ is the natural embedding.
Every $g\in G$ induces an isomorphism $K_\fP\rightarrow K_{g\fP}$, also denoted $g$ and
such that $g\circ \iota_\fP=\iota_{g\fP}\circ g$. Standard facts from Local Class Field theory imply that
$
(g\alpha,g\beta)_{K_{g\fP},\pnpo}=g(\alpha,\beta)_{K_{\fP},\pnpo}
$
in $K_{g\fP}^\times$ for any $\alpha,\beta\in K_\fP^\times$. It follows easily that
\beq\label{eq:Gaction1}
[g\alpha,g\beta]_{g\fP,n}=\kappa_n(g)[\alpha,\beta]_{\fP,n}\ \ \ \ \mbox{for all $\alpha,\beta\in K_\fP^\times,\ g\in G$}
\eeq
If $\fP$ divides $p$ then $\iota_\fP$ extends to a $\bbQ_p$-algebra map from $K_p:=K\otimes_{\bbQ} \bbQ_p$ to $K_\fP$.
Thus we obtain a
pairing
\displaymapdef
{[\cdot,\cdot]_{K,n}}{K_p^\times\times K_p^\times}{\bbZ/\pnpo\bbZ}
{(\alpha,\beta)}{\displaystyle \sum_{\fP|p}[\iota_\fP(\alpha),\iota_\fP(\beta)]_{\fP,n}}
Letting $G$ act on $K_p$ through $K$, we still have $\iota_\fP\circ g=g\circ \iota_{g\inv\fP}$  for any $g\in G$ and $\fP|p$,
so~(\ref{eq:Gaction1}) implies
\beq\label{eq:Gaction2}
[g\alpha,g\beta]_{K,n}=\kappa_n(g)[\alpha,\beta]_{K,n}\ \ \ \ \mbox{ for all $\alpha,\beta\in K_p^\times$,\ $g\in G$}
\eeq
The product map $\prod_{\fP|p}\iota_\fP: K_p\rightarrow\prod_{\fP|p} K_\fP$ is a  $G$-equivariant ring isomorphism
(where $g((x_\fP)_\fP)=(g x_{g\inv\fP})_\fP$ in $\prod_{\fP|p} K_\fP$). We shall regard this as an identification so that $\iota_\fP$ identifies with the projection
$\prod_{\fP|p} K_\fP\rightarrow K_\fP$.
Thus we identify the principal semilocal units $\prod_{\fP|p} U^1(K_\fP)$ with a  $\bbZ G$-submodule of $K_p^\times$ and denote it $U^1(K_p)$.
Regarding each $U^1(K_\fP)$ as a finitely generated $\bbZ_p$-module, $U^1(K_p)$ becomes a finitely generated $\bbZ_p G$-module.

{\em From now on we assume that  $p$ is  odd}.
Consider the unique ring automorphism of $(\ZpnpoZ) G$ sending $g\in G$ to $\kappa_n(g)g\inv$. Since $\kappa_n(c)=-1$,
this restricts to a ring isomorphism from
$(\ZpnpoZ) G^+$ to  $(\ZpnpoZ) G^-$. Composing with $\bar{2}\inv\nu_{K/K^+}:(\ZpnpoZ) \bar{G}\rightarrow(\ZpnpoZ) G$, we obtain a
ring isomorphism $\bar{\kappa}_n^\ast=\bar{\kappa}_{K,n}^\ast:(\ZpnpoZ) \bar{G}\rightarrow(\ZpnpoZ) G^-$. Explicitly, if $h\in \bar{G}$ and $g\in G$ then
\beq\label{eq:formulae for kappanbarstar}
\bar{\kappa}_n^\ast(h)=\bar{2}\inv\sum_{\th\in G\atop\piKKp(\th)=h}
\kappa_n(\th)\th\inv\ \ \mbox{and hence}\ \
\bar{\kappa}_n^\ast(\piKKp(g))=e^-\kappa_n(g)g\inv
\eeq
Given a set $S\supset S^0$  as in previous sections, any $u\in U^1(K_p)$ defines homomorphism
$f_u\in\Hom_\bbZ(U_S(K^+),\ZpnpoZ)$ by setting $f_u(\eps)=[\eps, u]_{K,n}$ (by abuse, we write $\eps$ for $\eps\otimes 1\in K_p^\times$).
Using the `$\tilde{\Delta}$' notation of the last section we may now define a map
\displaymapdef{H_{K/k,S,n}}{\Zbp\Lambda_{0,S}(K^+/k)\times U^1(K_p)^d}{\ZpnpoZG^-}{(\eta;\listsub{u}{d})}
{2^d\bar{\kappa}_n^\ast(\tDelta_{f_{u_1},\ldots f_{u_d}}(\eta))}
\begin{prop}\label{prop:pairing_H}Suppose $\eta\in \Zbp\Lambda_{0,S}(K^+/k)$ and $\listsub{u}{d}\in U^1(K_p)^d$
\begin{enumerate}
\item \label{part1:pairing_H} For any $x\in\bbZ \bar{G}$ we have
$H_{K/k,S,n}(x\eta;\listsub{u}{d})=\bar{\kappa}_n^\ast(\bar{x})H_{K/k,S,n}(\eta;\listsub{u}{d})$ where
$\bar{x}$ denotes the image of $x$ in $(\ZpnpoZ)\bar{G}$
\item\label{part2:pairing_H} $H_{K/k,S,n}$ is $\bbZ G$-multilinear (hence $\bbZ_p G$-multilinear) and alternating as a function of \listsub{u}{d}
\end{enumerate}
\end{prop}
\bPf\ Part~\ref{part1:pairing_H} follows from part~\ref{part1:properties of tilde R} of Proposition~\ref{prop:properties of tilde R}.
The $\bbZ$-multilinearity in part~\ref{part2:pairing_H}
follows from part~\ref{part2:properties of tilde R} of Proposition~\ref{prop:properties of tilde R} so it suffices to prove that
replacing $u_i$ by $gu_i$ (for $g\in G$) multiplies $H_{K/k,S,n}(\eta;\listsub{u}{d})$ by $g$, or indeed by $e^-g$ since it lies in the minus part.
But if we write $h$ for $\pi_{K/K^+}(g)\in\bar{G}$ then
Equation~(\ref{eq:Gaction2}) and Proposition~\ref{prop:properties of tilde R} part~\ref{part2:pairing_H} give
\begin{eqnarray*}
\bar{\kappa}_n^\ast(\tDelta_{f_{u_1},\ldots,f_{gu_i},\ldots, f_{u_d}}(\eta))&=&
\bar{\kappa}_n^\ast(\tDelta_{f_{u_1},\ldots,\kappa_n(g)f_{u_i}\circ h\inv,\ldots, f_{u_d}}(\eta))\\
&=&\bar{\kappa}_n^\ast(\tDelta_{f_{u_1},\ldots f_{u_d}}(\eta)\kappa_n(g){h}\inv)\\
&=&\bar{\kappa}_n^\ast(\kappa_n(g){h}\inv)\bar{\kappa}_n^\ast(\tDelta_{f_{u_1},\ldots f_{u_d}}(\eta))\\
&=&e^-g\bar{\kappa}_n^\ast(\tDelta_{f_{u_1},\ldots f_{u_d}}(\eta))\\
\end{eqnarray*}
by~(\ref{eq:formulae for kappanbarstar}) and the result follows.\ePf
\noindent By part~\ref{part2:pairing_H} of the Proposition, $H_{K/k,S,n}$ defines a unique pairing (also
denoted $H_{K/k,S,n}$) from
$\Zbp\Lambda_{0,S}\times {\textstyle \bigwedge_{\bbZ_p G}^d}U^1(K_p)$ to $(\ZpnpoZ)G^-$. By $\bbZ_p G$-linearity in the second variable, it
factors through the projection on
${\textstyle \bigwedge_{\bbZ_p G}^d}U^1(K_p)^-$ .
An important and simple special  case is when
$\eta$ equals $(1\otimes\eps_1)\wedge\ldots\wedge(1\otimes\eps_d)\in \overline{\bigwedge^d_{\bbZ \bar{G}} U_{S}(K^+)}$.
Using  Proposition~\ref{prop:properties of tilde R}~\ref{part3:properties of tilde R} and equation~(\ref{eq:formulae for kappanbarstar}) and tracing through the definitions,
we find that for all $\listsub{u}{d}$ in $U^1(K_p)$
\begin{eqnarray}
H_{K/k,S,n}((1\otimes\eps_1)\wedge\ldots\wedge(1\otimes\eps_d),\wedgesub{u}{d})&=&
  2^d\bar{\kappa}_n^\ast(\Delta_{f_{u_1},\ldots,f_{u_d}}(\wedgesub{\eps}{d}))\nonumber\\
&=&\bar{\kappa}_n^\ast\left(\det\left(2\sum_{h\in\bar{G}}[h\inv\eps_i,u_t]_{K,n}h\right)_{i,t=1}^d\right)\nonumber\\
&=&\det\left(\bar{\kappa}_n^\ast\left(\sum_{g\in G}[g\inv\eps_i,u_t]_{K,n}\piKKp(g)\right)\right)_{i,t=1}^d\nonumber\\
&=&\det\left(e^-\sum_{g\in G}\kappa_n(g)[g\inv\eps_i,u_t]_{K,n}g\inv\right)_{i,t=1}^d\nonumber
\end{eqnarray}
But  $\sum_{g\in G}\kappa_n(g)[g\inv\eps_i,u_t]_{K,n}g\inv$ clearly lies in the minus part and
Equation~(\ref{eq:Gaction2}) allows us to rewrite it as
$\sum_{g\in G}[\eps_i,gu_t]_{K,n}g\inv$. Thus we obtain simply
\beq\label{eq:simple formula for H_n}
H_{K/k,S,n}((1\otimes\eps_1)\wedge\ldots\wedge(1\otimes\eps_d),\wedgesub{u}{d})=
\det\left(\sum_{g\in G}[\eps_i,gu_t]_{K,n}g\inv\right)_{i,t=1}^d
\eeq
This shows in particular that on $\bigwedge^d_{\bbZ \bar{G}} U_{S}(K^+)\times {\textstyle \bigwedge_{\bbZ_p G}^d}U^1(K_p)$,
the pairing $H_{K,n}(\alpha(\cdot ),\cdot )$ agrees with
that defined by the pairing $\cH_{K,n}(\cdot ,\cdot )$ of~\cite{twizo}.\vertsp\\
\rem\label{rem: Lambdas as S increases} If $S\supset S'\supset S^0$ we shall always view the natural injection
${\textstyle \bigwedge^d_{\bbQ \bar{G}}\bbQ U_{S'}(K^+)}\rightarrow{\textstyle \bigwedge^d_{\bbQ \bar{G}}\bbQ U_{S}(K^+)}$
as an inclusion. It is then a simple exercise to check `compatibility of the pairings as $S$ varies' in the sense
that $\Lambda_{0,S}$ contains $\Lambda_{0,S'}$ and
$H_{K/k,S,n}$ agrees with $H_{K/k,S',n}$ on $\Zbp\Lambda_{0,S'}\times {\textstyle \bigwedge_{\bbZ_p G}^d}U^1(K_p)$. For this reason,
we shall usually omit the reference to $S$ and write simply $H_{K/k,n}$.
\subsection{The Map $\fs_{K/k,S}$}
For the time being we drop Hypothesis~\ref{hyp:S_contains_d+1_places} and the assumption that
$K$ contains $\mu_{\pnpo}$. We use the element $\aKkSm$ to define
a generalisation  of the map $\fs_{K/k}$ of~\cite{twizo} (slightly modified).
Let $j$ be any embedding of $\barbbQ$ into a fixed algebraic closure $\barbbQ_p$ of $\bbQ_p$. For each
$i=1,\ldots,d$, the composite $j\tau_i:\barbbQ\rightarrow\barbbQ_p$
defines a prime ideal $\fP_i$  of $\cO_K$ dividing $p$, namely
$\fP_i=\{a\in\cO_K\,:\,|j\tau_i(a)|_p<1\}$. (Of course the ideals
$\fP_1,\ldots,\fP_d$ are not in general distinct.) So $j\tau_i$
gives rise to an isometric embedding $K_{\fP_i}\rightarrow
\barbbQ_p$ (with the appropriately normalised $\fP_i$-adic metric on
$K_{\fP_i}$) whose image is  the topological closure
$\overline{j\tau_i(K)}$. This embedding is also denoted $j\tau_i$, by abuse.
There is a composite homomorphism of $\bbQ_p$-algebras
\[
\delta_i=\delta^{(j)}_i:=j\tau_i\circ \iota_{\fP_i}:K_p\rightarrow\barbbQ_p
\]
where $\iota_{\fP_i}: K_p\rightarrow K_{\fP_i}$ is as in the previous section.
It follows in particular that if $u$ lies in $U^1(K_p)\subset K_p$
then $|\delta^{(j)}_i(u)-1|_p<1$ for all $i$, hence the element
$\log_p(\delta^{(j)}_i(u))$ of $\overline{j\tau_i(K)}$ is given by
the usual logarithmic series. In
Proposition/Definition~\ref{propdefn:determinantal maps} we take
$\cR=\bbZ_p$, $\cS=\barbbQ_p$, $M=U^1(K_p)$, $H=G$, $l=d$ and set
$f_i(u):=\log_p(\delta^{(j)}_i(u))\ \forall\,u\in U^1,\
i\in\{1,\ldots,d\}$ to get a {\em $p$-adic regulator map}
$R_{K/k,p}:=\Delta_{\listsub{f}{d}}:{\textstyle
\bigwedge^d_{\bbZ_pG}U^1(K_p)}\rightarrow\barbbQ_p G$.
(We will denote it
$R_{K/k,p}^{(j)}$ or  $R_{K/k,p}^{(j;\tau_1,\ldots,\tau_d)}$ if we need to indicate the dependence on $j$ and/or $\tau_1,\ldots,\tau_d)$.)
For any
abelian group $H$ and commutative ring $\cR$ we define an involutive
automorphism $\underline{\ }^\ast$ of $\cR H$ by setting $(\sum
a_hh)^\ast=\sum a_hh\inv$. The element $a^-_{K/k,S}$ lies in $\barbbQ G^-$
by~(\ref{eq:generating QmufGminus}), hence so does
$a^{-,\ast}_{K/k,S}$ and applying $j$ to the coefficients we obtain
an element $j(a^{-,\ast}_{K/k,S})$ of $\barbbQ_p G^-$.
\begin{defn}\label{defn:s_Kk} For any $\theta\in \textstyle \bigwedge^d_{\bbZ_pG}U^1(K_p)$ we define $\fs_{K/k,S}(\theta)=\fs_{K/k,S,p}(\theta)$ to be the
product $j(a^{-,\ast}_{K/k,S})R_{K/k,p}^{(j)}(\theta)$ in $\barbbQ_p G$.
\end{defn}
\rem\label{rem:dependence of sKkS on taus}
It is easy to see that permuting the $\tau_i$ can only change the sign of the regulator $R_{K/k,p}^{(j)}$ and hence of the map \sKkS\
and that if $\tau_i$ is replaced by $\tau_i\tau$ for some $\tau\in\Gal(\barbbQ/k)$ then both are
multiplied by $\tau|_K\in G$. If clarity demands it we shall indicate this (simple) dependence on the $\tau_i$
by writing $\sKkS^{\dotlist{\tau}{1}{d}}$ instead of $\sKkS$. \vertsp\\
If $\fs_{K/k}$ denotes the map introduced in Definition~3.1 of~\cite{twizo}
then~(\ref{eq:aKkS_intermsof_PhiKko}) gives
\beq\label{eq:sKkS_intermsof_PhiKko_old_sKk}
\fs_{K/k,S}(\theta)=e^-\prod_{\fq\in S\setminus S^0} \left(1-N\fq^{-1}\sigma_\fq\right)j(\sqrt{d_k}\Phi_{K/k}(0)^\ast)R_{K/k,p}^{(j)}(\theta)=
e^-\prod_{\fq\in S\setminus S^0} \left(1-N\fq^{-1}\sigma_\fq\right)\fs_{K/k}(\theta)
\eeq
and if $k\neq \bbQ$ then we can even drop the factor $e^-$.
Equation~(\ref{eq:sKkS_intermsof_PhiKko_old_sKk}) and Proposition~3.4 of~\emph{ibid.} imply the important
\begin{prop}\label{prop:about sKkS}
$\fs_{K/k,S}(\theta)$ lies in $\bbQ_p G^-$ for every $\theta\in \textstyle \bigwedge^d_{\bbZ_pG}U^1(K_p)$. Moreover it
is independent of the choice of $j$.\ePf
\end{prop}
In~\cite{twizo}, $\fs_{K/k}$ was considered as a ($\bbZ_pG$-linear) map from $\bigwedge^d_{\bbZ_pG}U^1(K_p)$ to $\bbQ_pG$.
But because of the factor $e^-$ in~(\ref{eq:sKkS_intermsof_PhiKko_old_sKk}), we now have
$\fs_{K/k,S}(e^-\theta)=e^-\fs_{K/k,S}(\theta)=\fs_{K/k,S}(\theta)$.
For this reason, we prefer to consider $\fs_{K/k,S}$
as a $\bbZ_p G$-linear map from $\textstyle \bigwedge^d_{\bbZ_pG}U^1(K_p)^-$ to $\bbQ_p G^-$.
\begin{prop}\label{prop:ker and im of sKks}
The kernel of $\fs_{K/k,S}$ is precisely the ($\bbZ_p$-) torsion submodule of $\textstyle \bigwedge^d_{\bbZ_pG}U^1(K_p)^-$ which is finite. The image
of $\fs_{K/k,S}$ is a fractional ideal of $\bbQ_pG^-$ (\ie\ a finitely generated $\bbZ_p G$-submodule of $\bbQ_pG^-$ which spans it over $\bbQ_p$).
\end{prop}
\bPf\ In Remark~3.2 of~\cite{twizo} it was shown that
$\ker(R_{K/k,p}^{(j)})$ is finite and that $\im(R_{K/k,p}^{(j)})$ spans $\barbbQ_p G$ over $\barbbQ_p$. Also,
Equation~(\ref{eq:generating QmufGminus}) implies that $j(a^{-,\ast}_{K/k,S})$
is a unit of the ring $\barbbQ_pG^-$. It follows that $\ker(\fs_{K/k,S})$ lies in $\ker(R_{K/k,p}^{(j)})$ and hence in $\left(\textstyle \bigwedge^d_{\bbZ_pG}U^1(K_p)^-\right)_{\rm tor}$. The reverse inclusion is clear, since
$\bbQ_p G^-$ is torsion-free. For the second statement, finite-generation  follows from that of $U^1(K_p)$ and we have
$
\barbbQ_p\im(\fs_{K/k,S})=\barbbQ_pG^-\im(\fs_{K/k,S})=\barbbQ_pG^-\im(R_{K/k,p}^{(j)})=\barbbQ_pG^-
$
Since $\bbQ_pG^-$ contains a $\barbbQ_p$-basis of $\barbbQ_pG^-$, it follows that
$\bbQ_p\im(\fs_{K/k,S})=\bbQ_pG^-$.\ePf
\begin{defn}\label{defn:S_Kk} We set $\fS_{K/k,S}=\fS_{K/k,S,p}:=\im(\fs_{K/k,S,p})\subset\bbQ_p G^-$. (Proposition~\ref{prop:about sKkS}
and Remark~\ref{rem:dependence of sKkS on taus} show that $\SKkS$ is independent
of $j$ and the choice and ordering of the $\tau_i$'s.)
\end{defn}
Thus
\beq\label{eq:SKkS_intermsof_old_SKk}\fS_{K/k,S}=
e^-\prod_{\fq\in S\setminus S^0} \left(1-N\fq^{-1}\sigma_\fq\right)\fS_{K/k}
\eeq
where $\fS_{K/k}=\im(\fs_{K/k})$ as in~\cite{twizo}, and if $k\neq \bbQ$ then we can drop the factor $e^-$.
Finally, the dependence of $\fs_{K/k,S}$ and $\fS_{K/k,S}$ on $S$ is clear:
if $S\supset S'\supset S^0$ then~(\ref{eq:thetaS'toS}) and the definition of $a^-_{K/k,S}$
give
\beq\label{eq:sKkS'to_sKkS_and SKkS'to_SKkS}
\fs_{K/k,S}=\prod_{\fq\in S\setminus S' }\left(1-N\fq^{-1}\sigma_\fq\right)\fs_{K/k,S'}\ \ \ \mbox{and}\ \ \
\fS_{K/k,S}=\prod_{\fq\in S\setminus S' }\left(1-N\fq^{-1}\sigma_\fq\right)\fS_{K/k,S'}
\eeq
\section{Statements of the Conjectures}\label{section:conjectures}
Let us write $S_p$ for $S_p(k)$ and $S^1=S^1(K/k)$ for $S_p\cup S^0=S_p\cup \Sram(K/k)\cup S_\infty$.
\begin{hyp}\label{hyp:S contains S1}
$S$ contains $S^1$
\end{hyp}
Henceforth, the three conditions $p\neq 2$, Hypothesis~\ref{hyp:K_CM_k_totreal} and Hypothesis~\ref{hyp:S contains S1} will be referred to as {\em `the standard hypotheses'}
and {\em  will be assumed to hold unless it is explicitly stated otherwise.}
Our `Integrality Conjecture' (IC) reads:
\begin{conjalt}{$IC(K/k,S,p)$ }\label{conj:IC}
$\SKkS\subset\bbZ_pG^-$.
\end{conjalt}
\rem\ By using~\cite[Cor.~2.1]{twizo} and estimates of $\log_p$ one can find explicit values of $N$ such that $\SKkS\subset  p^{-N}\bbZ_pG^-$ ({\em cf.} the proof of Prop.~4.2, {\em ibid.}).
The conjecture says we can take $N=0$. Fixing $K/k$ but letting $p$ (hence $S^1$) vary, one can also show that $\fS_{K/k,S^1,p}=\bbZ_pG^-$ for all but finitely many $p\neq 2$. In fact, this follows easily from
Theorem~\ref{thm:chi(SKkS) when p ndiv G}.
\vspace{1ex}\\
\rem\label{rem:new_IC_from_old_IC}\
Equation~(\ref{eq:SKkS_intermsof_old_SKk}) gives
\[
\fS_{K/k,S^1}=
e^-\prod_{\fp\in S_p\setminus\Sram} \left(1-N\fp^{-1}\sigma_\fp\right)\fS_{K/k}=
e^-\left(\prod_{\fp\in S_p\setminus\Sram}N\fp\right)\inv\fS_{K/k}
\]
(For the second equality, observe that if $\fp$ lies in $S_p\setminus\Sram$ then $N\fp-\sigma_\fp$ is a unit of $\bbZ_pG$
.) If $k\neq \bbQ$ we may,
as usual, drop the factor $e^-$ in the last equation. It follows in particular that if $k\neq \bbQ$ then $IC(K/k,S^1,p)$
is equivalent to Conjecture~5.2 of~\cite[\S~5.2]{twizo}. If $k=\bbQ$ the latter conjecture was proven in~\emph{ibid.}.
$IC(K/\bbQ,S^1,p)$ follows on applying $e^-$  and will be re-proven in
Theorem~\ref{thm: CC and IC for K over Q}.\vspace{1ex}\\
Hypothesis~\ref{hyp:S contains S1} implies Hypothesis~\ref{hyp:S_contains_d+1_places} so that the conditions of
Conjecture $RSC(K^+/k,S;\Zbp)$ are met. Our `Congruence Conjecture' (CC) reads:
\begin{conjalt}{$CC(K/k,S,p,n)$ (Congruence Conjecture)}\label{conj:CC}\\
Suppose that Conjecture $IC(K/k,S,p)$ holds and that $RSC(K^+/k,S;\Zbp)$ holds with solution $\eta_{K^+/k,S}$.
If also $K\supset\mu_\pnpo$ for some $n\geq 0$ then, for all $\theta\in\WUoKpm$ we have
\beq\label{eq:the congruence}
\overline{\sKkS(\theta)}=\kappa_n(\tau_1\ldots\tau_d)H_{K/k,n}(\eta_{K^+/k,S},\theta)
\mbox{\ \ \ (in $(\ZpnpoZ) G^-$)}
\eeq
\end{conjalt}
\vspace{-2ex}
\rem\label{rem:CC 1}
The factor $\kappa_n(\tau_1\ldots\tau_d)$ means that the Congruence Conjecture is independent of the choice of $\tau_1,\ldots,\tau_d$. For example,
if we replace $\tau_i$ by $\tau_i\tau\inv$ for some $\tau\in\Gal(\barbbQ/k)$ then
Remark~\ref{rem: changing taus}, Proposition~\ref{prop:pairing_H}~\ref{part1:pairing_H} and~(\ref{eq:formulae for kappanbarstar})
show that the R.H.S. is multiplied by $\tau|_K\inv\in G$ and the same is true for the
L.H.S. by Remark~\ref{rem:dependence of sKkS on taus}.
\vspace{1ex}\\
\rem\label{rem:CC 2}
$CC(K/k,S,p,n)$ replaces Conjecture~5.4 of~\cite{twizo}.
The latter is essentially the special case of the CC in which $S=S^1$ and  $p$ splits in $k$ (so that $\mu_p\subset K$ forces $S^0=S^1$). In fact, it is a direct
consequence of this case provided one assumes (with no significant loss of generality) that $K$ is CM, $k\neq \bbQ$ and
one replaces ${\textstyle\bigwedge^d_{\bbZ G}K^\times}$ in  Conjecture~5.4 with
${\textstyle\bigwedge^d_{\bbZ \bar{G}}U_S(K^+)}$ as here.
The awkwardness in the formulation of
Conjecture~5.4 (using $\cI(\eta^+_{K/k})$, $\tilde{\eta}_x$ \etc.) has been avoided in the CC thanks to our `extension' of $H_{K,n}$ to $\Zbp\Lambda_{0,S}$.
\section{Evidence for the IC and the CC}\label{sec: evidence}
\subsection{The Results of~\cite{twizo}}
Conjecture~5.2 of~\cite{twizo} implies $IC(K/k,S,p)$ for $S=S^1$ (see Remark~\ref{rem:new_IC_from_old_IC}) and hence for all $S$ by~Prop.\ref{prop: IC under increasing S}. Therefore
Proposition~4.2 of~\cite{twizo} translates to
\begin{thm}\label{thm:oldcase 1_of_IC}
$IC(K/k,S,p)$ holds  whenever $p$ is unramified in $K/\bbQ$.\ePf
\end{thm}
Similarly, the main result, Theorem~4.1, of~\cite{twizo}
gives
\begin{thm}\label{thm:oldcase 2_of_IC}
$IC(K/k,S,p)$ holds whenever $p$ is splits completely in $k/\bbQ$ and either $\Sram(K/k)\not\subset S_p(k)$ or $\mu_p(K)=\{1\}$. \ePf
\end{thm}
\subsection{The IC and the ETNC}
Working on the original version of the IC in~\cite{twizo}, Andrew Jones has shown that a certain refinement would follow from the ETNC
(see the Introduction). Let  $\Cl_{\fm}(K)$ be the
ray-class group of $K$ corresponding to the cycle which is the formal product of the finite
places of $K$ above those in $S^1$ and write ${\rm Fitt}_{\bbZ G}(\Cl_{\fm}(K))$
for its (initial) Fitting ideal as a $\bbZ G$-module. In our notation, the first part of~\cite[Theorem 4.1.1]{Jones' thesis} then says that the relevant case of the ETNC
(namely~\cite[Conj.~4.(iv)]{Bu=[8] of Jones}  for the pair $(h^0{\rm Spec}(K)(1), e^-\bbZ G )$) implies
\beq\label{eq: Jones' main thm}
\fS_{K/k,S^1}\subset(\bbZ_p{\rm Fitt}_{\bbZ G}(\Cl_{\fm}(K)))^-\ (={\rm Fitt}_{\bbZ_p G^-}((\Cl_{\fm}(K)\otimes \bbZ_p)^-))
\eeq
for all odd primes $p$.
The inclusion~(\ref{eq: Jones' main thm}), hence the ETNC, clearly implies $IC(K/k,S,p)$  for $S=S^1$, hence for all $S$. Of course, it implies considerably more (for instance,
that $\fS_{K/k,S^1}$ annihilates the $p$-part of  $\Cl_{\fm}(K)$) and in this sense refines the IC in a different direction
from the CC. We note that the relevant case of the ETNC  has been proven in our set-up only when $K$ is an absolutely abelian field (see below).
\subsection{Strengthenings of the IC in the Case $p\ndiv |G|$ }
The second part of Jones' Theorem 4.1.1 states that if the above case of the ETNC holds and also $p\ndiv |G|$,
then we have the following strengthening of~(\ref{eq: Jones' main thm})
\beq\label{eq: Jones on ETNC and pndivG}
\fS_{K/k,S^1}=
\eitheror
{{\rm Fitt}_{\bbZ_p G^-}\left(\mu_{p^\infty}(K_p)^-/\mu_{p^\infty}(K)\right){\rm Fitt}_{\bbZ_p G^-}((\Cl_{\fm}(K)\otimes \bbZ_p)^-)}
{if $\Sram(K/k)\subset S_p$}
{{\rm Fitt}_{\bbZ_p G^-}\left(\mu_{p^\infty}(K_p)^-\right){\rm Fitt}_{\bbZ_p G^-}((\Cl_{\fm}(K)\otimes \bbZ_p)^-)}
{if $\Sram(K/k)\not\subset S_p$}
\eeq
Corollary~4.1.8 of~\cite{Jones' thesis}  also establishes~(\ref{eq: Jones on ETNC and pndivG})
when $p\ndiv |G|$ {\em without} assuming the ETNC but imposes a mild condition on the characters of $G$.
(The proof uses results of~\cite{W=[46] of Jones} and work of Bley, Burns and others on, roughly speaking,
the compatibility of the ETNC with the functional equations of $L$-functions.)

Independently, we used the functional equations themselves and more elementary, index-type arguments to give a different (and unconditional)
formula for $\fS_{K/k,S}$ whenever $p\ndiv|G|$. This is presented as Theorem~\ref{thm:chi(SKkS) when p ndiv G} in Section~\ref{sec: pndivG}.
Corollary~\ref{cor: IC when pndivG} shows how one may quickly deduce $IC(K/k,S,p)$ in this case. Of course, it would also
follow immediately from Jones' formula~(\ref{eq: Jones on ETNC and pndivG}). In fact, there is a direct link between the two formulae, explained in
Remark~\ref{rem: Jones' formula from ours}.

\subsection{The Case $k=\bbQ$}
When $k=\bbQ$, the IC follows from Corollary~4.1 of~\cite{twizo} (or indeed from the work of Jones, see below).
In Section~\ref{sec: k equals Q} we shall prove the CC in this case, re-proving the IC along the way:
\begin{thm}\label{thm: CC and IC for K over Q}
\begin{enumerate}
\item\label{part1: CC and IC for K over Q} Conjecture $IC(K/\bbQ,S,p)$ holds and
\item\label{part2: CC and IC for K over Q} If $K$ contains $\mu_\pnpo$ for some $n\geq 0$, then Conjecture $CC(K/\bbQ,S,p,n)$ holds.
\end{enumerate}
\end{thm}
\subsection{The Case of Absolutely Abelian $K$}\label{subsec: case of absab preview}
As noted in~\cite[Cor 4.1.7]{Jones' thesis}, the relevant case of the ETNC
follows from~\cite[Cor 1.2]{BF2} whenever
$K$ is absolutely abelian and $k$ is any totally real subfield (possibly but not necessarily equal to $\bbQ$). Thus the inclusion~(\ref{eq: Jones' main thm}) holds and in particular
\begin{thm}\label{thm: IC for K absab}
If $K$ is an abelian extension of $\bbQ$, then Conjecture $IC(K/k,S,p)$ holds.\ePf
\end{thm}
To state our result for the CC, in this case, we first define a set
of rational primes ${\rm Bad}(S):=\{q\in
\Sram(k/\bbQ)\,:\,S_q(k)\not\subset S\}$ and formulate
\begin{hyp}\label{hyp: badS etc}
$p\ndiv e_q(k/\bbQ)$ for all $q\in\badS$
\end{hyp}
In Section~\ref{sec: Kabsab} we shall show
\begin{thm}\label{thm: CC for Kabsab}
If $K$ is an abelian extension of $\bbQ$ containing $\mu_\pnpo$ for
some $n\geq 0$ and Hypothesis~\ref{hyp: badS etc} is satisfied, then
Conjecture $CC(K/\bbQ,S,p,n)$ holds.
\end{thm}
(At the same time we shall obtain a second proof of Theorem~\ref{thm: IC for K absab}
which assumes Hypothesis~\ref{hyp: badS etc} but is independent of the ETNC.)
The proof of~Theorem~\ref{thm: CC for Kabsab} uses induction formulae for $L$-functions to relate the situation for $K/k$ to that of
$F/\bbQ$ for various CM subfields $F$ of $K$, and this in two parallel applications. The first concerns $\fs_{K/k,S}$ and works at $s=1$. The second concerns $RSC(K^+/k;\Zbp)$
and works at $s=0$.
Popescu introduced the latter application in~\cite{Pop}. (In fact, he applied it to his own variant of Rubin's conjecture $B'$ which also implies $RSC(K^+/k,S;\bbZ)$.) He worked
under a hypothesis which implies $\badS=\varnothing$. This simplifies matters (we only need consider $F=K$)
but is rather restrictive (\eg\ ${\rm Bad}(S^1(K/k))\neq\varnothing$ whenever a rational prime $q\neq p$ ramifies in $k/\bbQ$ but not in $K/k$).
The elaboration of Popescu's techniques which allows us to conclude under our
weaker Hypothesis~\ref{hyp: badS etc} is one ingredient of Cooper's work on Popescu's Conjecture in~\cite{Cooper}.
Hypothesis~\ref{hyp: badS etc} holds, for example, whenever~$p\ndiv [k:\bbQ]$ (\eg\ $[k:\bbQ]$ is a power of $2$). Alternatively, suppose  $K=\bbQ(\xi_f)$ and
$k=K^+$ where $f=p^{n+1}f'\not\equiv 2\pmod{4}$, $n\geq 0$ and $p\ndiv f'$. If we take $S=S^1(K/k)=S_{\infty}\cup S_p$ then $\badS=\varnothing\Leftrightarrow f'=1$ but
Hypothesis~\ref{hyp: badS etc} holds provided only $p\ndiv q-1$ for all $q|f'$.
\subsection{Two `Trivial' Cases of the Congruence~(\ref{eq:the congruence})}
Suppose $K\supset \mu_\pnpo$ for some $n\geq 0$ and $S$ contains at least $d+2$ places and at least one finite place
$\fq$  that splits completely in $K^+$.
Equations~(\ref{eq:defn_rschi}) and~(\ref{eq:ordLKpluskSatzero}) imply $\Theta_{K^+/k,S}^{(d)}=0$
so that $RSC(K^+/k,S;\bbZ)$ holds
with $\eta_{K^+/k,S}=0$. The congruence~(\ref{eq:the congruence})
is thus equivalent to
$\sKkS(\theta)\in\pnpo \bbZ_p G^-$. The extension
$K/K^+=K^+(\mu_\pnpo)/K^+$ is unramified outside $p$, so  if $\fq$ does not divide $p$ then it
cannot lie in $S^1$ (which forces $|S|\geq d+2$). We can then apply the following result. (For a case with $\fq|p$, see the next subsection.)
\begin{prop}\label{prop: trivial CC 1}
Suppose $K\supset \mu_\pnpo$ and $\fq\in S\setminus S^1(K/k)$ splits in $K^+$.
If $IC(K/k,S\setminus\{\fq\},p)$ holds (\eg\ if $p\ndiv |G|$) then $\SKkS\subset\pnpo\bbZp G^-$.
In particular, $CC(K/k,S,p,n)$ holds.\ePf
\end{prop}
\bPf\ By equation~(\ref{eq:sKkS'to_sKkS_and SKkS'to_SKkS}) it
clearly suffices to show that $\pnpo$ divides
$(N\fq-\sigma_\fq)e^-$. Since $\fq$ splits in $K^+$, $\sigma_\fq$ is
either $1$ or $c$ and since $\fq\ndiv p$, it acts on $\mu_\pnpo$ by
$N\fq$. If $\sigma_\fq=1$, it also acts trivially, so $\pnpo$
divides $N\fq-1$ hence also $(N\fq-1)e^-=(N\fq-\sigma_\fq)e^-$. If
$\sigma_\fq=c$, it also acts by $-1$, so $\pnpo$ divides $N\fq+1$
hence also $(N\fq+1)e^-=(N\fq-\sigma_\fq)e^-$.\ePf \noindent  Next,
suppose $\theta$ is a $\bbZ_p$-torsion element in
$\bigwedge^d_{\bbZ_p G}U^1(K_p)^-$. Proposition~\ref{prop:ker and im
of sKks} implies that the  L.H.S. of~(\ref{eq:the congruence})
vanishes, so, assuming $RSC(K^+/k,S;\Zbp)$, this congruence is
equivalent to $H_{K/k,n}(\eta_{K^+/k,S},\theta)=0$. If also $p\ndiv
|G|$, then this is an immediate consequence of the following result,
to be proved in Section~\ref{sec: pndivG}. (The verification  seems
harder if $p||G|$ (and $d\geq 2$), not least because
$\left(\bigwedge^d_{\bbZ_p G}U^1(K_p)\right)_{\rm tor}$ is then
harder to characterise.)
\begin{prop}\label{prop: trivial CC 2}
Suppose $p\ndiv |G|$, $K\supset \mu_\pnpo$  and  $\eta$ is any element of $\Zbp\Lambda_{0,S}(K^+/k)$ satisfying the eigenspace condition w.r.t.\ $(S,d,\bar{G})$. Then
$
H_{K/k,n}(\eta,\theta)=0$
for all $\theta\in \left(\bigwedge^d_{\bbZ_p G}U^1(K_p)\right)_{\rm tor}$.
\end{prop}

\subsection{The Case $k=K^+$}
In this case $G=\{1,c\}$ so  $p\ndiv |G|$ and the IC holds for all admissible $S$.
For the CC, we assume  $K\supset \mu_\pnpo$ with  $n\geq 0$ so that
$K=k(\mu_\pnpo)$ and $S^1=S_\infty\cup S_p$. All places of $k$ split in $K^+$ so if $S\neq S^1$ then $CC(K/k,S,p,n)$ holds by
Proposition~\ref{prop: trivial CC 1}.
Also, if $|S_p(k)|\geq 2$ then $|S|\geq r+2$ so once again $CC(K/k,S^1,p,n)$ is equivalent to $\SKkS\subset\pnpo\bbZp G^-$ (see above). But this
will follow from equation~(\ref{eq: reformulation of thm pndivG}) in Section~\ref{sec: pndivG} (for the unique odd character $\phi$.
Indeed, the first term on the RHS of~(\ref{eq: reformulation of thm pndivG}) is clearly divisible by $(\pnpo)^{|S_p|-1}$.)
This leaves only the case $S=S_\infty\cup S_p$ with $|S_p(k)|=1$. Then $\eta_{k/k,S^1}$ is non-zero and
can be written explicitly in terms of a $\bbZ$-basis $\eps_1,\ldots,\eps_{d}$ for $U_{S^1}(k)/\{\pm 1\}$ and the
$S_p$-classnumber of $k$. In this case
$CC(K/k,S^1,p,n)$ reduces to an apparently novel identity in $\bbZ/\pnpo\bbZ$, relating a $p$-adic regulator of
elements of $U^1(K_p)^-$ to a determinant of their Hilbert symbols with the $\eps_i$.
In ongoing work, Malcolm Bovey has done extensive computations verifying this identity
and established some
partial results supporting it.
\subsection{Other Computational Results}

$RSC(K^+/k,S;\bbQ)$ is not currently known to hold non-trivially for any $S$ unless
either $K^+$ is absolutely abelian or all the characters $\chi\in\hat{\bar{G}}$ satisfying
$\ord_{s=0}L_{K^+/k, S}(s,\chi)=d$ are of order $1$ or $2$.
However, if $d$ is not too large, high-precision computation can
identify $\eta_{K^+/k,S}$ with virtual certainty as the unique
solution of~(\ref{eq:rubin_stark}) in $e_{S,d,\bar{G}}\textstyle \bigwedge^d_{\bbQ \bar{G}}\bbQ U_{S}(K^+)$. (This was done in~\cite{Solomon-Roblot1}.)
This makes it possible to check the CC (and simultaneously the IC) on a computer.
In~\cite{Solomon-Roblot2} we give details of such numerical verifications for more than $40$
cases of $CC(K/k,S^1,p,n)$ with $k$ real quadratic, $n=0$ or $1$ and varying $K$ and $p$.
\section{Changing $S$, $K$ and $n$}\label{sec: functoriality}
If $\fq$ is a prime ideal of $k$ not in $S_p$ then $\left(1-N\fq^{-1}\sigma_\fq\right)$ lies in $\bbZ_p G$.
Hence equation~(\ref{eq:sKkS'to_sKkS_and SKkS'to_SKkS}) gives
\begin{prop}\label{prop: IC under increasing S}
If $S\supset S'\supset S^1$ then $IC(K/k, S',p)$ implies $ IC(K/k, S,p)$. \ePf
\end{prop}
\rem\ For the converse implication one would need
$e^-\left(1-N\fq^{-1}\sigma_\fq\right)$ to be invertible in $\bbZ_pG^-$ for each $\fq\in S\setminus S'$.
But for any such $\fq$ one has an isomorphism of $\bbZ_p G$-modules
\beq\label{eq: galois iso for res fld}
\bbZ_p G/\left(1-N\fq^{-1}\sigma_\fq\right)\cong (\cO_K/\fq\cO_K)^\times\otimes_\bbZ\bbZ_p\cong \bigoplus_\fQ(\bbF_\fQ^\times\otimes_\bbZ\bbZ_p)
\eeq
where $\fQ$ runs through the primes dividing $\fq$ in $K$ and $\bbF_\fQ$ denotes $\cO_K/\fQ$. We deduce that if $c\nin D_\fq(K/k)$ (respectively, $c\in D_\fq(K/k)$) then
$e^-\left(1-N\fq^{-1}\sigma_\fq\right)$ lies in $(\bbZ_pG^-)^\times$ if and only if $p\ndiv|\bbF_\fQ^\times|$ (respectively, $p\ndiv|\bbF_\fQ^\times/(1+c)\bbF_\fQ^\times|$)
for one, hence any $\fQ$.
This fails in particular if $\mu_p\subset K$, in which case $IC(K/k, S,p)$ does not by itself imply $IC(K/k, S',p)$ for any
$S\supsetneq S'$.
\begin{prop}\label{prop: RSC under increasing S}
Suppose $S\supset S'\supset S^1$ and $RSC(K^+/k, S';\Zbp)$ holds with solution $\eta_{K^+/k, S'}$. Then $RSC(K^+/k, S;\Zbp)$ holds with solution
$\eta_{K^+/k, S}=\prod_{\fq\in S\setminus S'}(1-\sigma_{\fq,K^+}\inv)\eta_{K^+/k, S'}$.
\end{prop}
\bPf\ It follows easily from~(\ref{eq:thetaS'toS}) and
characterisation~\ref{part3: equivalent `eigenspace' conditions} of the eigenspace condition that
$\prod_{\fq}(1-\sigma_{\fq,K^+}\inv)\eta_{K^+/k, S'}$ is a solution of $RSC(K^+/k,S;\bbQ)$. The result follows since
$\Zbp\Lambda_{0,S'}$ is a $\bbZ\bar{G}$-submodule of $\Zbp\Lambda_{0,S}$.\ePf
\begin{prop}\label{prop: CC under increasing S}
If $K\supset\mu_{\pnpo}$ for some $n\geq 0$ and $S\supset S'\supset S^1$ then $CC(K/k,S',p,n)$ implies $ CC(K/k,S,p,n)$. \ePf
\end{prop}
\bPf\ We assume that $CC(K/k,S',p,n)$ holds  so also $IC(K/k,S',p)$ and $RSC(K^+/k,S';\Zbp)$. Thus  $IC(K/k,S',p)$ and
$RSC(K^+/k,S';\Zbp)$ hold by Props.~\ref{prop: IC under increasing S} and~\ref{prop: RSC under increasing S}.
Using the latter and~Proposition~\ref{prop:pairing_H}~\ref{part1:pairing_H} we find
that for any $\theta\in\WUoKpm$
\begin{eqnarray}
\kappa_n(\tau_1\ldots\tau_d)H_{K/k,n}(\eta_{K^+/k, S},\theta)&=&
\prod_{\fq\in S\setminus S'}(1-\bar{\kappa}_n^\ast(\sigma_{\fq,K^+}\inv))\kappa_n(\tau_1\ldots\tau_d)H_{K/k,n}(\eta_{K^+/k, S'},\theta)\nonumber\\
&=&\prod_{\fq\in S\setminus S'}(1-\bar{\kappa}_n^\ast(\sigma_{\fq,K^+}\inv))\overline{\fs_{K/k,S'}(\theta)}\mbox{\ in $(\ZpnpoZ) G^-$}\label{eq: H on increasing S}
\end{eqnarray}
For each $\fq\in S\setminus S'$, equation~(\ref{eq:formulae for kappanbarstar}) with $g=\sigma_\fq\inv=\sigma_{\fq,K}\inv$
gives $\bar{\kappa}_n^\ast(\sigma_{\fq,K^+}\inv)=e^-\kappa_n(\sigma_{\fq})\inv\sigma_{\fq}$ and since $\fq\ndiv p$\ it follows that
$\kappa_n(\sigma_{\fq})=\overline{N\fq}$ in $\ZpnpoZ$.
Thus $1-\bar{\kappa}_n^\ast(\sigma_{\fq,K^+}\inv)$ acts as $1-N\fq\inv\sigma_{\fq}$ on $(\ZpnpoZ) G^-$ and combining~(\ref{eq: H on increasing S})
with~(\ref{eq:sKkS'to_sKkS_and SKkS'to_SKkS}) gives~(\ref{eq:the congruence}), as required.\ePf
\noindent Now suppose that $F$ is any CM subfield of $K$ containing $k$. Then $p$, $F$ and $S$ satisfy the standard hypotheses.
We write $G_F$ for $\Gal(F/k)$ and $N_{K/F}$ for the norm map $K_p\rightarrow F_p$.
(If we identify $K_p$ with $\prod_{\fP|p} K_\fP$ and $F_p$ with $\prod_{\fp|p} F_\fp$, then $N_{K/F}$ sends $(x_\fP)_\fP$ to $(y_\fp)_\fp$, where
$y_\fp=\prod_{\fP|\fp}N_{K_\fP/F_\fp}x_\fP$). We shall also write $N_{K/F}$ for the $\bbZ_p$-linear map
${\textstyle \bigwedge^d_{\bbZ_pG}U^1(K_p)}\rightarrow {\textstyle \bigwedge^d_{\bbZ_pG_F}U^1(F_p)}$
sending $u_1\wedge\ldots\wedge u_d$ to $N_{K/F}u_1\wedge\ldots\wedge N_{K/F}u_d$.
One checks easily that $\pi_{K/F}\circ R^{(j)}_{K/k,p}=R^{(j)}_{F/k,p}\circ N_{K/F}$ and also that
$\pi_{K/F}\circ\Theta_{K/k,S}=\Theta_{F/k,S}$ (as meromorphic functions $\bbC\rightarrow\bbC G_F$)
so that $\pi_{K/F}(a^-_{K/k,S})=a^-_{F/k,S}$.
We deduce easily
\begin{prop}\label{prop: IC under descending K}
If $K\supset F\supset k$ as above then $\pi_{K/F}\circ \fs_{K/k,S}=\fs_{F/k,S}\circ N_{K/F}$. In particular, if
$N_{K/F}:{\textstyle \bigwedge^d_{\bbZ_pG}U^1(K_p)}^-\rightarrow{\textstyle \bigwedge^d_{\bbZ_pG_F}U^1(F_p)}^-$
is surjective then $\pi_{K/F}(\fS_{K/k,S})=\fS_{F/k,S}$, so $IC(K/k,S,p)$ implies $IC(F/k,S,p)$. \ePf
\end{prop}
\rem\label{rem: surjectivity}\ The surjectivity condition is certainly satisfied whenever $N_{K/F}(U^1(K_p))=U^1(F_p)$. The latter condition holds iff
$K/F$ is at most tamely ramified at each prime in $S_p(F)$ (by local class-field theory). Of course, it actually suffices that $N_{K/F}(U^1(K_p)^-)=U^1(F_p)^-$
which can be shown to be equivalent to the statement {\em `$K/F^+$ is at most tamely ramified at each  prime in $S_p(F^+)$ which splits in $F$'.}
One can also show that $\nu_{K/F}\circ\fs_{F/k,S}=|G|^{1-d}\fs_{K/k,S}\circ i_{K/F}$ where $i_{K/F}$ is the natural map
${\textstyle \bigwedge^d_{\bbZ_pG_F}U^1(F_p)}^-\rightarrow{\textstyle \bigwedge^d_{\bbZ_pG}U^1(K_p)}^-$, but for present purposes this is only helpful when $p\ndiv |G|$ or $d=1$.
\vertsp\\
Let $\bar{G}_F=\Gal(F^+/k)$. The norm  $N_{K^+/F^+}$  maps $U_S(K^+)$ into $U_S(F^+)$. We also use the symbol `$N_{K^+/F^+}$' to denote the map
$1\otimes N_{K^+/F^+}:\bbQ U_S(K^+)\rightarrow\bbQ U_S(F^+)$ and the $\bbQ$-linear map
${\textstyle \bigwedge^d_{\bbQ \bar{G}}\bbQ U_S(K^+)}\rightarrow {\textstyle \bigwedge^d_{\bbQ \bar{G}_F}\bbQ U_S(F^+)}$
sending $x_1\wedge\ldots\wedge x_d$ to $N_{K^+/F^+}x_1\wedge\ldots\wedge N_{K^+/F^+}x_d$.
\begin{prop}\label{prop: RSC(Q) under descending K}
Suppose  $K\supset F\supset k$ as above and $RSC(K^+/k, S; \bbQ)$ holds with solution $\eta_{K^+/k, S'}$. Then $RSC(F^+/k,S; \bbQ)$ holds with solution
$\eta_{F^+/k, S}=N_{K^+/F^+}\eta_{K^+/k, S}$.
\end{prop}
\bPf\ $\pi_{K^+/F^+}(N_{D_\fq(K^+/k)})$ is a $\bbZ$-multiple of $N_{D_\fq(F^+/k)}$ for all $\fq\in S\setminus S_\infty$. It follows easily from this
that form~\ref{part3: equivalent `eigenspace' conditions} of the eigenspace condition on $\eta_{K^+/k, S}$  (w.r.t.\ $(S,d,\bar{G})$) implies the same on
$N_{K^+/F^+}\eta_{K^+/k, S}$ (w.r.t.\ $(S,d,\bar{G}_F)$). Similarly, since $\pi_{K^+/F^+}\circ\Theta_{K^+/k,S}=\Theta_{F^+/k,S}$ and
$\pi_{K^+/F^+}\circ R_{K^+/k}=R_{F^+/k}\circ N_{K^+/F^+}$,
if we apply $\pi_{K^+/F^+}$ to condition~(\ref{eq:rubin_stark}) for $\eta_{K^+/k, S}$  then we get the equivalent condition on
$N_{K^+/F^+}\eta_{K^+/k, S}$.\ePf
\noindent Before attacking the Congruence Conjecture in this context, we need two Lemmas.
\begin{lemma}\label{lemma: Norming down Lambda}
If $d=1$ then $N_{K^+/F^+}(\Lambda_{0,S}(K^+/k))\subset \Lambda_{0,S}(F^+/k)$. If $d>1$, then
$N_{K^+/F^+}(\Lambda_{0,S}(K^+/k))$ is contained in
$e^{-d}\Lambda_{0,S}(F^+/k)$ where $e=\exp((U_S(K^+)/U_S(F^+))_{\rm tor})=1$ or $2$.
\end{lemma}
\bPf\ The first statement follows from that of  Prop.~\ref{prop:index of barLdU in Lambda}. Next,
$(U_S(K^+)/U_S(F^+))_{\rm tor}$ injects into ${\rm Hom}(\Gal(K^+/F^+),\mu(K^+))={\rm Hom}(\Gal(K^+/F^+),\{\pm 1\})$
(by sending $[\eps]$ to the map $g\mapsto g(\eps)/\eps$)
so $e=1$ or $2$. Let
$(U_S(K^+)/U_S(F^+))_{\rm tor}=V/U_S(F^+)$ where $U_S(K^+)\supset V\supset U_S(F^+)$.
Since $U_S(K^+)/V$ is torsionfree, $U_S(K^+)$ splits over $\bbZ$ as $V\oplus V'$. The sum $U_S(F^+)+V'$ is also direct and contains $U_S(K^+)^e$.
Therefore, any $f_1,\ldots,f_d$
lying in $\Hom_{\bbZ}(U_S(F^+),\bbZ)$ (considered as a subset of $\Hom_{\bbQ}(\bbQ U_S(F^+),\bbQ)$) extend to
$\hat{f}_1,\ldots,\hat{f}_d$ in $\Hom_{\bbZ}(U_S(F^+)+V',\bbZ)$ considered as a subset of $\Hom_{\bbQ}(\bbQ U_S(K^+),\bbQ)$ and
$e\hat{f}_i\in \Hom_{\bbZ}(U_S(K^+),\bbZ)$ for all $i$. It is easy to see from the definitions that
\beq\label{eq: Delta and pi}
\pi_{K^+/F^+}(\Delta_{\hat{f}_1,\ldots\hat{f}_d}(\eta))=
\Delta_{f_1,\ldots f_d}(N_{K^+/F^+}\eta)\ \ \ \mbox{for all $\eta\in {\textstyle \bigwedge^d_{\bbQ \bar{G}}\bbQ U_S(K^+)}$}
\eeq
 Hence, if $\eta\in\Lambda_{0,S}(K^+/k)$ then
$\Delta_{f_1,\ldots f_d}(e^d N_{K^+/F^+}\eta)=\pi_{K^+/F^+}(\Delta_{e\hat{f}_1,\ldots e\hat{f}_d}(\eta))$ lies in $\bbZ \bar{G}_F$. Letting the $f_i$ vary,
it follows that $N_{K^+/F^+}\eta$ lies in $e^{-d}\Lambda_{0,S}(F^+/k)$.\ePf
\noindent (The proof shows that $e=1$ if, for instance, $|\Gal(K^+/F^+)|=[K:F]$ is odd.) {\em Suppose now that $\mu_{\pnpo}\subset F$ for some $n\geq 0$} and
that $\fP\in S_p(K)$ lies above $\fp\in S_p(F)$, so we may regard $F_\fp$ as a subfield of
$K_\fP$. Basic properties of the Hilbert symbol show that $(a,b)_{K_\fP,\pnpo}=(a,N_{K_\fP/F_\fp}b)_{F_\fp,\pnpo}$ for all $a\in F_{\fp}^\times$ and $b\in K_{\fP}^\times$.
Regarding $F$ as a subset of $K_p$, it follows easily that
\beq\label{eq: H symbols and the norm}
[\alpha,\beta]_{K,n}=[\alpha, N_{K/F}\beta]_{F,n}\ \ \ \mbox{for all $\alpha\in F^\times$ and $\beta\in K_p^\times$.}
\eeq
\begin{lemma}\label{lemma: second for norming down} Let $\eta\in\Zbp\Lambda_{0,S}(K^+/k)$ and $\theta\in{\textstyle \bigwedge^d_{\bbZ_pG}U^1(K_p)}$. Then $N_{K^+/F^+}\eta$ lies in $\Zbp\Lambda_{0,S}(F^+/k)$ and
$\pi_{K/F}(H_{K/k,n}(\eta,\theta))=H_{F/k,n}(N_{K^+/F^+}\eta, N_{K/F}\theta)$.
\end{lemma}
By $\Zbp$-linearity in $\eta$ and the fact that $p\neq 2$, we may assume $\eta\in e\Lambda_{0,S}(K^+/k)$ with $e$ as in Lemma~\ref{lemma: Norming down Lambda}.
The latter then shows that $N_{K^+/F^+}(\eta)$ lies in $\Lambda_{0,S}(F^+/k)\subset\Zbp\Lambda_{0,S}(F^+/k)$.
Similarly, we may assume that $\theta=u_1\wedge\ldots\wedge u_d$ with $u_i\in U^1(K_p)^-\ \forall\,i$. We let $f_i$
be the map $[\,\cdot\,, u_i]_{K,n}\in\Hom_\bbZ(U_S(K^+),\ZpnpoZ)$ and choose a lift $\tilde{f}_i\in\Hom_\bbZ(U_S(K^+),\bbZ)$ for each $i$.
If $\tilde{g}_i$ denotes the restriction of $\tilde{f}_i$ to
$U_S(F^+)$ then equation~(\ref{eq: H symbols and the norm}) says that $\tilde{g}_i$ lifts
the map $g_i:=[\,\cdot\,, N_{K/F}u_i]_{F,n}\in\Hom_\bbZ(U_S(F^+),\ZpnpoZ)$. Just as for~(\ref{eq: Delta and pi}) we find
$
\pi_{K^+/F^+}(\Delta_{\tilde{f}_1,\ldots\tilde{f}_d}(\eta))=\Delta_{\tilde{g}_1,\ldots \tilde{g}_d}(N_{K^+/F^+}\eta)
$
and since both sides lie in $\bbZ G_F$, we can reduce modulo $\pnpo$ to get
$
\pi_{K^+/F^+}(\tilde{\Delta}_{f_1,\ldots f_d}(\eta))=\tilde{\Delta}_{g_1,\ldots g_d}(N_{K^+/F^+}\eta)
$.
We conclude by applying $2^d\bar{\kappa}_{F,n}^\ast$ to both sides and using $\bar{\kappa}_{F,n}^\ast\circ\pi_{K^+/F^+}=\pi_{K/F}\circ\bar{\kappa}_{K,n}^\ast$.\ePf
\begin{prop}\label{prop: CC under descending K}
Suppose $K\supset F\supset k$ as above and $N_{K/F}:{\textstyle \bigwedge^d_{\bbZ_pG}U^1(K_p)}^-\rightarrow{\textstyle \bigwedge^d_{\bbZ_pG}U^1(F_p)}^-$ is surjective.
If $F\supset\mu_{\pnpo}$ for some $n\geq 0$  then $CC(K/k,S,p,n)$ implies $CC(F/k,S,p,n)$.
\end{prop}
\bPf\
We assume that $CC(K/k,S,p,n)$ holds, so also $IC(F/k,S,p)$ holds and $RSC(F^+/k,S;\Zbp)$ holds with solution $\eta_{K^+/k,S}$, say. Prop.~\ref{prop: IC under descending K} implies
$IC(F/k,S,p)$. Moreover, Prop.~\ref{prop: RSC(Q) under descending K} and Lemma~\ref{lemma: second for norming down} imply
$RSC(F^+/k,S;\Zbp)$  and that for any $\theta\in{\textstyle \bigwedge^d_{\bbZ_pG}U^1(K_p)}^-$ we have
\begin{eqnarray*}
\kappa_n(\tau_1\ldots\tau_d)H_{F/k,n}(\eta_{F^+/k, S}, N_{K/F}\theta)&=&
\kappa_n(\tau_1\ldots\tau_d)H_{F/k,n}(N_{K^+/F^+}\eta_{K^+/k, S}, N_{K/F}\theta)\\
&=&\pi_{K/F}(\kappa_n(\tau_1\ldots\tau_d)H_{K/k,n}(\eta_{K^+/k, S},\theta))\\
&=&\pi_{K/F}(\overline{\fs_{K/k,S}(\theta)})\\
&=&\overline{\fs_{F/k,S}(N_{K/F}\theta)}
\end{eqnarray*}
The result now follows from the surjectivity condition.\ePf
\noindent Finally, if $n\geq n'\geq 0$ then $H_{K,n}(\eta,\theta)\equiv H_{K,n'}(\eta,\theta)\pmod{p^{n'+1}}$ for all $\eta\in\Zbp\Lambda_{0,S}(K^+/k)$ and
$\theta\in{\textstyle \bigwedge^d_{\bbZ_pG}U^1(K_p)}^-$. (The proof is an exercise using the definitions of the Hilbert symbol, $[\cdot,\cdot]_{K,n}$,
$\tilde{\Delta}$, $H_{k,n}$, $\kappa_n$ \etc\ and the fact that $\zeta_n^{p^{n-n'}}=\zeta_{n'}$!). One deduces easily
\begin{prop}\label{prop: CC under decreasing n}
If $K\supset\mu_{\pnpo}$ for some $n\geq 0$  then $CC(K/k,S,p,n)$ implies $CC(K/k,S,p,n')$ for all $n'$ with $n\geq n'\geq 0$. \ePf\end{prop}
\section{The Case $p\ndiv |G|$}\label{sec: pndivG}
Let $\cX_{\bbQ_p}$ denote the set of irreducible $\bbQ_p$-valued characters of $G$ which is
in natural bijection with $\Gal(\barbbQ_p/\bbQ_p)$-conjugacy
classes of absolutely irreducible characters $\phi\in\Hom(G,\bar{\bbQ}_p^\times)$. (Precisely, if  $\Phi$ lies in $\cX_{\bbQ_p}$ then
its idempotent $e_\Phi\in\bbQ_p G$ splits in $\barbbQ_p G$ as the sum of the idempotents $e_\phi$ where $\phi$ runs once through the conjugacy
class corresponding to $\Phi$). We shall say that the characters $\phi$ in this conjugacy class {\em belong} to $\Phi$ and we shall call $\Phi$ {\em odd} iff
one -- hence any -- such $\phi$ is odd (\ie\ $\phi(c)=-1$).
Henceforth we set $\fa:=\bbZ_p G$ and $\fa_\Phi:=e_\Phi\bbZ_p G$. Any $\phi$ belonging to $\Phi$ extends
$\bbQ_p$-linearly to to a homomorphism $\bbQ_p G\rightarrow F_\phi:=\bbQ_p(\phi)$ which in turn restricts to {\em iso}morphisms
from $e_\Phi\bbQ_p G$ to $ F_\phi$ and from $\fa_\Phi$ to $\cO_\phi:=\bbZ_p[\phi]$, the ring of valuation integers of $F_\phi$.
In particular, $\fa_\Phi$ is a complete d.v.r., hence a p.i.d.

For the rest of this section we suppose that {\em the prime $p$ does not divide $|G|$}. This means that the idempotent
$e_\Phi$ lies in $\bbZ_p G$
for each $\Phi\in \cX_{\bbQ_p}$ so that $\fa$ is a product $\prod_{\Phi\in\cX_{\bbQ_p}}\fa_\Phi$.
Any $\fa$-module $M$ splits as a corresponding direct sum $\bigoplus_{\Phi\in\cX_{\bbQ_p}}M_\Phi$,
where $M_\Phi$ denotes the $\fa_\Phi$-module $e_\Phi M$, and $M\mapsto M_\Phi$ is an exact functor.
Since any $\phi$ belonging to $\Phi$ has order prime to $p$, a uniformiser of
$\cO_\phi$ -- hence of $\fa_\Phi$ -- is given by $p$. The $\fa_\Phi$-{\em order ideal} $\oridPhi{N}$ of any finite (=finite length) $\fa_\Phi$-module $N$
is therefore $p^l\fa_\Phi$ where $l$ is the length of any $\fa_\Phi$-composition series for $N$. We shall assume the usual properties
of the order ideal, such as multiplicativity in exact sequences.  Each $p$-adic-valued character
$\phi\in\Hom(G,\bar{\bbQ}_p^\times)$  corresponds to a unique complex character $\chi\in\hat{G}$ such that $\phi=j\circ\chi$
where $j$ is the fixed embedding $\barbbQ\rightarrow\barbbQ_p$. We write
$\hat{\chi}$ and $\hat{\phi}=j\circ\hat{\chi}$ respectively for the associated complex and $p$-adic {\em primitive} ray-class characters, $\ff_\phi$ for $\ff_\chi$ and
$K^\phi$ for the field $K^{\ker(\phi)}=K^{\ker(\chi)}$ cut out by $\phi$, so that $\chi$ and $\phi$ factor through $G_\phi:=\Gal(K^\phi/k)$.
Work of Siegel~\cite{Siegel} and Klingen (see also
Shintani~\cite[Cor. to Thm. 1]{Shintani}) implies that $\Theta_{K^\phi/k,S_\infty}(0)$ lies in $\bbQ G_\phi$ so that $L(0,\hat{\chi}\inv)=\chi(\Theta_{K^\phi/k,S^0(K^\phi/k)}(0))$
lies in $\bbQ(\chi)$. Thus
$j(L(0,\hat{\chi}\inv))=\phi(\Theta_{K^\phi/k,S^0(K^\phi/k)}(0))$ lies in $F_\phi$ and is independent of $j$ so, by a slight abuse of notation, we write it simply
as $L(0,\hat{\phi}\inv)$.
\begin{thm}\label{thm:chi(SKkS) when p ndiv G}
If $p\ndiv|G|$ then,  for any odd $\Phi\in\cX_{\bbQ_p}$ and $\phi\in\Hom(G,\barbbQ_p^\times)$ belonging to $\Phi$, we have
\beq\label{eq: of thm pndivG}
\phi(\SKkS)=\phi(\oridPhi{(U^1(K_p)_{tor})_\Phi})\prod_{\fq\in S\setminus S_\infty\atop\fq\ndiv p\ff_\phi}
\left(1-N\fq\inv\hat{\phi}([\fq])\right)L(0,\hat{\phi}\inv)
\eeq
(an equality of
fractional ideals of $F_\phi$) where $L(0,\hat{\phi}\inv)$ is as defined above.
\end{thm}
Equation~(\ref{eq: of thm pndivG}) for each $\Phi$ clearly determines $\SKkS$. Before giving the proof, we reformulate it and deduce some consequences.
Firstly, $U^1(K_p)_{tor}$ is nothing but
$\mu_{p^\infty}(K_p)=\prod_{\fP|p}\mu_{p^\infty}(K_\fP)$. Next, for given $\phi$ as above we define a $\bbZ_p G_\phi$ submodule of $\bbQ_p G_\phi$ by
\[
J_\phi:={\rm ann}_{\bbZ_p G_\phi}(\mu_{p^\infty}(K^\phi))\Theta_{K^\phi/k,S^0(K^\phi/k)}(0)
\]
Since $p\ndiv [K:K^\phi]$, we have
\begin{eqnarray*}
\phi(\nu_{K/K^\phi}(J_\phi))&=&\phi({\rm ann}_\fa(\mu_{p^\infty}(K^\phi))\nu_{K/K^\phi}(\Theta_{K^\phi/k,S^0(K^\phi/k)}(0)))\\
  &=&\phi({\rm ann}_\fa(\mu_{p^\infty}(K^\phi)))[K:K^\phi]L(0,\hat{\phi}\inv)\\
  &=&\phi({\rm ann}_{\fa_\Phi}(\mu_{p^\infty}(K^\phi)_\Phi))L(0,\hat{\phi}\inv)\\
  &=&\phi({\rm ann}_{\fa_\Phi}(\mu_{p^\infty}(K)_\Phi))L(0,\hat{\phi}\inv)\\
  &=&\phi([\mu_{p^\infty}(K)_\Phi]_{\fa_\Phi})L(0,\hat{\phi}\inv)
\end{eqnarray*}
(the last equation because $\mu_{p^\infty}(K^\phi)_\Phi$ is cyclic
over $\bbZ$, so over $\fa_\Phi$). Thus we may reformulate~(\ref{eq:
of thm pndivG}) as
\beq\label{eq: reformulation of thm pndivG}
\phi(\SKkS)=\phi(\oridPhi{(\mu_{p^\infty}(K_p)/\mu_{p^\infty}(K))_\Phi})\prod_{\fq\in
S\setminus S_\infty\atop\fq\ndiv p\ff_\phi}
\left(1-N\fq\inv\hat{\phi}([\fq])\right)\phi(\nu_{K/K^\phi}(J_\phi))
\eeq But $J_\phi$ is spanned over $\bbZ_p$ by ${\rm ann}_{\bbZ
G_\phi}(\mu(K^\phi))\Theta_{K^\phi/k,S^0(K^\phi/k)}(0)$ which lies
in $\bbZ G_\phi$ by the well-known result of Deligne-Ribet and
(independently) Pi.~Cassou-Nogu\`es (see Th\'eor\`eme~6.1
of~\cite[p.~107]{Tate}).
Hence $J_\phi\subset \bbZ_p G_\phi$ and  so~(\ref{eq: reformulation
of thm pndivG}) implies that $\phi(\SKkS)\subset \cO_\phi$ for all
odd $\phi\in\Hom(G,\barbbQ_p^\times)$. Consequently,
\begin{cor}\label{cor: IC when pndivG} If $p\ndiv |G|$ then $IC(K/k,S,p)$ holds.
\end{cor}
\rem\label{rem: Jones' formula from ours}\ We explain the relation between Jones'
formula~(\ref{eq: Jones on ETNC and pndivG}) and our Theorem~\ref{thm:chi(SKkS) when p ndiv G}, recast as equation~(\ref{eq: reformulation of thm pndivG}) for all odd $\phi$:
The ray-class group $\Cl_\fm(K)$ appearing in~(\ref{eq: Jones on ETNC and pndivG}) fits into an
exact sequence of $\bbZ G$-modules:
\[
0\rightarrow\overline{\cO_K^\times}\longrightarrow\prod_{\fq\in S^1\setminus S_\infty}\prod_{\fQ|\fq}(\cO_K/\fQ)^\times
\longrightarrow\Cl_\fm(K)\longrightarrow\Cl(K)\rightarrow 0
\]
(where the first non-zero term is simply the image of $\cO_K^\times$ in the second). Now tensor this sequence with $\bbZp$ and take minus parts.
Using the fact $(\cO_K^\times\otimes\bbZ_p)^-=\mu_{p^\infty}(K)$ and
isomorphisms similar to~(\ref{eq: galois iso for res fld}) one finds with a little
work that~(\ref{eq: Jones on ETNC and pndivG}) is equivalent to the following for each odd $\phi$ as in
Theorem~\ref{thm:chi(SKkS) when p ndiv G}.
\[
\phi(\fS_{K/k,S^1})=\phi(\oridPhi{(\mu_{p^\infty}(K_p)/\mu_{p^\infty}(K))_\Phi})\prod_{\fq\in
S^1\setminus S_\infty\atop\fq\ndiv p\ff_\phi}
\left(1-N\fq\inv\hat{\phi}([\fq])\right)\phi((\oridPhi{\Cl(K)\otimes\bbZ_p})_\Phi)
\]
Since $p\ndiv[K:K^\phi]$, one sees that this in turn is equivalent to our~(\ref{eq: reformulation of thm pndivG}) (with $S=S^1$)
{\em if and only if} $\phi(J_\phi)=\phi((\oridPhi{\Cl(K^\phi)\otimes\bbZ_p})_\Phi)$
(where $\Phi$ and $\phi$ are now considered as odd characters of $G_\phi$). But Theorem~3 of~\cite{W=[46] of Jones} establishes the latter equality subject to a rather mild
condition (`$S_{\phi,p}=0$') on the character $\phi$.\vertsp\\
\noindent \textsc{Proof of Theorem~\ref{thm:chi(SKkS) when p ndiv G}}\
For each $i=1,\ldots,d$, we write $\fp_i$ for $\fP_i\cap k$ (namely the prime ideal in $S_p(k)$ which is defined
by the embedding $j\tau_i:\barbbQ\rightarrow\barbbQ_p$). The map $\{1,\ldots,d\}\rightarrow S_p(k)$ sending $i$ to $\fp_i$ is clearly surjective
so for any $\fp\in S_p(k)$ we write $I(\fp)$ for its fibre over $\fp$ and choose an element $i(\fp)\in I(\fp)$. Thus $\fP_{i(\fp)}\cap k=\fp_{i(\fp)}=\fp,\ \forall\,\fp\in S_p(k)$
and the extension $K_{\fP_{i(\fp)}}/k_\fp$ is Galois with group
$D_\fp(K/k)$ of order prime to $p$. It follows
(\eg\ by a theorem of E.~Noether, since $K_{\fP_{i(\fp)}}/k_\fp$ is tame) that we may choose
an element $b_\fp\in\cO_{K_{\fP_{i(\fp)}}}$ freely generating $\cO_{K_{\fP_{i(\fp)}}}$ over $\cO_{k_\fp}D_\fp(K/k)$.
Let $b$ be the element of $\cO_{K_p}:=\prod_{\fP\in S_p(K)}\cO_{K_\fP}$
whose component in $\cO_{K_\fP}$ is $b_\fp$ whenever $\fP=\fP_{i(\fp)}$ for some $\fp\in S_p(k)$ and is $0$ otherwise.
Then $b$ is a free generator for $\cO_{K_p}$ over $\cO_{k_p}G$, where
$\cO_{k_p}$ denotes the ring $\prod_{\fp\in S_p(k)}\cO_{k_\fp}$ which we identify
with $\cO_k\otimes_\bbZ\bbZ_p$. So if $\listsub{c}{d}$ is a $\bbZ$-basis of $\cO_k$ then $c_1\otimes 1,\ldots,c_d\otimes 1$ is a
$\bbZ_p$-basis of $\cO_{k_p}$ and
$a_1:=b(c_1\otimes 1),\ldots, a_d:=b(c_d\otimes 1)$ is a free basis for  $\cO_{K_p}$ over $\bbZ_p G=\fa$.

For any $\fP\in S_p(K)$ let $\hat{\fP}$ and $e_\fP$ denote
respectively the maximal ideal and the ramification index of
$K_\fP/\bbQ_p$. Clearly, $e_\fP$ depends only on $\fp$, the prime lying
below $\fP$ in $K$. The exponential series converges on
$p\cO_{K_\fP}=\hat{\fP}^{e_\fP}$ for each $\fP\in S_p(K)$ and defines
a $\bbZ_p D_\fp(K/k)$-isomorphism to $U^{e_\fP}(K_\fP)$. To shorten
notation, we write $U^1$ for $U^1(K_p)$ and  $U^{\underline{e}}$ for
$\prod_{\fP\in S_p(K)}U^{e_\fP}(K_\fP)\subset U^1$. It follows from
the above that the map ${\rm Exp}_p=\left(\prod_{\fP\in
S_p(K)}\exp_p\right):p\cO_{K_p}\rightarrow U^{\underline{e}}$ is an
$\fa$-isomorphism and hence that $U^{\underline{e}}$ is free over
$\fa$ with basis $w_1:={\rm Exp}_p(pa_1),\ldots,w_d:={\rm
Exp}_p(pa_d)$. It is also of finite index in $U^1$ and since $\fa$
is a product of the p.i.d.'s $\fa_\Phi$, it follows that
$U^1/U^1_{tor}$ must also be $\fa$-free of rank $d$, so, in an
additive notation, we get \beq\label{eq: basis for Uonemod tor}
U^1=U^1_{tor}\oplus\bigoplus_{i=1}^d\fa u_i\ \ \ \mbox{where
$\listsub{\bar{u}}{d}$ is any free $\fa$ basis of $U^1/U^1_{tor}$}
\eeq Now let $\phi$ and $\Phi$ be as in the statement of the Theorem
and let $M\in M_d(\fa_\Phi)$ be the matrix representing
$\listsub{e_\Phi\bar{w}}{d}$ in terms of the $\fa_\Phi$-basis
$\listsub{e_\Phi\bar{u}}{d}$ of $(U^1/U^1_{tor})_\Phi$. The determinant of $M$ has two different interpretations.
On the one hand,
if we write $\overline{U^{\underline{e}}}$ for the isomorphic image
of $U^{\underline{e}}$ in $U^1/U^1_{tor}$ then the general theory of
p.i.d.'s and order ideals gives
\[
\det(M)\fa_\Phi=\oridPhi{(U^1/U^1_{tor})_\Phi/(\overline{U^{\underline{e}}})_\Phi}=
\oridPhi{U^1_\Phi/U^{\underline{e}}_\Phi}\oridPhi{(U^1_{tor})_\Phi}\inv
\]
Now, for each $\fp\in S_p(k)$, $\fP\in S_p(K)$ above $\fp$ and $l\geq 1$, there is a well-known $\bbZ_p D_\fp(K/k)$-isomorphism
$U^l(K_\fP)/U^{l+1}(K_\fP)\rightarrow \hat{\fP}^l/\hat{\fP}^{l+1}$ induced by $x\mapsto x-1$. This gives an $\fa$-isomorphism after taking products of both sides over the
$\fP$ above $\fp$. Applying $e_\Phi$ and letting $\fp$ and $l$ vary,
a simple argument with exact sequences shows that $U^1_\Phi/U^{\underline{e}}_\Phi$  has the same $\fa_\Phi$-order ideal as $\fM_\Phi/(p\cO_{K_p})_\Phi$  where
$\fM$ denotes $\prod_{\fP\in S_p(K)}\hat{\fP}\subset K_p$.
Therefore
\begin{eqnarray}
\det(M)\fa_\Phi&=&
\oridPhi{\fM_\Phi/(p\cO_{K_p})_\Phi}\oridPhi{(U^1_{tor})_\Phi}\inv\nonumber\\
&=&
\oridPhi{(\cO_{K_p})_\Phi/(p\cO_{K_p})_\Phi}\oridPhi{(\cO_{K_p})_\Phi/\fM_\Phi}\inv\oridPhi{(U^1_{tor})_\Phi}\inv\nonumber\\
&=&
p^d\oridPhi{(\cO_{K_p})_\Phi/\fM_\Phi}\inv\oridPhi{(U^1_{tor})_\Phi}\inv\label{eq:first eq with det M}
\end{eqnarray}
since $\cO_{K_p}$ is free of rank $d$ over $\fa$. On the other hand, Proposition~\ref{prop:ker and im of sKks},
Equation~(\ref{eq: basis for Uonemod tor}) and the definition of
$\sKkS$ give
\begin{eqnarray}
\phi(\det(M))\phi(\SKkS)&=&\phi(\det(M))\phi(e_\Phi\sKkS(\wedgesub{u}{d})\fa)\nonumber\\
  &=&\phi(\det(M)\sKkS(\wedgesub{e_\Phi u}{d}))\cO_\phi\nonumber\\
  &=&\phi(\sKkS(\wedgesub{e_\Phi w}{d}))\cO_\phi\nonumber\\
  &=&j(\chi(a^{-,\ast}_{K/k,S}))\phi(R^{(j)}_{K/k,p}(\wedgesub{w}{d}))\cO_\phi\label{eq:second eq with det M}
\end{eqnarray}
where $\phi=j\circ \chi$. But tracing through the definitions we have
\[
R_{K/k,p}^{(j)}(\wedgesub{w}{d})=
\det\left(\sum_{g\in G}\log_p(\delta_i^{(j)}(g\inv{\rm Exp}_p(pa_t)))g\right)_{i,t=1}^d
\]
and
\begin{eqnarray*}
\log_p(\delta_i^{(j)}(g\inv{\rm Exp}_p(pa_t)))&=&\log_p(j\tau_i\circ\iota_{\fP_i}{\rm Exp}_p(g\inv pa_t))\\
  &=&\log_p(j\tau_i\exp_p(\iota_{\fP_i}g\inv pa_t))\\
  &=&\delta_i^{(j)}(g\inv pa_t)\\
  &=&p\delta_i^{(j)}(g\inv b)j\tau_i(c_t)\\
\end{eqnarray*}
so that
\beq\label{eq:R in terms of prod of delta b's}
R_{K/k,p}^{(j)}(\wedgesub{w}{d})=p^d\prod_{i=1}^d\left(\sum_{g\in G}\delta_i^{(j)}(g\inv b)g\right)\det(j\tau_i(c_t))_{i,t=1}^d
=\pm p^dj(\sqrt{d_k})\prod_{i=1}^d\delta_i^{(j),G}(b)
\eeq
Applying $\phi$ to  Equations~(\ref{eq:first eq with det M}) and (\ref{eq:R in terms of prod of delta b's})
and combining them with~(\ref{eq:second eq with det M}) gives
\begin{eqnarray}
\phi(\SKkS)&=&\phi\left(\oridPhi{(\cO_{K_p})_\Phi/\fM_\Phi}\right)\phi\left(\oridPhi{(U^1_{tor})_\Phi}\right)j(\sqrt{d_k}\chi(a^{-,\ast}_{K/k,S}))\prod_{i=1}^d
\phi(\delta_i^{(j),G}(b))\label{eq:next stage 1}
\end{eqnarray}
Now fix $\fp\in S_p(k)$ and write $D_\fp$ for $D_\fp(K/k)$ and $T_\fp$
for $T_\fp(K/k)$. Considering $\prod_{\fP|\fp}(\cO_{K_\fP}/\hat{\fP})$ as an $\fa$-submodule of $\cO_{K_p}/\fM$,
we have natural $\fa$-isomorphisms:
\[
\prod_{\fP|\fp}(\cO_{K_\fP}/\hat{\fP})\cong
\fa \otimes_{\bbZ_p D_\fp}(\cO_{K_{\fP_{i(\fp)}}}/\hat{\fP}_{i(\fp)})\cong
\fa \otimes_{\bbZ_p D_\fp}(\bbZ_p D_\fp\otimes_{\bbZ_p T_\fp}(\cO_k/\fp))\cong
\fa \otimes_{\bbZ_p T_\fp}(\cO_k/\fp)
\]
(where the action  on $\cO_k/\fp$ is trivial and the second isomorphism is from the normal basis theorem in the {\em residue field} extension of  $K_{\fP_{i(\fp)}}/k_\fp$). It follows easily that
$(\prod_{\fP|\fp}(\cO_{K_\fP}/\hat{\fP}))_\Phi$ is trivial unless $T_\fp\subset \ker(\phi)$ (\ie\ $\fp\ndiv \ff_\phi$) in which case it has order ideal $(N\fp) \fa_\Phi$. Taking the product over
all $\fp\in S_p(k)$, it follows that
\beq\label{eq:next stage 2}
\phi\left(\oridPhi{(\cO_{K_p})_\Phi/\fM_\Phi}\right)=\left(\prod_{\fp\in S_p(k)\atop\fp\ndiv\ff_\phi}N\fp\right)\cO_\phi
\eeq
Furthermore, equations~(\ref{eq:defn_aKkS}),~(\ref{eq:thetaStozetaandL}) and~(\ref{eq:Lfns}) give:
\begin{eqnarray}
\sqrt{d_k}\chi(a^{-,\ast}_{K/k,S})&=&\prod_{\fq\in S\setminus S_\infty \atop \fq\ndiv \ff_\chi}\left(1-N\fq^{-1}\hat{\chi}([\fq])\right)\sqrt{d_k}(i/\pi)^d L(1,\hat{\chi})\nonumber\\
  &=&\prod_{\fq\in S\setminus S_\infty \atop \fq\ndiv \ff_\chi}\left(1-N\fq^{-1}\hat{\chi}([\fq])\right)(-1)^d\tau(\chi)\inv L(0,\hat{\chi}\inv)\label{eq:next stage 3}
\end{eqnarray}
The second equality follows from Hecke's functional equation for the $L$-function. (To be perfectly precise, we are using the version stated on p.~36
of~\cite{Frohlich}, taking $s=0$  and Fr\"ohlich's complex character `$\bar{\theta}$' on ${\rm Id}(k)$ -- the id\`ele group of $k$ -- to be the one
obtained by composing $\chi$ with the map ${\rm Id}(k)\rightarrow G$ coming from class-field theory.)
Applying $j$ to and~(\ref{eq:next stage 3}) and combining with Equations~(\ref{eq:next stage 1}) and (\ref{eq:next stage 2}) gives
\begin{eqnarray*}
\phi(\SKkS)&=&\phi\left(\oridPhi{(U^1_{tor})_\Phi}\right)\prod_{\fp\in S_p(k)\atop\fp\ndiv\ff_\phi}N\fp
\prod_{\fq\in S\setminus S_\infty \atop \fq\ndiv \ff_\phi}\left(1-N\fq^{-1}\hat{\phi}([\fq])\right)L(0,\hat{\phi}\inv)j(\tau(\chi))\inv\prod_{i=1}^d\phi(\delta_i^{(j),G}(b))\\
  &=&\phi\left(\oridPhi{(U^1_{tor})_\Phi}\right)
\prod_{\fq\in S\setminus S_\infty \atop \fq\ndiv p\ff_\phi}\left(1-N\fq^{-1}\hat{\phi}([\fq])\right)L(0,\hat{\phi}\inv)j(\tau(\chi))\inv\prod_{i=1}^d\phi(\delta_i^{(j),G}(b))
\end{eqnarray*}
where we have used the facts that every prime ideal $\fp $ in $S_p(k)$ is contained in $S$ and that if, in addition,
it does not divide $\ff_\phi$ then $N\fp(1-N\fp^{-1}\hat{\phi}([\fq]))=(N\fp-\hat{\phi}([\fq]))$
lies in $\cO_\phi^\times$. The argument so far  shows that $j(\tau(\chi))\inv\prod_{i=1}^d\phi(\delta_i^{(j),G}(b))$ lies in $F_\phi$.
The theorem will follow if we can prove that it too lies in $\cO_\phi^\times$, \ie\ that
\[
j(\tau(\chi))\sim\prod_{i=1}^d\phi(\delta_i^{(j),G}(b))
\]
where `$a\sim b$' means that $a,b\in\barbbQ_p^\times$ have the same $p$-adic absolute value.
Recall that Fr\"{o}hlich defines $\tau(\chi)$ as the product $\prod_{\fq\nin S_\infty }\tau(\chi_\fq)$ where $\chi_\fq:k_\fq^\times\rightarrow\barbbQ^\times$
is the $\fq$-component of the complex
id\`ele character associated to $\chi$ and $\tau(\chi_\fq)$ is the `local Gauss sum' (which equals $1$ unless
the $\fq|\ff_\chi$, so the product is finite). For definitions and basic properties of the algebraic integers $\tau(\chi_\fq)$ see see~\cite[p. 34-35]{Frohlich}
or~\cite[II-\S 2]{Martinet}. In particular, Eq.~(5.7) on~\cite[p. 34]{Frohlich} shows that $j(\tau(\chi))\sim 1$ unless $\fq\in S_p(k)$.
Hence $j(\tau(\chi))\sim\prod_{\fp\in S_p(k)}j(\tau(\chi_\fp))$ and since $\{1,\ldots,d\}$ is the disjoint union  $\bigcup_{\fp\in S_p(k)}I(\fp)$ it suffices to show that for any
$\fp$ in $S_p(k)$
\beq\label{eq: pf of pndivG reduces to this}
j(\tau(\chi_\fp))\sim\prod_{i\in I(\fp)}\phi(\delta_i^{(j),G}(b))
\eeq
But this is essentially (a special case of) Theorem~23 of~\cite{Frohlich}: Take $F:=\overline{j\tau_{i(\fp)}(k)}$, $L:=\overline{j\tau_{i(\fp)}(K)}$ as subfields of
$\barbbQ_p$,  isomorphic via
$j\tau_{i(\fp)}$ to $k_\fp$ and $K_{\fP_{i(\fp)}}$ respectively. The extension
$L/F$ is thus abelian  with Galois group $\Gamma$ which we identify via $j\tau_{i(\fp)}$ with $D_\fp$.  We take
Fr\"ohlich's  character `$\chi$' to be our $\chi_\fp:k_\fp^\times\rightarrow\barbbQ^\times$ which factors through the local reciprocity map $k_\fp^\times\rightarrow D_\fp$
and so may also be regarded as $\chi$ restricted to  $D_\fp=\Gamma$. Thus Fr\"ohlich's `$\chi^j$' may similarly be identified with our $\phi$ restricted to  $\Gamma$.
Since $\Gamma$ of order prime to $p$,
$L/F$ is {\em tame} so Theorem~23 applies to give (with these identifications)
\[
j(\tau(\chi_\fp))\sim\cN_{F/\bbQ_p}(j\tau_{i(\fp)}(b_\fp)|\phi)
\]
where the R.H.S. is the {\em norm resolvent} (see below) associated to the free generator $j\tau_{i(\fp)}(b_\fp)$ of $\cO_L$ over $\cO_F\Gamma$.
Thus~(\ref{eq: pf of pndivG reduces to this}) and hence our Theorem will follow from
\beq\label{eq: claim}
\prod_{i\in I(\fp)}\phi(\delta_i^{(j),G}(b))\sim\cN_{F/\bbQ_p}(j\tau_{i(\fp)}(b_\fp)|\phi)
\eeq
The proof of~(\ref{eq: claim}) is largely a matter of unravelling our definitions and comparing with Fr\"ohlich's, so we only sketch it. For any $i\in I(\fp)$  we can
choose $g_i\in G$ such that $g_i\fP_i=\fP_{i(\fp)}$ and then
$\sigma_i\in \Gal(\barbbQ_p/\bbQ_p)$ such that $\sigma_i j\tau_{i(\fp)}(x)=j\tau_i g_i\inv(x)$ for any $x\in K_{\fP_{i(\fp)}}$. Then
\begin{eqnarray*}
\phi(\delta_i^{(j),G}(b))&=&\sum_{g\in G}j\tau_i\iota_{\fP_i}(g\inv b)\phi(g)\\
 &=&\sum_{h\in D_\fp}j\tau_ig_i\inv h\inv (b_\fp)\phi(hg_i)\\
 &=& \phi(g_i)\sigma_i\left(\sum_{\gamma\in \Gamma}\gamma\inv (j\tau_{i(\fp)}(b_\fp)\sigma_i\inv(\phi(\gamma))\right)\\
 &\sim&\sigma_i\left(j\tau_{i(\fp)}(b_\fp)|\sigma_i\inv\circ \phi\right)
\end{eqnarray*}
where $\left(j\tau_{i(\fp)}(b_\fp)|\sigma_i\inv\circ \phi\right)$ denotes the {\em resolvent} defined for example in~\cite[Eq.~(4.4), p. 29]{Frohlich}.
Equation~(\ref{eq: claim}) now follows on taking the product over $i\in I(\fp)$, using the definition of the norm resolvent in~\cite[Eq.~(1.4), p. 107]{Frohlich} and the fact (which the reader can easily check) that
as $i$ runs through $I(\fp)$, so $\sigma_i$ runs once through a set of
left coset representatives for $\Gal(\barbbQ_p/F)$ in $\Gal(\barbbQ_p/\bbQ_p)$. (Fr\"ohlich uses right cosets because of his exponential notation for Galois action). This completes the proof of
Theorem~\ref{thm:chi(SKkS) when p ndiv G}. \ePf
Some of the facts used in the above proof will also be useful in the\vertsp\\
\noindent {\sc Proof of Proposition~\ref{prop: trivial CC 2}}\
Since $p\ndiv |G|$, we can use equation~(\ref{eq: basis for Uonemod tor}) to show that any $\theta\in\bigwedge^d_{\bbZ_p G}U^1$ may be expressed as the sum of $xu_1\wedge\ldots\wedge u_d$
(for some $x\in\fa$) and finitely many terms of form $z\wedge v_2\wedge\ldots\wedge v_d$ with $z\in U^1_{\rm tor}$ and $v_i\in U^1$ for $i=2,\ldots,d$.
Since we are assuming that $\theta$ is $\bbZ_p$-torsion, so also is its image $x(\bar{u}_1\wedge\ldots\wedge \bar{u}_d)$ in
$\bigwedge^d_{\bbZ_p G}(U^1/U^1_{\rm tor})$ and since $\bar{u}_1\wedge\ldots\wedge \bar{u}_d$ freely generates the latter over $\fa$, it follows that $x=0$. Thus, by linearity,
we may assume that $\theta=z\wedge v_2\wedge\ldots\wedge v_d$.
On the other hand, $p\ndiv |G|$ also implies $\Zbp\Lambda_{0,S}=\Zbp\overline{\bigwedge^d_{\bbZ\bar{G}}U_S(K^+)}$ by Prop.~\ref{prop:index of barLdU in Lambda}.
If we write $\te$ for $|\bar{G}|e_{S,d,\bar{G}}\in \bbZ \bar{G}$ then it
follows from the eigenspace condition on $\eta$ that it equals $|G|^{-d}(2\te)^d\eta$ and so may be written as a $\Zbp$-linear combination of terms
of form $(1\otimes \te\eps_1)\wedge\ldots\wedge (1\otimes \te\eps_d)$ with $\eps_i\in U_S(K^+)^2\ \forall\,i$. By $\Zbp$-linearity and equation~(\ref{eq:simple formula for H_n}), it therefore
suffices
to show that $[\te\eps,z]_{K,n}=0$ for any $\eps\in U_S(K^+)^2$ and any $z\in U^1_{\rm tor}=\mu_{p^\infty}(K_p)$, say
$z=(z_{\fP})_{\fP}$ with $z_{\fP}\in\mu_{p^\infty}(K_\fP)$ for each $\fP\in S_p(K)$.
By the definitions of $[\cdot,\cdot]_{K,n}$, $[\cdot,\cdot]_{\fP,n}$ and $(\cdot,\cdot)_{K_\fP,\pnpo}$ this reduces further
to the statement that $\sigma_{\iota_\fP(\te\eps), K_\fP}(\zeta_\fP)=\zeta_\fP$ for each $\fP$, where $\zeta_\fP:=z_\fP^{1/\pnpo}$ is a $p$-power root of unity in $(K_\fP)^{\rm ab}$.
But $\zeta_\fP$ actually lies in $\bbQ_p^{\rm ab}$, so local class field theory tells us that $\sigma_{\iota_\fP(\te\eps), K_\fP}(\zeta_\fP)=\sigma_{a_\fP,\bbQ_p}(\zeta_\fP)$ where
$a_\fP:=N_{K_\fP/\bbQ_p}\iota_\fP(\te\eps)=N_{k_\fp/\bbQ_p}\iota_\fP(N_{D_\fp(K/k)}\te\eps)$ and $\fp\in S_p(k)$ lies below $\fP$.
But the image of $N_{D_\fp(K/k)}$ in $\bbZ \bar{G}$ is $N_{D_\fp(K^+/k)}$ or $2N_{D_\fp(K^+/k)}$.
If $|S|>d+1$ then, since $\fp$ lies in $S$, the formula~(\ref{eq: formula for eSdGbar})
shows that $N_{D_\fp(K^+/k)}\te=0$ in $\bbZ \bar{G}$, so $a_\fP=1$ for all $\fP$ and the result follows. Finally, if $|S|=d+1$ then we must have
$S=S_\infty(k)\cup S_p(k)=S_\infty(k)\cup\{\fp\}$ and~(\ref{eq: formula for eSdGbar}) now implies $N_{D_\fp(K^+/k)}\te=|D_\fp(K^+/k)|N_{\bar{G}}$. Hence
$a_\fP$ is a power of $N_{k_\fp/\bbQ_p}\iota_\fP(N_{K^+/k}\eps)$ which equals $\iota_\fP(N_{k/\bbQ}N_{K^+/k}\eps)=N_{K^+/\bbQ}\eps$ since
$S_p(k)=\{\fp\}$.
But $\eps\in U_S(K^+)^2$ implies that $N_{K^+/\bbQ}\eps$, hence also $a_\fP$, is a power of $p$ and the result follows from the well known fact that
$\sigma_{p,\bbQ_p}(\zeta_\fP)=\zeta_\fP$
(Indeed, $p=N_{\bbQ_p(\zeta_\fP)/\bbQ_p}(1-\zeta_\fP)$ implies that $\sigma_{p,\bbQ_p}$ restricts to the identity
on $\bbQ_p(\zeta_\fP)$).\ePf
\section{The Case $k=\bbQ$}\label{sec: k equals Q}
The following lemmas will be used in the proof of Theorem~\ref{thm: CC and IC for K over Q}.
Let $p$ be an odd prime and $f$ a positive integer. We write
$f$ as $f'p^{m+1}$ for some $m\geq -1$ and $f'$ prime to $p$.
We shall abbreviate $\Qxi{f}$ to $K_f$, $\Gal(\Qxi{f}/\bbQ)$ to $G_f$ and $\Gal(\Qf^+/\bbQ)$ to $\bar{G}_f$.
For any
$\bar{a}\in (\ZmodZ{f})^\times$ we write $\sigma_a$ for the element of $G_f$ sending
$\xi_f$ to $\xi_f^a$.
\begin{lemma}\label{lemma: formula for theta plus zero and a minus in case kequals Q}
Let $S=\{\infty\}\cup S_f(\bbQ)$ which contains $S^0(K_f/\bbQ)$. Then, with the above notations,
\begin{enumerate}
\item\label{part1: formula for theta plus zero and a minus in case kequals Q}
$
\Theta^{(1)}_{K_f^+/\bbQ,S}(0)=-\half{\displaystyle\sum_{\bar{g}\in \bar{G}_f}}\log|\bar{g}((1-\xi_f)(1-\xi_f^{-1}))|\bar{g}\inv
$ and
\item\label{part2: formula for theta plus zero and a minus in case kequals Q}
$
a^-_{K_f/\bbQ,S}=e^-.\frac{1}{f}{\displaystyle \sum_{g\in G_f}}
g(\xi_f/(1-\xi_f))g\inv
$
\end{enumerate}
\end{lemma}
\bPf\ For part (i), see~\eg~\cite[p. 203]{Stark}.
A rather indirect proof of the equation in~(ii) uses Prop.~1 of~\cite{Shintani} to calculate $\Phi_{K_f/\bbQ}(0)$ as outlined in~\cite[Example 3.1]{twizo} and returns to
$s=1$ with~(\ref{eq:aKkS_intermsof_PhiKko}).
In principle one can also work `$\chi$-by-$\chi$', calculating $\chi(L.H.S.)$ in~(ii) from the usual formula for
for $L(1,\phi)$  when $\phi$ is an odd {\em primitive} Dirichlet character. (See \eg\ in~\cite[Theorem~67 (b)]{F-T}).
However, the imprimitivity of our $\chi$ and presence of a Gauss sum in the formula make the relation to $\chi$(R.H.S.)
surprisingly difficult. We therefore sketch a direct and very elementary proof of~(ii), similar in some respects to
that of~\cite[Theorem~67]{F-T}: Equation~(\ref{eq:thetaStozetaandL}) shows that $\Theta^-_{K_f/\bbQ,S}(1)=\sum_{a=1,(a,f)=1}^{f-1} t_a\sigma_a\inv$ where
$t_a=\half\lim_{s\rightarrow 1}(\zeta_{K_f/\bbQ,S}(s,\sigma_a)-\zeta_{K_f/\bbQ,S}(s,\sigma_{-a}))$ and
$\zeta_{K_f/\bbQ,S}(s,\sigma_a)=\sum_{n\geq 1, n\equiv a\ (f)}n^{-s}$ for ${\rm Re}(s)>1$.
For any $1\leq c\leq f-1$ the function $Z(s,c):=\sum_{n\geq 1}\xi_f^{cn}n^{-s}$ converges
(conditionally) to a continuous function of $s\in(0,\infty]$. For each $1\leq a\leq f-1$ with $(a,f)=1$ and any $s\in(1,\infty]$ we find easily that
\begin{eqnarray}
\zeta_{K_f/\bbQ,S}(s,\sigma_a)-\zeta_{K_f/\bbQ,S}(s,\sigma_{-a})&=&{\frac{1}{f}}\sum_{b=1}^{f-1}(\xi_f^{-ab}-\xi_f^{ab})Z(s,b)\nonumber\\
&=&{\frac{1}{2f}}\sum_{b=1}^{f-1}(\xi_f^{-ab}-\xi_f^{ab})(Z(s,b)-Z(s,f-b))\label{eq: first elementary}
\end{eqnarray}
But if $\log$ denotes the principal branch of logarithm, then Abel's Lemma and some Euclidean geometry show that
$Z(1,c)=-\log(1-\xi_f^c)=-\log|1-\xi_f^c|+i\pi(\half -\frac{c}{f})$. So, letting
$s\rightarrow 1+$ in~(\ref{eq: first elementary}), substituting for $Z(1,b)$, $Z(1,f-b)$ and using the identity
$\sum_{b=1}^{f-1}b\xi_f^{ab}=-\frac{f}{1-\xi_f^a}$ for $(a,f)=1$, we find after rearranging that
$t_a=-\frac{i\pi}{2f}(\frac{\xi_f^a}{1-\xi_f^a}-\frac{\xi_f^{-a}}{1-\xi_f^{-a}})$, which implies~(ii).
\ePf\noindent
Let us write $\hat{K}_f$ for the field $\bbQ_p(\mu_f)\subset\barbbQ_p$.
The proof of Theorem~\ref{thm: CC and IC for K over Q} depends crucially on the following cyclotomic explicit reciprocity law
due to Coleman. (The case $f'=1$ was proved much earlier by
Artin and Hasse in~\cite{AH}).
\begin{lemma}\label{lemma: Coleman's rec. law}
Let $\hat{\xi}_f$ be any primitive $f$th root of unity in $\hat{K}_f$ and let $v\in U^1(\hat{K}_f)$. Then
$ b(\hat{\xi}_f, v):= \frac{1}{f}{\rm
Tr}_{\hat{K}_f/\bbQ_p}((\hat{\xi}_f/(1-\hat{\xi}_f))\log_p(v)) $
lies in $\bbZ_p$. Furthermore
\[
(1-\hat{\xi}_f,v)_{\hat{K}_f,p^{m+1}}=(\hat{\xi}_f^{f'})^{-b(\hat{\xi}_f,
v)}
\]
(The R.H.S. makes sense because $\hat{\xi}_f^{f'}$ is a primitive
$p^{m+1}$th root of unity.)
\end{lemma}
\bPf This follows from Corollary~15 of~\cite{Coleman Dilogarithm}.
We first write $\hat{\xi}_f$ uniquely as
$\hat{\xi}_f=\hat{\zeta}_{p^{m+1}}\hat{\zeta}_{f'}$ where
$\hat{\zeta}_{p^{m+1}}$ and $\hat{\zeta}_{f'}$ are generators
respectively of $\mu_{p^{m+1}}$ and $\mu_{f'}$ in $\hat{K}_f$. We also write $H$ for $\hat{K}_{f'}$, an unramified extension of $\bbQ_p$, and $\cO_H$ for its ring of integers.
Now $1-\hat{\xi}_f=h(u_m)$ where $h(T)$ denotes the linear
polynomial $1-\hat{\zeta}_{f'}(1-T)\in\cO_H[T]$ and
$u_m:=1-\hat{\zeta}_{p^{m+1}}$. The Frobenius element  $\varphi$ of
$\Gal(H/\bbQ_p)$ may be extended to an automorphism of $\cO_H[T]$ (\emph{resp.} of
$H(\mu_{p^\infty})\subset\barbbQ_p$) by acting trivially on $T$
(\emph{resp.} on $\mu_{p^\infty}$)).  Suppose $l\geq 1$ and $l\geq
i\geq 0$. Since $\varphi(\hat{\zeta}_{f'})=\hat{\zeta}_{f'}^p$, one
verifies easily that
\[
\varphi^{l-i}h(1-(1-T)^{p^{l-i}})=\prod_{\hat{\zeta}\in \mu_{p^{l-i}}}(1-\hat{\zeta}_{f'}\hat{\zeta}(1-T))
\]
Substituting $T=u_l=(1-\hat{\zeta}_{p^{l+1}})$ for any generator $\hat{\zeta}_{p^{l+1}}$ of $\mu_{p^{l+1}}$, it is easy to see that the R.H.S.
becomes the norm from $H(\mu_{p^{l+1}})$ to
$H(\mu_{p^{l-i+1}})$ of $h(u_l)$. Thus
$h(T)$ lies in the subgroup of $\cO_H((T))^\times$ denoted ${\mathcal M}^{(l)}$ by Coleman, and this for \emph{any} $l\geq 1$.
Indeed, this follows from the equation at the foot of p.~376 of~\cite{Coleman Dilogarithm} after correcting the misprint `$\phi^{n-i}$'
to read `$\phi^{i-n}$' (which is necessary for consistency with Coleman's equation~(1) on p.~377).
Now we can apply Coleman's Corollary~15, p.~396, after first correcting another obvious misprint: the meaningless `$\lambda(\alpha)$' in the main
equation should
be replaced by $\lambda(1-\alpha)$ ($=-\log_p(\alpha)$). If we take Coleman's `$n$' to be our $m$, his `$u$' to be our $u_m$
(so that his `$H_n$' is our $\hat{K}_f$)
his $\alpha$ to be our $v$ and his $g$ to be our $h$
(so that $\delta h(T)=(1-T)h(T)\inv dh(T)/dT=\hat{\zeta}_{f'}(1-T)/(1-\hat{\zeta}_{f'}(1-T))$) then
the R.H.S. of the main equation in his Corollary~15 equals $-f'b(\hat{\xi}_f,v)$. The Corollary implies that this lies in $\bbZ_p$ and
(taking into account Coleman's definitions of `${\rm Ind}_{u_m}$' and of
`$(x,y)_m$', the latter agreeing with our $(y,x)_{\hat{K}_f,p^{m+1}}$) that
$(1-\hat{\xi}_f,v)_{\hat{K}_f,p^{m+1}}=(1-u_m)^{-f'b(\hat{\xi}_f,v)}$, from which our Lemma follows immediately.\ePf

\noindent {\sc Proof of Theorem~\ref{thm: CC and IC for K over Q}}
Let $K$ be an absolutely abelian CM field and suppose that
$f=f'p^{m+1}$ is the conductor of $K$ \ie\ the smallest positive
integer such that $K\subset K_f$. Then
$\Sram(K/\bbQ)=\Sram(K_f/\bbQ)=S_f(\bbQ)$. Since $p$ is odd,
$\mu_\pnpo\subset K$ implies (\eg\ by ramification) that $n\leq m$.
Therefore, if $m=-1$ then the Congruence Conjecture
doesn't apply and $IC(K/\bbQ,S,p)$ follows from Theorem~\ref{thm:oldcase 1_of_IC}.
So we may assume $m\geq 0$. By Props.~\ref{prop: IC under increasing S} and~\ref{prop: CC under increasing S} we may further assume
that $S=S^1(K/\bbQ)=\{\infty\}\cup S_f(\bbQ)$ (which is also equal
to $S^0(K/\bbQ)$ and to $S^1=S^0(K_f/\bbQ)$).
If $m=0$ then
$K_f/\bbQ$ is tamely ramified at $p$. If $m\geq 1$ then (since
$p\neq 2$) the ramification group $T_p(K_f/\bbQ)=\Gal(K_f/\Qxi{f'})$
has a unique minimal subgroup of order $p$, namely
$\Gal(K_f/\Qxi{f/p})$. This cannot be contained in $\Gal(K_f/K)$ by
minimality of the conductor $f$. Thus, in any case,
$K_f/K$ is at most
tamely ramified above $p$.
So by Remark~\ref{rem: surjectivity} it suffices to prove
$CC(K_f/\bbQ,S^1,p,m)$ and apply Props.~\ref{prop: CC under decreasing n},
and~\ref{prop: CC under descending K}.

We start with $RSC(K_f^+/\bbQ,S^1;\Zbp)$ (see also~\cite[p. 79]{Tate}).
The algebraic integer $(1-\xi_f)(1-\xi_f^{-1})=(1-\xi_f)^{1+c}$ lies in
$K^+_f$ and the norm
relations for cyclotomic numbers (see for example~\cite[Lemma 2.1]{SolGalRel}) show that, for any $q\in
S_f(\bbQ)$, the number $N_{D_\fq(K_f/\bbQ)}(1-\xi_f)$ equals $p$ or
$1$ according as $f'=1$ or $f'\neq 1$. It follows firstly that $(1-\xi_f)^{1+c}$ lies in $U_{S^1}(K_f^+)$
(and even in $U_{S_\infty}(K_f^+)$ for $f'\neq 1$) and secondly, using Proposition~\ref{prop: equivalent `eigenspace' conditions}, that
$\eta_f:=-\half\otimes (1-\xi_f)^{1+c}\in\half\overline{\textstyle \bigwedge_{\bbZ \bar{G}_f}^1 U_{S^1}(K^+_f)}=
\half\overline{U_{S^1}(K^+_f)}$ satisfies the eigenspace condition w.r.t.\ $(S^1,1,\bar{G}_f)$.
Moreover
$R_{K^+_f/\bbQ}=\lambda_{K_f^+/\bbQ,1}$, so {\em taking $\tau_1$ to
be the identity}, Lemma~\ref{lemma: formula for theta plus zero and a minus in case kequals Q}~\ref{part1: formula for theta plus zero and a minus in case kequals Q}
shows that $\eta_f$ is the unique solution $\eta_{K_f^+/\bbQ,S^1}$ of
$RSC(K_f^+/\bbQ,S^1;\bbZ)$. For any $u\in\bigwedge_{\bbZ_p G_f}^1
U^1(K_f)^-=U^1(K_f)^-$, it follows from~(\ref{eq:simple formula for
H_n}) and~(\ref{eq:Gaction2}) that
\begin{eqnarray}
H_{K_f/\bbQ,m}(\eta_f,u)&=&-\bar{2}\inv \sum_{g\in G_f}[(1-\xi_f)^{(1+c)},gu]_{K_f,m}g\inv\nonumber\\
  &=&-\bar{2}\inv \sum_{g\in G_f}[1-\xi_f,gu^{(1-c)}]_{K_f,m}g\inv\nonumber\\
  &=&-\sum_{g\in G_f}[1-\xi_f,gu]_{K_f,m}g\inv\nonumber\\
  &=&\sum_{g\in G_f}\left(\sum_{\fP\in S_p(K_f)}-[1-\xi_f,\iota_\fP(gu)]_{\fP,m}\right)g\inv\label{eq: H eta in case k= Q}
\end{eqnarray}
(where $(1-\xi_f)$ is identified with its natural images $(1-\xi_f)\otimes 1$ and $\iota_\fP(1-\xi_f)$ in $K_p$ and $K_{\fP}$ respectively).
Next we need to calculate the map $\fs_{K_f/\bbQ,S^1}$.
Fix a choice of $j:\barbbQ\rightarrow\barbbQ_p$. It follows easily from the definitions of $R^{(j)}_{K_f/\bbQ,p}$ and $\fs_{K_f/\bbQ,S^1}$ and from
Lemma~\ref{lemma: formula for theta plus zero and a minus in case kequals Q}~\ref{part2: formula for theta plus zero and a minus in case kequals Q} that
for any $u\in U^1(K_f)^-$
\beq\label{eq: defnition of ag(u)}
\fs_{K_f/\bbQ,S^1}(u)=\sum_{g\in G_f}a_g(u)g\inv\ \ \ \mbox{where}\ \ \ a_g(u)=\frac{1}{f}\sum_{h\in G_f}jh(\xi_f/(1-\xi_f))\log_p(\delta_1^{(j)}(hgu))
\eeq
Let $D$ denote $D_\fP(K_f/\bbQ)$, identified with $\Gal(K_{f,\fP}/\bbQ_p)$, for one hence any $\fP\in S_p(K_f)$.
Recall that $\fP_1\in S_p(K_f)$ is the ideal defined by the embedding $j=j\tau_1$ which therefore
gives rise to an isomorphism (also denoted $j$) from $K_{f,\fP_1}$ to $\overline{j(K_f)}=\hat{K}_f$. This  in turn
induces an isomorphism from $D$ to
$\hat{D}:=\Gal(\hat{K}_f/\bbQ_p)$ sending $d\in D$ to $\hat{d}$ say, where $jd=\hat{d}j$. For each $\fP\in S_p(K_f)$ we choose $h_\fP\in G_f$ such that $h_\fP(\fP)=\fP_1$ so that
$h_\fP$ extends to an isomorphism from $K_{f,\fP}$ to $K_{f,\fP_1}$. Thus
$G=\bigcup_\fP h_\fP D$ and for any $d$ in $D$, $jh_\fP d=\hat{d}jh_\fP$ defines an isomorphism from $K_{f,\fP}$ to $\hat{K}_f$. It follows that if
$u\in U^1(K_{f,p})^-$ and $g\in G$, then
$\log_p(\delta_1^{(j)}(h_\fP dgu))=\log_p(j\iota_{\fP_1}(h_\fP dgu))=\log_p(jh_\fP d\iota_\fP(gu))=\hat{d}\log_p(jh_\fP \iota_\fP(gu))$. Consequently, we find
\begin{eqnarray}
a_g(u)&=&\sum_{\fP\in S_p(K_f)}\frac{1}{f}\sum_{d\in D}jh_\fP d(\xi_f/(1-\xi_f))\log_p(\delta_1^{(j)}(h_\fP dgu))\nonumber\\
  &=&\sum_{\fP\in S_p(K_f)}\frac{1}{f}{\rm Tr}_{\hat{K}_f/\bbQ_p}(jh_\fP (\xi_f/(1-\xi_f))\log_p(jh_\fP \iota_\fP(gu))\nonumber\\
  &=&\sum_{\fP\in S_p(K_f)}b(\hat{\xi}_{f,\fP}, v_{g,\fP})\label{eq: formula for ag(u)}
\end{eqnarray}
where, for each $\fP\in S_p(K_f)$, we have set $\hat{\xi}_{f,\fP}:=jh_\fP(\xi_f)$ and $v_{g,\fP}:=jh_\fP \iota_\fP(gu)$ and where
$b(\hat{\xi}_{f,\fP}, v_{g,\fP})$ is as defined in Lemma~\ref{lemma: Coleman's rec. law}.
The first statement of this Lemma therefore shows that $a_g(u)\in \bbZ_p$ for all $u\in U^1(K_{f,p})^-$ and $g\in G$, \ie\ that $IC(K_f/\bbQ,S^1,p)$ holds. Also,
the definition of the pairing
$[\cdot,\cdot]_{\fP,m+1}$ gives
\[
\iota_\fP(\xi_f^{f'})^{[1-\xi_f,\iota_\fP(gu)]_{\fP,m+1}}=(1-\xi_f,\iota_\fP(gu))_{K_\fP,p^{m+1}}
\]
Applying $jh_\fP$ to both sides and using the second statement of Lemma~\ref{lemma: Coleman's rec. law}, we get
\[
(\hat{\xi}_{f,\fP}^{f'})^{[1-\xi_f,\iota_\fP(gu)]_{\fP,m+1}}=(1-\hat{\xi}_{f,\fP}, v_{g,\fP})_{\hat{K}_f,p^{m+1}}=(\hat{\xi}_{f,\fP}^{f'})^{-b(\hat{\xi}_{f,\fP}, v_{g,\fP})}
\]
which implies that $b(\hat{\xi}_{f,\fP}, v_{g,\fP})\equiv -[1-\xi_f,\iota_\fP(gu)]_{\fP,m+1}\pmod{p^{m+1}}$. Summing this congruence over all $\fP\in S_p(K_f)$ and combining with
equations~(\ref{eq: formula for ag(u)}),~(\ref{eq: defnition of ag(u)}) and~(\ref{eq: H eta in case k= Q}), we obtain
$\fs_{K_f/\bbQ,S^1}(u)\equiv H_{K_f/\bbQ,m}(\eta_f,u)\pmod{p^{m+1}}$
for any $u\in U^1(K_f)^-$ giving $CC(K_f/\bbQ,S^1,p,m)$.\ePf
\section{The Case of $K$ Absolutely Abelian}\label{sec: Kabsab}
If $L/M$ is any Galois extension of number fields and $\phi$ any complex character of $\Gal(L/M)$, the \emph{$T$-truncated Artin $L$-function} $L_{L/M,T}(s,\phi)$
is defined for any finite set
$T$ of places of $M$ containing $S_\infty(M)$ but \emph{not necessarily} $\Sram(L/M)$. If  $\Gal(L/M)$ is abelian and $\phi$ is irreducible
(\ie\ $\phi\in\widehat{\Gal(L/M)}$) then, as noted in Remark~\ref{rem: artin L functions}, the definition agrees with the third member in~(\ref{eq:Lfns}).
In particular, there is no conflict with previous notation in the case $T\supset\Sram(L/M)$ and we always have
\beq\label{eq: extended defn of L fn}
L_{L/M, T}(s,\phi)=\prod_{\fq\in T\setminus S_\infty(M)\atop \fq\ndiv \ff_\phi}\left(1-N\fq^{-s}\hat{\phi}([\fq])\right)L(s,\hat{\phi})=
\prod_{\fq\nin T\cup \Sram(L^\phi/M)}(1-N\fq^{-s}\hat{\phi}(\fq))\inv
\eeq
where  $\hat{\phi}$  denotes the associated primitive ray-class
character modulo $\ff_\phi$,  $L^\phi=L^{\ker(\phi)}$  and the infinite product converges only for ${\rm Re}(s)>1$.

\begin{lemma}\label{lemma: induction of Artin L-functions}
Suppose $L/M$ and $T$ are as above, with $\Gal(L/M)$ abelian, and suppose $l$  is any intermediate field, $L\supset l\supset M$.
Then for any $\chi\in\widehat{\Gal(L/l)}$, we have an identity of meromorphic functions on $\bbC$:
\[
L_{L/l, T(l)}(s,\chi)=\prod_{\phi\in\widehat{\Gal(L/M)}\atop \phi|_{\Gal(L/l)}=\chi} L_{L/M, T}(s,\phi)
\]
\end{lemma}
\bPf\ This follows from the usual induction and `additivity' properties
for Artin $L$-functions (see~\cite[p.~15]{Tate}) and the fact (\eg\ by Frobenius Reciprocity)
that ${\rm Ind}_{\Gal(L/l)}^{\Gal(L/M)}\chi=\sum_{\phi\in\widehat{\Gal(L/M)}\atop \phi|_{\Gal(L/l)}=\chi}\phi$.\ePf
\begin{lemma}\label{lemma: dets on group-rings}
Let $B$ be a finite abelian group, $C$ any subgroup of $B$ and $x$ any element of $\bbC B$.
We write $x|\bbC B$ for the endomorphism of $\bbC B$, considered as a free $\bbC C$-module, determined by multiplication by $x$.
For any $\chi\in\hat{C}$, we have
\[
{\textstyle \chi(\det_{\bbC C}(x|\bbC B))}=\prod_{\phi\in\hat{B}\atop \phi|_{C}=\chi}\phi(x)
\]
(All characters extended linearly to homomorphisms from the complex group-rings to $\bbC$).
\end{lemma}
\bPf\ Choose any $\bbC C$-basis $\cB=\{y_1,\ldots,y_n\}$ for $\bbC
B$ (where $n=|B:C|$) and let $T=(t_{ij})_{i,j}\in M_n(\bbC C)$ be
the matrix of $x|\bbC B$ w.r.t.\ $\cB$. If $e_{\chi,C}$ denotes the
idempotent attached to $\chi$ in $\bbC C$, then  $x|\bbC B$ acts on
the submodule $e_{\chi,C}\bbC B$ and its matrix w.r.t.\ the
$\bbC$-basis $\{e_{\chi,C}y_1,\ldots,e_{\chi,C}y_n\}$ of the latter
is clearly $\chi(T):=(\chi(t_{ij}))_{i,j}\in M_n(\bbC)$. Hence $
\chi(\det_{\bbC C}(x|\bbC
B))=\chi(\det(T))=\det(\chi(T))=\det_\bbC(x|{e_{\chi,C}\bbC B}) $.
On the other hand, $e_{\chi,C}\bbC B$ also has a $\bbC$-basis
consisting of the $\bbC B$-idempotents $e_{\phi,B}$ for the
characters $\phi\in\hat{B}$ such that $\phi|_{C}=\chi$. (This follows easily from the fact that
$e_{\chi,C}$  is the sum of the corresponding $e_{\phi,B}$'s). The
result follows, since $xe_{\phi,B}=\phi(x)e_{\phi,B}$. \ePf
\noindent For the rest of this section, we fix $K/k$, $S$ and $p$ satisfying the standard hypotheses
with $K$  absolutely abelian. Thus $G=\Gal(K/k)$ is a
subgroup of the abelian group $\Gamma:=\Gal(K/\bbQ)$.
We define a set of places $S_\bbQ$ of $\bbQ$ by
\[
S_\bbQ=\{\infty\}\cup\{\mbox{$q$  prime  such that $S_q(k)\subset S$}\}
\]
Thus $p\in S_\bbQ$ and $S_\bbQ(k)$ is the maximal $\Gal(k/\bbQ)$-stable subset of $S$. The definition of ${\rm Bad(S)}$ in Subsection~\ref{subsec: case of absab preview} gives
\[
\Sram(K/\bbQ)={\rm Bad}(S)\cupdot(\Sram(K/\bbQ)\cap S_\bbQ)
\]
(disjoint union).
Let us
write $A$ for the subgroup $\prod_{q\in{\rm Bad}(S)}T_q(K/\bbQ)$ of
$\Gamma$ (trivial if ${\rm Bad}(S)=\varnothing$). If $F$ is any subfield of $K$, it follows that
\beq\label{eq: equiv condits for F in T to the A}
F\subset K^A\Leftrightarrow\mbox{all primes $q\in {\rm Bad}(S)$ are unramified in $F$}\Leftrightarrow\Sram(F/\bbQ)\subset S_\bbQ
\eeq
We denote by  $\cX_{\bbQ}(A)$  the set of irreducible $\bbQ$-valued characters of $A$. Each $\cA\in \cX_{\bbQ}(A)$ corresponds to a $\Gal(\barbbQ/\bbQ)$-conjugacy
class of characters $\alpha\in\hat{A}$ which `belong' to $\cA$. We set $\ker(\cA):=\ker(\alpha)$ for one (hence any) such $\alpha$  and
$K^\cA:=K^{\ker(\cA)}\supset K^A$. We define
\[
S_\cA=S_\bbQ \cup \Sram(K^\cA/\bbQ)\supset S^1(K^\cA/\bbQ)
\]
and write $\tilde{\nu}_{\cA}$ for the `averaged corestriction' map $|\ker(\cA)|\inv\nu_{K/K^\cA}$ which is a (non-unital)
homomorphism from $\bbC \Gal(K^A/\bbQ)$ to $\bbC \Gamma$. Finally, let $e_\cA$
denote the idempotent of $\bbQ A$ corresponding to $\cA$. With these notations, we define a meromorphic function
\displaymapdef{x_{K/k,S}}{\bbC}{\bbC \Gamma}{s}
{\displaystyle\sum_{\cA\in\cX_{\bbQ}(A)}e_\cA\tilde{\nu}_{\cA}(\Theta_{K^\cA/\bbQ,S_\cA}(s))}
\begin{prop}\label{prop: Theta K over k as a det}
With the above hypotheses and notations
\beq\label{eq: theta equals det x}
\Theta_{K/k,S_\bbQ(k)}(s)={\textstyle\det_{\bbC G}}(x_{K/k,S}(s)|\bbC \Gamma)
\eeq
(as $\bbC G$-valued meromorphic functions of $s\in \bbC$).
\end{prop}
\bPf\ By meromorphic continuation, it suffices to prove
$\chi({\rm L.H.S})=\chi({\rm R.H.S})$ in~(\ref{eq: theta equals det x}), for ${\rm Re}(s)>1$
and for all $\chi\in \hat{G}$. Equation~(\ref{eq:thetaStozetaandL}) and
Lemma~\ref{lemma: induction of Artin L-functions}
give
\[
\chi(\mbox{L.H.S of (\ref{eq: theta equals det x})})=L_{K/k,S_\bbQ(k)}(s,\chi\inv)=
\prod_{\phi\in\hat{\Gamma}\atop \phi|_G=\chi} L_{K/\bbQ, S_\bbQ}(s,\phi\inv)
\]
and evaluating $\chi(\mbox{R.H.S of \ref{eq: theta equals det x}})$ \emph{via} Lemma~\ref{lemma: dets on group-rings},
it suffices to show that $L_{K/\bbQ, S_\bbQ}(s,\phi\inv)=\phi(x_{K/k,S}(s))$ for any $\phi\in\hat{\Gamma}$. Suppose
$\alpha_\phi:=\phi|_A$ belongs to
$\cA_\phi\in\cX_\bbQ(A)$ so that $\ker(\alpha_\phi)=A\cap\ker(\phi)$ and
$K^{\cA_\phi}=K^AK^\phi$.
On the one hand, this means that  $\phi$ factors through $\Gal(K^{\cA_\phi}/\bbQ)$ and
$\phi(e_{\cA_\phi}\tilde{\nu}_{\cA_\phi}(y))=\phi(y)$ for all $y\in\bbC\Gal(K^{\cA_\phi}/\bbQ)$, while $\phi(e_\cA)=0$ for any  $\cA\neq\cA_\phi$. On the other, it implies
that $\Sram(K^{\cA_\phi}/\bbQ)=\Sram(K^A/\bbQ)\cup \Sram(K^\phi/\bbQ)$. Now, crucially for our argument, (\ref{eq: equiv condits for F in T to the A}) implies that
$\Sram(K^A/\bbQ)\subset S_\bbQ$ so $S_{\cA_\phi}=S_\bbQ\cup\Sram(K^\phi/\bbQ)$.
Putting this together,~(\ref{eq:thetaStozetaandL}),~(\ref{eq:Lfns})
and~(\ref{eq: extended defn of L fn}) give (for ${\rm Re}(s)>1$):
 \begin{eqnarray*}
\phi(x_{K/k,S}(s))
=\phi(e_{\cA_\phi}\tilde{\nu}_{\cA_\phi}(\Theta_{K^{\cA_\phi}/\bbQ,S_{\cA_\phi}}(s)))
&=&L_{K^{\cA_\phi}/\bbQ,S_{\cA_\phi}}(s,\phi\inv)\\
&=&\prod_{q\nin S_\bbQ\cup \Sram(K^{\phi}/\bbQ)}(1-q^{-s}\hat{\phi}\inv(q))\\
&=&L_{K/\bbQ, S_\bbQ}(s,\phi\inv)\hspace{10em}\Box\\
\end{eqnarray*}

\noindent Let us write $\cX^-_{\bbQ}(A)$ for the set  $\{\cA\in\cX_{\bbQ}(A)\,:\,c\nin\ker(\cA)\}$ ($=\cX_{\bbQ}(A)$ if $c\nin A$) and $x^-_{K/k,S}$ for the function
$e^-x_{K/k,S}:\bbC\rightarrow\bbC\Gamma^-$. If $\cA\in\cX_{\bbQ}(A)$ lies in $\cX^-_{\bbQ}(A)$ then $K^\cA$ is CM. Otherwise $e^-e_\cA=0$. Therefore
$
x^-_{K/k,S}(s)$ equals $\sum_{\cA\in\cX^-_{\bbQ}(A)}e_\cA\tilde{\nu}_{\cA}(
\Theta^-_{K^\cA/\bbQ,S_\cA}(s)
)
$ and is entire as a function of $s$.
Now take $s=1$, multiply  by $i/\pi$ and apply the involution $\_^\ast:\bbC \Gamma\rightarrow\bbC\Gamma$  (which fixes each $e_\cA$) to get
\beq\label{eq: one of two}
\left(\frac{i}{\pi}\right)x^-_{K/k,S}(1)^\ast=\sum_{\cA\in\cX^-_{\bbQ}(A)}e_\cA\tilde{\nu}_{\cA}(
a_{K^\cA/\bbQ,S_\cA}^{-,\ast})
\eeq
which lies in $\barbbQ\Gamma$ by~(\ref{eq:generating QmufGminus}). On the other hand, multiplying $\Theta_{K/k,S}(s)$ by $(i/\pi)^de^-=((i/\pi)e^-)^{|\Gamma:G|}$ in the previous Proposition
and letting $s\rightarrow 1$ implies that $a^-_{K/k,S_\bbQ(k)}$ is the $\bbC G$-determinant of $(i/\pi)x^-_{K/k,S}(1)$ acting on $\bbC \Gamma$. It follows easily from this that
\beq\label{eq: two of two}
a^{-,\ast}_{K/k,S_\bbQ(k)}={\textstyle\det_{\barbbQ G}}\left(\left(\frac{i}{\pi}\right)x^-_{K/k,S}(1)^\ast|\barbbQ\Gamma\right)
\eeq
For each $\cA\in\cX^-_{\bbQ}(A)$ the data $K^\cA/\bbQ$, $S_\cA$ and $p$ satisfy the standard hypotheses. In particular we have a well-defined $\bbZ_p \Gal(K^\cA/\bbQ)$-linear map
$\fs^{\rm id}_{K^\cA/\bbQ,S_\cA}$ from $U^1(K^\cA_p)^-$ to $\bbQ_p \Gal(K^\cA/\bbQ)^-$ (where `${\rm id}$' denotes the identity element of $\Gal(\barbbQ/\bbQ)$)
Both the norm map
$N_{K/K^\cA}:U^1(K_p)\rightarrow U^1(K^\cA_p)$ and the averaged corestriction
$\tilde{\nu}_\cA:\bbQ_p \Gal(K^\cA/\bbQ)\rightarrow\bbQ_p \Gamma$ take minus parts to minus parts. The
automorphism $\tau_i\in\Gal(\barbbQ/\bbQ)$ restricts to an element $\gamma_i:=\tau_i|_K$ of $\Gamma$ for $i=1,\ldots,d$
such that $\{\gamma_1\inv,\ldots,\gamma_d\inv\}$  is a set of coset representatives for
$G$ in $\Gamma$, hence also a basis for
$\cR\Gamma$ over $\cR G$, for any commutative ring $\cR$. We can now state:
\begin{thm}\label{thm: description of fs for K abs ab} With the above hypotheses and notations, suppose that $u_1,\ldots,u_d$ are any elements of $U^1(K_p)^-$.
Then
\[
\fs^{\tau_1,\ldots,\tau_d}_{K/k,S_\bbQ(k)}(u_1\wedge\ldots\wedge u_d)=\det(c_{i,l})_{i,l=1}^d
\]
where $c_{i,l}\in\bbQ_p G^-$ is the coefficient of $\gamma_i\inv$ when the element
${\displaystyle \sum_{\cA\in\cX^-_{\bbQ}(A)}}e_\cA\tilde{\nu}_{\cA}(\fs^{\rm id}_{K^\cA/\bbQ,S_\cA}(N_{K/K^\cA}u_l))$
of $\bbQ_p \Gamma^-$ is expressed in the $\bbQ_p G$-basis $\{\gamma_1\inv,\ldots,\gamma_d\inv\}$ of $\bbQ_p \Gamma$.
\end{thm}
\bPf\ Choose an embedding $j:\barbbQ\rightarrow\barbbQ_p$ inducing a prime ideal $\fP\in S_p(K)$, say, and write
$\lambda_p$ for the ($1\times 1$) regulator $R_{K/\bbQ,p}^{(j;{\rm id})}:U^1(K_p)\rightarrow\barbbQ_p\Gamma$.
If $u\in U^1(K_p)$ then, by definition,
\[
\lambda_p(u)
=\sum_{i=1}^d\sum_{g\in G}\log_p(j\circ\iota_\fP(\gamma_ig\inv u))(g\gamma_i\inv)
=\sum_{i=1}^d\left(\sum_{g\in G}\log_p(j\tau_i\circ\iota_{\fP_i}(g\inv u))g\right)\gamma_i\inv
\]
where $j\tau_i$ induces $\fP_i\in S_p(K)$. Now $\bigwedge_{\barbbQ_p G}^d\barbbQ_p \Gamma$ is $\barbbQ_p G$-free of rank one
on $\gamma_1\inv\wedge\ldots\wedge\gamma_d\inv$ and it follows easily from the last equation and the definition of $R^{(j;\tau_1,\ldots,\tau_d)}_{K/k,p}$ that
\beq\label{eq: wedgelambdapus}
\lambda_p(u_1)\wedge\ldots\wedge\lambda_p(u_d)
=R^{(j;\tau_1,\ldots,\tau_d)}_{K/k,p}(u_1\wedge\ldots\wedge u_d)\gamma_1\inv\wedge\ldots\wedge\gamma_d\inv\ \ \ \mbox{in $\bigwedge_{\barbbQ_p G}^d\barbbQ_p \Gamma$}
\eeq
On the other hand, $R_{K^\cA/\bbQ,p}^{ (j;{\rm id})}\circ N_{K/K^\cA}=\pi_{K/K^\cA}\circ\lambda_p$
for each $\cA\in\cX_\bbQ^-(A)$ so that
$\fs_{K^\cA/\bbQ,S_\cA}^{\rm id}(N_{K/K^\cA}u_l)=j(a^{-,\ast}_{K^\cA/\bbQ,S_\cA})\pi_{K/K^\cA}(\lambda_p(u_l))$ for each $\cA$ and $l$. It follows that
\[
\sum_{i=1}^d c_{i,l}\gamma_i\inv
=\left(\sum_{\cA\in\cX^-_{\bbQ}(A)}e_\cA\tilde{\nu}_{\cA}(j(a^{-,\ast}_{K^\cA/\bbQ,S_\cA}))\right)\lambda_p(u_l)
=j((i/\pi)x^-_{K/k,S}(1)^\ast)\lambda_p(u_l)
\]
by equation~(\ref{eq: one of two}). Using equation~(\ref{eq: two of two}), we deduce easily that
\begin{eqnarray*}
\det(c_{i,l})_{i,l}\gamma_1\inv\wedge\ldots\wedge\gamma_d\inv
&=&(j((i/\pi)x^-_{K/k,S}(1)^\ast)\lambda_p(u_1))\wedge\ldots\wedge(j((i/\pi)x^-_{K/k,S}(1)^\ast)\lambda_p(u_l))\\
&=&j(a^{-,\ast}_{K/k,S_\bbQ(k)})\lambda_p(u_1)\wedge\ldots\wedge\lambda_p(u_l)\\
\end{eqnarray*}
and combining this with equation~(\ref{eq: wedgelambdapus}), the result follows from the definition of $\fs^{\tau_1,\ldots,\tau_d}_{K/k,S_\bbQ(k)}$.\ePf
\noindent {\sc Proof of Theorem~\ref{thm: IC for K absab} under Hypothesis~\ref{hyp: badS etc}}
By Prop.~\ref{prop: IC under increasing S} it suffices to prove $IC(K/k,S_\bbQ(k),p)$, \ie\ that
$\fs_{K/k,S_\bbQ(k)}(u_1\wedge\ldots\wedge u_d)$ lies in $ \bbZ_p G$ for all $u_i,\ldots,u_d\in U^1(K_p)^-$ and this will clearly follow from
Theorem~\ref{thm: description of fs for K abs ab} provided we can show
\beq\label{eq: suff for IC if K abs ab}
e_\cA\tilde{\nu}_{\cA}(\fs_{K^\cA/\bbQ,S_\cA}(N_{K/K^\cA}u_l))\in \bbZ_p \Gamma\ \ \ \ \forall\,l,\ \forall\,\cA\in\cX^-_{\bbQ}(A)
\eeq
But
Theorem~\ref{thm: CC and IC for K over Q}~\ref{part1: CC and IC for K over Q} implies that $\fs_{K^\cA/\bbQ,S_\cA}(N_{K/K^\cA}u_l)$ lies in $\bbZ_p \Gal(K^\cA/\bbQ)$. Furthermore,
if $q\in{\rm Bad}(S)$ then $|T_q(K/\bbQ)|=e_q(k/\bbQ)$ and Hypothesis~\ref{hyp: badS etc} implies that this is prime to $p$ for all such $q$, hence that $p\ndiv|A|$. It follows that
$p\ndiv[K:K^\cA]$ for every $\cA$ so that $e_\cA\in\Zbp A$  and $\tilde{\nu}_{\cA}(\fs_{K^\cA/\bbQ,S_\cA}(N_{K/K^\cA}u_l))\in\bbZ_p\Gamma$,
establishing~(\ref{eq: suff for IC if K abs ab}).\ePf
\noindent Turning to the Congruence Conjecture, we suppose from now on that \emph{$K$ contains $\mu_{\pnpo}$ for some $n\geq 0$}. Since $\Sram(\bbQ(\mu_{\pnpo})/\bbQ)=\{p\}\subset S_\bbQ$,
it follows from~(\ref{eq: equiv condits for F in T to the A}) that $K^A$ contains $\bbQ(\mu_{\pnpo})$ so is  CM and  $\cX^-_{\bbQ}(A)=\cX_{\bbQ}(A)$.
We write $\bar{\Gamma}$ for $\Gal(K^+/\bbQ)$.
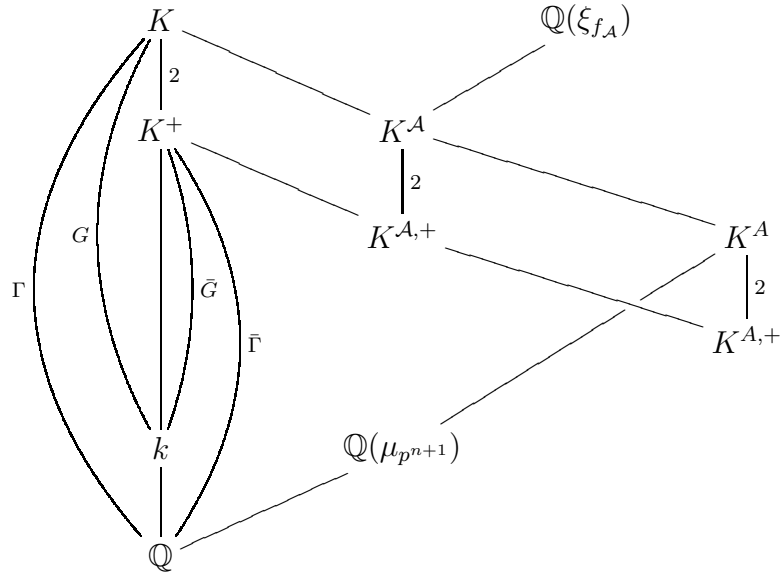
\[
\xymatrix{
K\ar@{-}[d]^2\ar@{-}[drr]
\ar@/_4pc/@{-}[ddddd]_{\Gamma}\ar@/_2pc/@{-}[dddd]_{G}&&&\bbQ(\xi_{f_\cA})\\
K^+\ar@{-}[ddd]\ar@{-}[drr]
\ar@/^2.5pc/@{-}[dddd]^{\bar{\Gamma}}
\ar@/^1pc/@{-}[ddd]^{\bar{G}}&&K^\cA\ar@{-}[d]^2\ar@{-}[drr]\ar@{-}[ru]\\
&&K^{\cA,+}
\ar@{-}[drr]&&K^A\ar@{-}[d]^2\ar@{-}|!{[ll];[d]}\hole[ddll]\\
&&&&K^{A,+}\\
k\ar@{-}[d]&&{\bbQ(\mu_\pnpo)}\ar@{-}[dll]\\
\bbQ
}
\]
Now $RSC(K^{\cA,+}/\bbQ,S_\cA;\bbQ)$ holds for each $\cA\in\cX_\bbQ(A)$. Indeed, let us write $f_\cA$ for the conductor of $K^\cA$
so that $\pnpo|f_\cA$ and $S^1(K^\cA/\bbQ)=S^1(\bbQ(\xi_{f_\cA})/\bbQ)=\{\infty\}\cup S_{f_\cA}(\bbQ)$.
Then, the determination of $\eta_{\bbQ(\xi_{f_\cA})^+/\bbQ, S^1(\bbQ(\xi_{f_\cA})/\bbQ)}$ in the proof
of Theorem~\ref{thm: CC and IC for K over Q}, together with~Props.~\ref{prop: RSC(Q) under descending K}, and~\ref{prop: RSC under increasing S}
imply that the solution $\eta_{K^{\cA,+}/\bbQ,S_\cA}$ of $RSC(K^{\cA,+}/\bbQ,S_\cA;\bbQ)$ (with
$\tau_1={\rm id}$) is
\begin{eqnarray*}
\eta_\cA&:=&\prod_{q\in S_\cA\setminus S^1(K^\cA/\bbQ)}(1-\sigma_{q,K^{\cA,+}/\bbQ}\inv)
N_{\bbQ(\xi_{f_\cA})^+/K^{\cA,+}}(\eta_{\bbQ(\xi_{f_\cA})^+/\bbQ, S^1(\bbQ(\xi_{f_\cA})/\bbQ)})\\
&=&\prod_{q\in S_\cA\setminus S^1(K^\cA/\bbQ)}(1-\sigma_{q,K^{\cA,+}/\bbQ}\inv)
N_{\bbQ(\xi_{f_\cA})^+/K^{\cA,+}}({\textstyle\frac{1}{2}}\otimes(1-\xi_{f_\cA})^{1+c})
\end{eqnarray*}
In fact, $(1-\xi_{f_\cA})$ lies in $ U_{\{\infty\}}(\bbQ(\xi_{f_\cA}))$ unless $f_\cA$ is a power of a prime,
necessarily $p$, in which case it lies in $ U_{\{\infty,p\}}(\bbQ(\xi_{f_\cA}))$. Thus for all $\cA\in\cX_\bbQ(A)$, the element
$\eta_\cA$
lies in
${\textstyle\frac{1}{2}}\overline{ U_{\{\infty, p\}}(K^{\cA,+})}\subset \bbQ U_{\{\infty, p\}}(K^{\cA,+})$.
Let us write $i_\cA$ for the natural injection $\bbQ U_{\{\infty, p\}}(K^{\cA,+})\rightarrow\bbQ U_{\{\infty, p\}}(K^+)$ and
$\tilde{i}_{\cA}$ for $|\ker(\cA)|\inv i_\cA$. We define
\[
\alpha_{S_\bbQ(k)}:=\sum_{\cA\in\cX_{\bbQ}(A)}e_\cA\tilde{i}_\cA(\eta_\cA)
\ \ \ \mbox{and}\ \ \ \eta_{S_\bbQ(k)}:=\gamma_1\inv\alpha_{S_\bbQ(k)}\wedge\ldots\wedge\gamma_d\inv\alpha_{S_\bbQ(k)}
\]
from which it is clear that $\alpha_{S_\bbQ(k)}$ lies in $|A|^{-2}{\textstyle\frac{1}{2}}\overline{ U_{\{\infty, p\}}(K^+)}$ and
$\eta_{S_\bbQ(k)}$ in $(|A|^{-2}{\textstyle\frac{1}{2}})^d\overline{\bigwedge_{\bbZ \bar{G}}^d U_{\{\infty, p\}}(K^+)}$.

\begin{prop}\label{prop: explicit solution of WRSC for Kabsab} With the above hypotheses,
$RSC(K^+/k,S_\bbQ(k);\bbQ)$ holds with solution $\eta_{K^+/k,S_\bbQ(k)}=\eta_{S_\bbQ(k)}$.
\end{prop}
We defer the proof.
The final ingredient in the Proof of Theorem~\ref{thm: CC for Kabsab} is
\begin{lemma}\label{lemma: calculating H for absab}
Suppose $\alpha$ is an element of $\overline{ U_{\{\infty, p\}}(K^+)}$ so that
$\gamma_1\inv\alpha\wedge\ldots\wedge\gamma_d\inv\alpha$ lies in the subset $\overline{{\textstyle\bigwedge^d_{\bbZ \bar{G}}}U_{S_\bbQ(k)}(K^+)}$ of
${\textstyle\bigwedge^d_{\bbQ \bar{G}}}\bbQ U_{S_\bbQ(k)}(K^+)$. Then, for
any $u_1,\ldots,u_d\in U^1(K_p)^-$ we have
\[
\kappa_n(\tau_1\ldots\tau_d)H_{K/k,n}(\gamma_1\inv\alpha\wedge\ldots\wedge\gamma_d\inv\alpha,u_1\wedge\ldots\wedge u_d)=
\det(d_{i,l})_{i,l=1}^d
\]
where $d_{i,l}\in(\bbZ/\pnpo\bbZ)G^-$ is the coefficient of $\gamma_i\inv$ when the element
$H_{K/\bbQ,n}(\alpha,u_l)$
of $(\bbZ/\pnpo\bbZ)\Gamma^-$ is expressed in the $(\bbZ/\pnpo\bbZ) G$-basis $\{\gamma_1\inv,\ldots,\gamma_d\inv\}$ of $(\bbZ/\pnpo\bbZ)\Gamma$.
\end{lemma}
\bPf\ If $\alpha=1\otimes\eps$ for some $\eps\in  U_{\{\infty, p\}}(K^+)$ then $\gamma_i\inv\alpha=1\otimes \gamma_i\inv\eps$ with
$\gamma_i\inv\eps\in  U_{\{\infty, p\}}(K^+)\subset
U_{S}(K^+)$ for all $i$.
Equations~(\ref{eq:simple formula for H_n})
and~(\ref{eq:Gaction2}) applied to $K/\bbQ$ give
\[
H_{K/\bbQ,n}(\alpha,u_l)=
\sum_{i=1}^d\left(\sum_{g\in G}[\eps,\gamma_i g u_l]_{K,n}g\inv\right)\gamma_i\inv
=\sum_{i=1}^d\left(\kappa_n(\tau_i)\sum_{g\in G}[\gamma_i\inv\eps, g u_l]_{K,n}g\inv\right)\gamma_i\inv
\]
since $\gamma_i=\tau_i|_K$. Thus $d_{i,l}=\kappa_n(\tau_i)\sum_{g\in
G}[\gamma_i\inv\eps, g u_l]_{K,n}g\inv$. Now use~(\ref{eq:simple
formula for H_n}) for  $K/k$.\ePf
\noindent \textsc{Proof of
Theorem~\ref{thm: CC for Kabsab}} By  Prop.~\ref{prop: CC under increasing S} it suffices to establish $CC(K/k,S_\bbQ(k),p,n)$ under Hypothesis~\ref{hyp: badS etc}.
But the latter has already been shown  to imply  $IC(K/k,S_\bbQ(k),p)$ and that $p\ndiv|A|$. Therefore
$\eta_{S_\bbQ(k)}$ lies in $\Zbp\overline{\bigwedge_{\bbZ \bar{G}}^d U_{\{\infty, p\}}(K^+)}\subset\Zbp\Lambda_{0,S_\bbQ(k)}$
and so is the solution of $RSC(K^+/k,S_\bbQ(k);\Zbp)$ by Proposition~\ref{prop: explicit solution of WRSC for Kabsab}. It remains to prove the
congruence~(\ref{eq:the congruence})
holds with $\eta_{K^+/k,S}=\eta_{S_\bbQ(k)}$ and  $\theta=u_1\wedge\ldots\wedge u_d$ with
$u_i\in U^1(K_p)^-\ \forall\,i$. (Such $\theta$ generate $\bigwedge_{\bbZ_p G}^dU^1(K_p)^-$.) For each $\cA\in\cX_{\bbQ}(\cA)$ we may write
$2\eta_\cA$ as $1\otimes\eps_\cA$ where
$\eps_\cA$ lies in $U_{\{\infty,p\}}(K^{\cA,+})$.
From equation~(\ref{eq: H symbols and the norm}) with $F=K^\cA$
and~(\ref{eq:simple formula for H_n}) (with $d=1$!) it follows easily that
\[
H_{K/\bbQ,n}(i_\cA(2\eta_\cA),u_l)=\nu_{K/K^\cA}(H_{K^\cA/\bbQ,n}(2\eta_\cA, N_{K/K^\cA}u_l))
\]
Therefore, using the $\bbZ \Gamma$-linearity of $H_{K/\bbQ,n}(\cdot,\cdot)$ in the first variable and the fact $|A|e_{\cA}\in\bbZ \Gamma$, we have, for each $l$:
\begin{eqnarray*}
H_{K/\bbQ,n}(2|A|^2 \alpha_{S_\bbQ(k)},u_l)
&=&\sum_{\cA\in\cX_{\bbQ}(A)}(|A||A:\ker(\cA)|e_\cA) H_{K/\bbQ,n}(i_\cA(2\eta_\cA), u_l)\\
  &=&\sum_{\cA\in\cX_{\bbQ}(A)}(|A||A:\ker(\cA)|e_\cA) \nu_{K/K^\cA}(H_{K^\cA/\bbQ,n}(2\eta_\cA, N_{K/K^\cA}u_l))\\
  &\equiv&\sum_{\cA\in\cX_{\bbQ}(A)}2(|A||A:\ker(\cA)|e_\cA) \nu_{K/K^\cA}(\fs^{\rm id}_{K^{\cA}/\bbQ,S_\cA}(N_{K/K^\cA}u_l))\\
  &\equiv&\sum_{\cA\in\cX_{\bbQ}(A)}2|A|^2 e_\cA \tilde{\nu}_\cA(\fs^{\rm id}_{K^{\cA}/\bbQ,S_\cA}(N_{K/K^\cA}u_l))\\
  &\equiv&\sum_{i=1}^d 2|A|^2c_{i,l}\gamma_i\inv\pmod{p^{n+1}}
\end{eqnarray*}
where $c_{i,l}$ is precisely as defined in Theorem~\ref{thm: description of fs for K abs ab}. Note that the first congruence above comes from
Theorem~\ref{thm: CC and IC for K over Q} which also shows that the last three expressions lie in $\bbZ_p\Gamma$.
It follows from Proposition~\ref{prop: explicit solution of WRSC for Kabsab} and
Lemma~\ref{lemma: calculating H for absab} that
\begin{eqnarray*}
\lefteqn{\overline{(2|A|^2)^d}\kappa_n(\tau_1\ldots\tau_d)H_{K/k,n}(\eta_{S_\bbQ(k)},u_1\wedge\ldots\wedge u_d)}\hspace*{10ex}&&\\
&=&\kappa_n(\tau_1\ldots\tau_d)H_{K/k,n}(\gamma_1\inv(2|A|^2 \alpha_{S_\bbQ(k)})\wedge\ldots\wedge\gamma_d\inv(2|A|^2 \alpha_{S_\bbQ(k)}),u_1\wedge\ldots\wedge u_d)\\
&=&\det(\overline{2|A|^2c_{i,l}})_{i,l}
\end{eqnarray*}
in $(\ZmodZ{\pnpo})G$. Since $p\ndiv 2|A|$, we may cancel the factor $\overline{(2|A|^2)^d}\in (\ZmodZ{\pnpo})^\times$ on both sides above and combining with
Theorem~\ref{thm: description of fs for K abs ab} we obtain
\[
\overline{\fs_{K/k,S_\bbQ(k)}^{\tau_1,\ldots,\tau_d}(u_1\wedge\ldots\wedge u_d)}=
\det(\overline{c_{i,l}})_{i,l}=
\kappa_n(\tau_1\ldots\tau_d)H_{K/k,n}(\eta_{S_\bbQ(k)},u_1\wedge\ldots\wedge u_d)
\]
as required. \ePf \rem\label{rem: Burns and Cooper} Burns has proven
Conjecture~$B'$ of~\cite{Rubin} whenever $K$ is absolutely abelian
(see~\cite[Theorem A]{Burns}). It follows from
Remark~\ref{rubin_stark_remark 2} that
$\eta_{S_\bbQ(k)}=\eta_{K^+/k, S_\bbQ(k)}$ must also lie in
$\half\Lambda_{0,S_{\bbQ}(k)}(K^+/k)$, although this is not obvious
from our expression for $\eta_{S_\bbQ(k)}$ and Burns' results do not
appear to provide an explicit expression. On the other hand Cooper
obtains essentially our expression
$\gamma_1\inv\alpha_{S_\bbQ(k)}\wedge\ldots\wedge\gamma_d\inv\alpha_{S_\bbQ(k)}$
in~\cite{Cooper}. (Indeed, we adapt his methods in the proof of
Proposition~\ref{prop: explicit solution of WRSC for Kabsab} below.)
By manipulating it cleverly and using the norm relations for
cyclotomic numbers, he shows explicitly that if $A$ is cyclic and
of odd order, then $\eta_{S_\bbQ(k)}$ lies in $2^{-d}\overline{\bigwedge^d_{\bbZ \bar{G}}U_{\{\infty,p\}}(K^+)}$. (This follows from~\cite[Theorem 5.2.2]{Cooper}).\vertsp\\
\textsc{Proof of Proposition~\ref{prop: explicit solution of WRSC for Kabsab}}
The arguments are mostly familiar by now: Applying $e^+$ to~(\ref{eq: theta equals det x}) one deduces that $\Theta_{K^+/k,S_\bbQ(k)}(s)$ is the determinant of
$\sum_\cA e_\cA|\ker(\cA)|\inv\nu_{K^+/K^{\cA,+}}(\Theta_{K^{\cA,+}/\bbQ, S_\cA}(s))$ acting on $\bbC\bar{\Gamma}$, where $\bar{\Gamma}:=\Gal(K^+/\bbQ)$ and we are identifying
$A$ with $\Gal(K^+/K^{A,+})$ by restriction. Now $\Theta_{K^{\cA,+}/\bbQ, S_\cA}(0)=0$ and $\Theta^{(1)}_{K^{\cA,+}/\bbQ, S_\cA}(0)$ has real coefficients for each $\cA$, so
\beq\label{eq: ThetaKpk as det}
\Theta_{K^+/k,S_\bbQ(k)}^{(r)}(0)={\textstyle\det_{\bbR\bar{G}}}\left(\sum_\cA e_\cA|\ker(\cA)|\inv\nu_{K^+/K^{\cA,+}}(\Theta_{K^{\cA,+}/\bbQ, S_\cA}^{(1)}(0))|\bbR\bar{\Gamma}\right)=
\det(e_{i,l})_{i,l=1}^d
\eeq
say, where $(e_{i,l})_{i,l=1}^d$ is the matrix of multiplication by $\sum_\cA\ldots$ on $\bbR\bar{\Gamma}$ w.r.t.\ the $\bbR\bar{G}$-basis ${\bar{\gamma}_l\inv:i=1,\ldots,d}$
(where $\bar{\gamma}_l:=\gamma_l|_{K^+}=\tau_l|_{K^+}$). Fix $l$ and take $\tau_1={\rm id}$. Using~(\ref{eq:rubin_stark}) for each extension $K^{\cA,+}/\bbQ$ and the relation
$\nu_{K^+/K^{\cA,+}}\circ\lambda_{K^{\cA,+}/\bbQ,1}=\lambda_{K^{+}/\bbQ,1}\circ i_{K^+/K^{\cA,+}}$, we get
\begin{eqnarray*}
\sum_\cA e_\cA|\ker(\cA)|\inv\nu_{K^+/K^{\cA,+}}(\Theta_{K^{\cA,+}/\bbQ, S_\cA}^{(1)}(0))\bar{\gamma}_l\inv&=&
\sum_\cA e_\cA|\ker(\cA)|\inv\nu_{K^+/K^{\cA,+}}(\lambda_{K^{\cA,+}/\bbQ,1}(\eta_\cA))\bar{\gamma}_l\inv\\
&=&
\lambda_{K^{+}/\bbQ,1}(\alpha_{S_\bbQ(k)})\bar{\gamma}_l\inv
\end{eqnarray*}
But for any element $\alpha=a\otimes\eps$ of $\bbQ U_S(K^+)$ (with $a\in\bbQ$)  we have
\[
\lambda_{K^{+}/\bbQ,1}(\alpha)\bar{\gamma}_l\inv=
\lambda_{K^{+}/\bbQ,1}(\bar{\gamma}_l\inv\alpha)=
a\sum_{i=1}^d\left(\sum_{\bar{g}\in\bar{G}}\log|\bar{\gamma}_i\bar{g}\bar{\gamma}_l\inv\eps|\bar{g}\inv\right)\bar{\gamma}_i\inv=
\sum_{i=1}^d\lambda_{K^{+}/k,i}(\bar{\gamma}_l\inv\alpha)\bar{\gamma}_i\inv
\]
and combining with the previous equation, we find $e_{i,l}=\lambda_{K^{+}/k,i}(\bar{\gamma}_l\inv\alpha_{S_\bbQ(k)})$. Substituting this
in~(\ref{eq: ThetaKpk as det}), it follows that $\eta_{S_\bbQ(k)}$ satisfies condition~(\ref{eq:rubin_stark})
for $K/k^+$ and $S_\bbQ(k)$.

To show that $\eta_{S_\bbQ(k)}$ satisfies the eigenspace condition w.r.t.\ $(S_\bbQ(k),d,\bar{G})$, one could adapt the argument of~\cite{Cooper}
(based on~\cite[Prop.~3.1.2]{Pop}) using condition~\ref{part4: equivalent `eigenspace' conditions} of Prop.~\ref{prop: equivalent `eigenspace' conditions}.
We sketch a more `algebraic' argument based on the equivalent condition~\ref{part3: equivalent `eigenspace' conditions}: Suppose $\fq\in S_\bbQ(k)\setminus S_\infty(k)$ lies above
$q\in S_\bbQ\setminus\{\infty\}$, write $D$ for $D_q(K/\bbQ)$ and $\fD$ for $D_\fq(K/k)=D\cap G$.
Let $\rho_1,\ldots,\rho_t$ be a set of representatives
for $D$ mod $\fD$, hence for $DG$ mod $G$, and let $\sigma_1,\ldots,\sigma_m$ be a set of representatives
for $\Gamma$ mod $DG$. Then $d=mt$ and both $\{\sigma_a\rho_b\}_{a,b}$ and $\{\gamma_i\inv\}_i$ are sets of representatives for $\Gamma$ mod $G$. Writing also
$\eta$ and $\alpha$ for $\eta_{S_\bbQ(k)}$ and $\alpha_{S_\bbQ(k)}$ respectively, it follows that
$\eta=\pm g\bigwedge_{a=1}^m\bigwedge_{b=1}^t\sigma_a\rho_b\alpha$ for some $g\in G$.
(The unordered `wedge product' (over $\bbQ\bar{G}$) on the RHS is defined only up to sign.)
Since $N_{D_\fq(K^+/k)}\eta$ equals $\half N_\fD\eta$ or $N_\fD\eta$, condition~\ref{part3: equivalent `eigenspace' conditions}
for $m=\eta$ and $S=S_\bbQ(k)$ will follow if we can show that $N_\fD\eta$ is fixed by $G$ (hence by $\bar{G}$) and is zero if $|S_\bbQ(k)|>d+1$. But
\beq\label{eq: formula for NfDeta}
N_\fD\eta=\pm |\fD|^{1-d}g\bigwedge_{a=1}^m\bigwedge_{b=1}^t\sigma_aN_\fD\rho_b\alpha=
\pm |\fD|^{1-d}g\bigwedge_{a=1}^m\left(\sigma_a N_D\alpha\wedge\sigma_a N_\fD\rho_2\alpha\wedge\ldots\wedge\sigma_a N_\fD\rho_t\alpha\right)
\eeq
(the second equality since $\sigma_a N_D\alpha=\sum_{b=1}^t\sigma_aN_\fD\rho_b\alpha$ for each $a$). If $|S_\bbQ|>2$ then $|S_\cA|>2$ for each $\cA\in\cX_\bbQ(A)$ so
the eigenspace condition on $\eta_\cA$ as a solution of $RSC(K^{\cA,+}/\bbQ,S_\cA;\bbQ)$
implies that it is annihilated by $N_{D_q(K^{\cA,+}/\bbQ)}$, hence by $N_D$. It follows that $N_D\alpha=0$ hence
$N_\fD\eta=0$ by~(\ref{eq: formula for NfDeta}). Otherwise, $|S_\bbQ|=2$, $S_\bbQ=\{\infty,q\}$ (so $q=p$) and $|S_\bbQ(k)|$ is precisely $d+m$. In this case, the
eigenspace condition on $\eta_\cA$ still shows that  $N_{D_q(K^{\cA,+}/\bbQ)}\eta_\cA$ is fixed by $\Gal(K^{\cA,+}/\bbQ)$ for all $\cA$
and it follows as above that $N_D\alpha$ is fixed by $G$. So equation~(\ref{eq: formula for NfDeta}) implies that $N_\fD\eta$ is fixed by $G$ and, if $m>1$, that it is zero, since
$\sigma_1 N_D\alpha=\sigma_2 N_D\alpha$.\ePf

\end{document}